\documentclass[12pt,a4paper]{amsart}
\usepackage{graphicx}
\usepackage{epstopdf}
\usepackage{amssymb}
\usepackage{amsthm}
\usepackage{amsmath}
\usepackage{pstcol, pst-node}
\usepackage{mathrsfs}
\usepackage{amscd}
\usepackage{xypic}
\usepackage{comment}

\includecomment{comment}
\usepackage[hidelinks]{hyperref}
\usepackage[all]{xy}
\newtheorem{thm}{Theorem}

\newtheorem{cor}[thm]{Corollary}
\newtheorem{prop}[thm]{Proposition}
\newtheorem{lem}[thm]{Lemma}
\newtheorem{rem}[thm]{Remark}

\theoremstyle{definition}
\newtheorem{defn}[thm]{Definition}

\newtheorem{prop-def}[thm]{Proposition-Definition}

\newtheorem{example}[thm]{Example}
\newtheorem{problem}[thm]{Problem}

\newcommand{\Spec}{\mathop{\mathrm{Spec}}\nolimits}

\theoremstyle{remark}

\begin{document}
\title{Deformation of singular curves on surfaces}

\author{Takeo Nishinou}
\date{}
\thanks{email : nishinou@rikkyo.ac.jp}
\address{Department of Mathematics, Rikkyo University,
	3-34-1, Nishi-Ikebukuro, Toshima, Tokyo, Japan} 
\subjclass[2000]{}
\maketitle
\begin{abstract}
In this paper, we consider deformations of singular complex curves on complex surfaces.   
Despite the fundamental nature of the problem, little seems to be known for curves on general surfaces.
Let $C\subset S$ be a complete integral curve on a smooth surface.
Let $\tilde C$ be a partial normalization of $C$, and 
 $\varphi\colon \tilde C\to S$ be the induced map.
In this paper, we consider deformations of $\varphi$.
The problem of the existence of deformations will be reduced to solving a certain explicit system of polynomial equations.
This system is universal in the sense that it is determined solely by simple local data of the singularity of $C$, 
 and does not depend on the global geometry of $C$ or $S$.
Under a relatively mild assumption on the properties of these equations, we will show that the map $\varphi$ has virtually
 optimal deformation property.
\end{abstract}

\section{Introduction}
The study of algebraic surfaces and curves on them has a long history.
With the exception of curves, it is one of the subjects which has been
 studied the longest and most frequently in algebraic geometry.
It has been the first and most deeply investigated in various developments in algebraic geometry,
 including birational geometry, various moduli theories,
 and more recent advances such as gauge theory and Gromov-Witten type invariants. 

However, knowledge of how curves on algebraic surfaces actually behave is surprisingly limited.
The most famous result known in this direction would be the so-called Severi's problem
 asking the irreducibility of the moduli space of nodal plane curves of the given genus and degree,
 solved affirmatively by Harris \cite{Ha}.
There are studies in this direction for other surfaces of non-positive Kodaira dimension, see \cite[Theorem B]{DS}.
See also \cite{GLS, GLS2} for extensive study from multiple points of view.
In spite of these studies, very few have been known about
 the behavior of curves on surfaces of positive Kodaira dimension, especially 
 those on surfaces of general type (see, for example, \cite{CC,CS} for results in this direction).
 
The most naive questions to be asked in the study in this direction would be the following.
\begin{problem}\label{problem:1}(\cite[Problem A]{DS})
Let $D$ be an integral curve on a smooth complete algebraic surface $X$.
Is it possible to deform $D$ in $X$ into a nodal curve while preserving its geometric genus?
\end{problem}
A weaker version of it asks the following.
\begin{problem}\label{problem:2}
Given $D$ and $X$ as above, let $C$ be the normalization of $D$ and $\varphi\colon C\to D$ be
 the natural map.
Is it possible to deform $\varphi$  into an immersion?
\end{problem}
These problems are almost completely open for surfaces of positive Kodaira dimension.
One of the main reasons for this would be the difficulty of applying powerful techniques of
 moduli theory to these problems due to the large potential obstruction.
In this paper, we attempt to fill this gap and obtain information about the behavior of curves 
 in a way independent of the nature of ambient surfaces.
Roughly speaking, our result claims that under a certain mild condition, 
 the deformation property of singular curves is 
 almost as optimal as possible on any surface.
In particular, assuming that condition, Problem \ref{problem:2} also has an almost optimal answer.\\

We study the so-called equigeneric (equivalently, equinormalisable) deformations of curves on surfaces
 from parametric point of view (see \cite{DS,GLS}).
In fact, we will deal with more general situations, but in the introduction, we restrict ourselves to 
 equigeneric deformations.
Specifically, given a complete integral curve $\overline C$ on a smooth algebraic surface $X$,
 we study the deformation of the map $\varphi\colon  C\to X$, 
 where $C$ is the normalization of $\overline C$ and $\varphi$ is the naturally induced map.
We assume that the map $\varphi$ satisfies the semiregularity condition 
 of Definition \ref{def:sr}, which is a natural generalization of the classical semiregularity 
 for embedded curves \cite{S1, S2, B,KS} and will be satisfied if the class of $\overline C$ is sufficiently ample (see the paragraph
 after Definition \ref{def:sr}).
As in these classical case, this is a natural assumption and there seems to be little hope to 
 control the deformation theory without it.

Let $p\in C$ be a point where the map $\varphi$ is not regular.
Then, taking a suitable analytic coordinate $s$ on 
 a neighborhood of $p$ in $C$ and coordinates $z, w$ 
 on a neighborhood of $\varphi(p)$ in $X$, the map $\varphi$ can be represented as
\[
(\varphi^*z, \varphi^*w) = (s^a, s^b + s^{b+1}g_0(s)),
\]
 where $a, b$ are integers satisfying $a<b$, and $g_0(s)$ is a convergent series 
 (see \cite[Chapter I, Corollary 3.8]{GLS}).
For notational ease, we will write $(\varphi^*z, \varphi^*w)$ and its analogues
 simply by $(z, w)$ from now on.
Up to a reparameterization of $C$, a $k$-th order deformation of $\varphi$ 
 on a neighborhood of $p$ is written in the form
\[
(z, w) = (s^a + \sum_{i=0}^{a-2}c_{a-i}s^i , \; s^b + s^{b+1}g_0(s)+\sum_{j=1}^kt^jg_{j}(s)),
\]
 where $k$ is a positive integer, $c_{a-i}\in t\Bbb C[t]/t^{k+1}$, $g_j(s)$ is a convergent series, 
 and $t$ is a generator of $\Bbb C[t]/t^{k+1}$.
It turns out that taking $c_{i}\in t^{i}\Bbb C[t]/t^{k+1}$ will be convenient and we assume this.
Given a global $k$-th order deformation $\varphi_k$ of $\varphi$, 
 the obstruction to deforming $\varphi_k$ one step further can be calculated by
 a \v{C}ech cocycle obtained as the difference of local $k+1$-th order deformations of $\varphi_k$, 
 which takes values in the normal sheaf of $\varphi$.
However, usually calculating this is very difficult and there is little hope to achieve it in general.

To overcome this difficulty, we choose a special type of deformations of $\varphi$ around each 
 singular point.
Namely, let $S$ be a local parameter on a punctured neighborhood of $p$
 on $C$ defined over $\Bbb C[[t]]$
 which satisfies
\[
S^a = s^a + \sum_{i=0}^{a-2}c_{a-i}s^i,
\]
 where we regard $c_i$ as an element of $t^i\Bbb C[[t]]$ 
 by taking the coefficients of $t^l$, $l>k$ to be zero.
Explicitly, we can choose $S$ by solving this equation, so that 
\[
S = s(1 + \sum_{i=1}^{\infty}\prod_{j=0}^{i-1}(\frac{1}{a}-j)\frac{1}{i!}(\sum_{l=2}^{a}\frac{c_l}{s^l})^i).
\]
Then, consider the deformation of $\varphi$ on the punctured neighborhood $\mathring U_p$ of $p$
 given by
\[
(z, w) = (S^a, S^b+S^{b+1}g_0(S)).
\]
Note that this reduces to the original parameterization $(z, w) = (s^a, s^b+s^{b+1}g_0(s))$
 of $\varphi$ over $\Bbb C[t]/t$.
In particular, though a priori the parameter $S$ is defined only on the punctured neighborhood 
 $\mathring U_p$, 
 it extends to the whole neighborhood $U_p$ of $p$ over $\Bbb C[t]/t$.
An important point is that this is true even up to some positive order of $t$.
Namely, though $S$ contains singular terms with respect to $s$, 
 $z = S^a = s^a + t\sum_{i=0}^{a-2}c_{a-i}s^i$ never contains such singular terms, 
 and singular terms in $w = S^b+S^{b+1}g_0(S)$ originate from terms of the form 
 $s^b(\sum_{l=2}^{a}\frac{c_l}{s^l})^i$ for some $i$.
Since we have $c_i\in t^i\Bbb C[[t]]$, singular terms do not appear 
 until we consider deformations of order $b+1$ with respect to $t$.
It follows that, up to this order, the curve defined on the punctured neighborhood $\mathring U_p$
 by the parameterization $(z, w)=(S^a, S^b+S^{b+1}g_0(S))$ actually extends to $p$.
Moreover, the extended curve has the same image as the original $\varphi$ around $p$.
Thus, we can glue these locally defined curves around the singular points of $\varphi$ 
 and the image of $\varphi$ away from the singular points
 into a global curve.
However, since $S$ is a singular parameter around $p$, the domain curve $C$
 must be
 nontrivially deformed in general.
Therefore, we have a nontrivial deformation of $\varphi$ while the image remains the same.

\begin{rem}
The same type of deformations appeared in \cite{AC, AC2} in the case of first order deformations
 and played crucial role in \cite{DS}.
In some sense, our argument extends it as far as possible.
\end{rem}

As we noted above, as long as there is no singular term in $(S^a, S^b+S^{b+1}g_0(S))$, 
 there is no obstruction to glue locally defined deformations, if we deform the domain curve $C$
 suitably.
Our study of obstructions begins when a singular term appears in $(S^a, S^b+S^{b+1}g_0(S))$.
Assume that we have constructed a $k$-th order deformation $\varphi_k$ of $\varphi$ in the above way, 
 and that a singular term appears in $(S^a, S^b+S^{b+1}g_0(S))$ at the order $t^{k+1}$
 at some singular point $p$ of $\varphi$.
In this case, although we cannot extend the curve defined on a punctured neighborhood $\mathring U_p$
 to the whole neighborhood $U_p$ of $p$, we can construct a $k+1$-th order deformation of 
 $\varphi_k$ on $U_p$ simply by removing those singular terms (see Proposition \ref{prop:obstrep}). 
On the other hand, away from singular points of $\varphi$, the map $\varphi_k$ is locally the same as 
 $\varphi$, so we can take $\varphi$ itself as a $k+1$-th order local deformation.
The advantage of this construction is that we can explicitly compute the obstruction cocycle.
For even higher order deformations, this construction allows us to compute the leading terms of the obstruction
 completely explicitly.
However, the actual computation is rather subtle, since when we compare local deformations at higher orders, 
 full non-linearity of various coordinate changes comes into play. 
This will be done in Section \ref{subsec:calculation of o_{N+1}}.
The main result here is Proposition \ref{prop:obsteval}.
It asserts that if a certain system of polynomial
 equations 
 written in terms of the polynomials $f_{b+i}^{(b)}$, $i = 1, \dots, a-1$ (see below, or Lemma \ref{lem:coeff}), 
 has a solution, 
 then, at each order of deformation,
 the obstruction can be taken sufficiently small so that we can apply perturbation method explained below
 to eliminate the obstruction.
The system of polynomial equations mentioned above is the equations $\{\star_{\eta}\}$ in Definition \ref{def:star}, 
 where $\eta$ parameterizes a basis of the dual space of the obstructions.
In particular, each $\eta$ naturally couples with the obstruction cocycle, and if all of these couplings are zero, 
 the map deforms.
The equations $\{\star_{\eta}\}$ imply that the leading term of such a coupling vanishes.

Now, assume that we have constructed an $l$-th order deformation $\varphi_l$ of $\varphi$
 with $l\geq k$.
The obstruction to deforming $\varphi_l$ one step further can be computed as in the last paragraph, 
 and usually it does not vanish.
This means that the map $\varphi_l$ does not deform any more.
Therefore, we need to find another path.
So, we return to some lower order deformation $\varphi_{l'}$, $l'<l$, and attempt to
 reconstruct the deformations $\overline{\varphi}_{m}$, $m>l'$, in a different way so that 
\begin{itemize}
\item we can deform it at least up to the order $t^l$, and
\item the new map $\overline{\varphi}_{l}$ has vanishing obstruction, so that we can continue the deformation
 up to higher orders. 
\end{itemize}
We will achieve this in a way explained below.

As we noted above, 
 the leading terms of the obstruction are expressed in terms of polynomials 
 of the coefficients $c_l$ introduced above.
These polynomials are written as $f_{b+i}^{(b)}$, $i = 1, \dots, a-1,$ in the main text (see Lemma \ref{lem:coeff}).
Thus, the study of the obstruction is reduced to controlling the value of these polynomials.
For this, in addition to the conditions $\{\star_{\eta}\}$ above, 
 we also need to assume a suitable transversality property for these polynomials,
 referred to as the condition ${\rm (T)}$ in Definition \ref{def:perturb}.
Specifically, we require a solution $\{c_l\}$ to the system of equations $\{\star_{\eta}\}$
 that also fulfills the condition ${\rm (T)}$.
In fact, this is the assumption of the main theorem (Theorem \ref{thm:main}).
Fortunately, this condition turns out to be relatively mild, as explained in the argument following Theorem \ref{thm:1}.

When we control the values of $f_{b+i}^{(b)}$ using the above transversality assumption, 
 we need to change the value of $c_l\in t^l\Bbb C[[t]]$.
Here, to change the values of $f_{b+i}^{(b)}$ at some order of $t$, we need to change $c_l$
 at lower orders of $t$.
This amounts to changing the map $\varphi_{l'}$ for some $l'<l$ using the notation in the last paragraph.
The difficulty is that during the reconstruction of deformations $\overline{\varphi}_{m}$, $m>l'$, 
 the same problem might happen.
Namely, some non-vanishing obstruction appears, and to deal with it, we might need to return
 to even lower deformation $\varphi_{l''}$, $l''<l'$.
Also, even if we could construct a new map $\overline{\varphi}_l$, 
 it is unclear that the obstruction to deforming it vanishes.

We can overcome this difficulty again using deformations which do not change the image.
Namely, we will show that we can modify the coefficients $c_l$ so that the resulting maps 
 $\overline{\varphi}_{m}$ has the same image as $\varphi_m$. 
Then, since the maps $\varphi_m$, $l'\leq m \leq l$ already exists, it is easy to see that 
 the obstruction to deforming $\overline{\varphi}_{m}$ ($l'\leq m < l$) one step further vanishes.
Thus, we obtain a map $\overline{\varphi}_{l}$, which also has the same image as 
 $\varphi_l$.
In this situation, we can compare the obstructions to deforming $\varphi_l$ and
 $\overline{\varphi}_{l}$, and since the change of the coefficients $c_l$ was taken so that
 it cancels the original obstruction to deforming $\varphi_l$, the obstruction to 
 deforming $\overline{\varphi}_{l}$ vanishes.

Again, the computation is rather subtle.
The primary reason is that when we calculate the obstruction cocycle, we compare the local deformations
 of $\overline\varphi_m$, which are defined using local coordinates.
However, altering the values of $c_l$ also changes the coordinate $S$ introduced above, making
 the determination of the collect values of $c_l$ a highly nontrivial task.

In any case, we obtain a new deformation $\overline{\varphi}_{l+1}$ of $\varphi$, which 
 coincides with $\varphi_l$ over $\Bbb C[t]/t^{l+1}$.
To deform $\overline{\varphi}_{l+1}$, again we need to return to some $\overline{\varphi}_{l''}$
 and repeat the argument.
In this way, we can construct a deformation $\varphi_N$ of $\varphi$ up to any high order 
 of $t^N$.
Finally, we will show that this can be done in a way that if we need to return to a map $\varphi_{a(N)}$, 
 $a(N)<N$, to deform $\varphi_N$, we have $a(N)\to \infty$ as $N\to\infty$.
Thus, eventually we obtain a projective system $\{\varphi_i\}$ of deformations of $\varphi$
 in which $\varphi_j$ reduces to $\varphi_i$ over $\Bbb C[t]/t^{i+1}$ for $j>i$.
Summarizing, we have the following (see Theorem \ref{thm:main} for the precise statement, in which
 the domain curve is not assumed to be smooth).

\begin{thm}\label{thm:1}
Let $\varphi\colon C\to X$ be a map from a regular complete curve to a regular surface which is 
 birational to the image.
Assume the map $\varphi$ is semiregular.
Let $\{p_1, \dots, p_r\}$ be the singular points of the map $\varphi$.
Then, if there is a set of solutions to the system of equations $(\star_{\eta})$ which also satisfies
 the condition ${\rm (T)}$, 
 there is a deformation of $\varphi$ which deforms the singularity of $\varphi$ at each $p_j$ non-trivially.
\end{thm}

The set of equations $\{\star_{\eta}\}$ depends on the geometry of the curve, 
 and checking the condition of Theorem \ref{thm:main} (the existence of solutions of $\{\star_{\eta}\}$
 satisfying the transversality condition ${\rm (T)}$) for each individual curve is a cumbersome task.
However, it appears that in most cases, a significantly stricter condition holds true, and we can,
 to a large extent, disregard the 
 individual properties of curves.
Here, we outline this point. 
See Section \ref{sec:examples} for details.
First, we note that the dual of the obstruction space is given by $H^0(C, \varphi^*\omega_X(Z))$, 
 where $Z$ is the ramification divisor of the map $\varphi$, as discussed in Sections  \ref{subsec:localobst} and 
 \ref{subsec:surfacecase}.
Due to the semiregularity assumption, sections of $H^0(C, \varphi^*\omega_X)$ do not pair with 
 obstruction classes non-trivially.
Therefore, to calculate obstructions, we need to assess how they interact when paired with 
 meromorphic sections of the pull back of $\omega_X$ to $C$.
In particular, the set $\{\eta\}$, which parameterizes the equations $\{\star_{\eta}\}$,
 forms a basis for the quotient 
 $H^0(C, \varphi^*\omega_X(Z))/H^0(C, \varphi^*\omega_X)$.

Let $p\in C$ be a singular point of the map $\varphi$.
In terms of the notation above, $a$ is the multiplicity of the singularity of the image $\varphi(p)$.
In particular, $a-1$ is the coefficient of $p$ in the ramification divisor $Z$.
On the other hand, $a-1$ is the same as the number of the freedom for deforming the singularity of 
 $\varphi$ at $p$, which is parameterized by $c_l$, $l = 1, \dots, a-1$, introduced above.
Thus, if there are $a-1$ 
 sections of $H^0(C, \varphi^*\omega_X(Z))$ which are singular only at the point $p$ as 
 sections of $H^0(C, \varphi^*\omega_X)$, and form a basis of the quotient
 space $H^0(C, \varphi^*\omega_X((a-1)p))/H^0(C, \varphi^*\omega_X)$, 
 then virtually we cannot expect that the singularity of the map $\varphi$ at $p$ to deform.
To ensure a nontrivial deformation, it is reasonable to assume that 
 the space $H^0(C, \varphi^*\omega_X(Z))$ does not contain such a set of sections for each singular point of
 $\varphi$.
This is what we refer to as condition ${\rm (D)}$ in Definition \ref{def:D}. 
 
On the other hand, we introduce the condition ${\rm (G)}$ in Definition \ref{def:(G)}, and prove
 that under the condition ${\rm (D)}$, the condition ${\rm (G)}$ implies that 
 there is a set of solutions of the system $(\star_{\eta})$ which satisfies the condition ${\rm (T)}$.
Combined with Theorem \ref{thm:1}, this implies the following (which is the same as Theorem
 \ref{thm:gen2}).
\begin{thm}\label{thm:2}
If the conditions ${\rm (D)}$ and ${\rm (G)}$ hold at each singular point of $\varphi$, the map $\varphi$ deforms.
\end{thm}

The advantage of the condition ${\rm (G)}$ is, contrary to the condition ${\rm (T)}$ (combined with the equations
 $(\star_{\eta})$), 
 it does not depend on the properties of the map $\varphi$.
In fact, the condition ${\rm (G)}$ depends only on the pair of positive integers $(a, b)$ and is
 independent of all the other geometry.
Thus, the problem of the existence of deformations of $\varphi$ is reduced to the cohomological calculation
 of checking the condition ${\rm (D)}$ (which is considerably easier than the deformation problem), and
 the condition ${\rm (G)}$, which is completely independent of the deformation problem.
As we noted above, the condition ${\rm (D)}$ is virtually the minimal requirement for the existence of 
 deformations of $\varphi$.
Theorem \ref{thm:2} claims that, if the condition ${\rm (G)}$ holds, then any semiregular map satisfying 
 the condition ${\rm (D)}$ deforms.
In other words, assuming the condition ${\rm (G)}$,
 these maps have almost optimal possible deformation property.

Usually, the condition ${\rm (G)}$ is much stronger than the condition required to apply 
 Theorem \ref{thm:1}.
Fortunately, however, the condition ${\rm (G)}$ seems to hold in almost all cases, though at present
 we do not know how to prove it in general.
We have checked it by computer calculation for the values of $(a, b)$ roughly up to $a+b<30$
 (see Table \ref{table:1}), and found that 
 the condition ${\rm (G)}$ holds except only one case of $a=4, b=6$.
Even in this exceptional case, 2 out of 3 cases of the condition ${\rm (G)}$ hold (see the argument
 at the end of Section \ref{subsec:(G)}).
The condition ${\rm (G)}$ holds if the system of polynomials $f_{b+i}^{(b)}$, $i=1, \dots, a-1$ has
 a transversality property similar to that of a generic system of polynomials.
The exceptional case of $a=4, b=6$ seems to be caused by an accidental factorization of some
 of the polynomials $f_{b+i}^{(b)}$ due to the smallness of the degree.

Another advantage of this result is that we use only a very small part of the data of singularities.
Namely, as we noted above, the condition ${\rm (G)}$ only depends on 
 the pair of numbers $a$ and $b$.
In particular, for $a>2$, once one shows the condition ${\rm (G)}$ for some $a$ and $b$, 
 it applies to infinitely many types of singularities.

On the other hand, if the singularities are of multiplicity two (in other words, double points), 
 we can even deduce the necessary and sufficient condition for the existence of deformations.
In this case, we do not need the conditions ${\rm (D)}$ and ${\rm (G)}$, and we can deduce 
 results considerably stronger than Theorems \ref{thm:1} and \ref{thm:2}.
The reason for this is that for $a=2$, the function $f_{b+i}^{(b)}$ is simply a power of $c_2$, 
 and this fact allows us to manipulate the obstruction much more efficiently than other cases, 
 even without transversality condition.
 
So, assume that the map $\varphi$ has singular points $\{p_1, \dots, p_l\}$ each of which is a double point.
Then, we have the following, which is essentially the same as 
 Theorem \ref{thm:doublepoint} but expressed in a slightly different manner.
\begin{thm}\label{thm:doublepoint}
	The semiregular map $\varphi$ deforms if and only if either one of the following conditions
	holds.
	\begin{enumerate}
		\item There is at least one $p_i$ such that there is no section of 
		$H^0(C, \varphi^*\omega_X(p_i))$ which is not contained in $H^0(C, \varphi^*\omega_X)$.
		\item The set $H^0(C, \bar{\mathcal N}_{\varphi})$ is not zero.
	\end{enumerate}
\end{thm}
Here, $\bar{\mathcal N}_{\varphi}$ is the non-torsion part of the normal sheaf of $\varphi$, 
 see the paragraph following Definition \ref{def:N_varphi}.
 
Now, let us return to Problem \ref{problem:2} mentioned at the beginning of the introduction, 
The above results imply that, if the singularities of $\varphi$ are double points or
 the ones for which the condition ${\rm (G)}$ is checked, 
 we can deform $\varphi$ until the condition ${\rm (D)}$ is violated, or a singularity which does not satisfy
 the condition ${\rm (G)}$ (presumably it happens only when $a=4, b=6$) appears.
Since the condition ${\rm (D)}$ is virtually the minimal requirement for the existence of deformations
 as we mentioned above, it follows that on any surface, Problem \ref{problem:2} has an almost optimal answer under
 condition ${\rm (G)}$.
 
\subsection{Notation}\label{subsec:notation}
We will work in the complex analytic category.
In the body of the paper, we will study non-constant maps $\varphi\colon C\to X$ from a curve $C$
 to a smooth complex surface $X$ and their deformations.
Here, the domain $C$ will deform but the target $X$ is fixed.
Usually, we assume the curve $C$ is integral.
A deformation of $\varphi$ over $\Spec\Bbb C[t]/t^{k+1}$ will be written as 
 $\varphi_k\colon C_k\to X\times\Spec\Bbb C[t]/t^{k+1}$, or often as
 $\varphi_k\colon C_k\to X$ for notational simplicity.
Let $p\in C$ be a regular point and $s$ be a local parameter around $p$.
Let $\{z, w\}$ be a local coordinate system on a neighborhood of $\varphi(p)$ on $X$.
Then, the pull back of the coordinates $z$ and $w$ by $\varphi$ are functions $z(s), w(s)$ of $s$, 
 and we call $(z(s), w(s))$ a local parameterization of the map $\varphi$.
By the image of a map $\varphi$ or $\varphi_k$, we mean the analytic locally ringed space
 with the annihilator structure, see \cite[Chapter I, Definition 1.45]{GLS}.
That is, if $U$ is an open subset of $C_k$ with the induced structure of an analytic locally ringed space, 
 and $V$ is an open subset of $X$ such that 
 $\varphi_k(U)$ is closed in $V$, we associate the structure sheaf 
\[
\mathcal O_V/\mathcal{A}nn_{\mathcal O_V}((\varphi_k)_*\mathcal O_U)
\]
 to the image $\varphi_k(U)$.

\section{Localization of obstructions}

\subsection{Meromorphic differential forms and cohomological pairings}\label{subsec:paring}
We begin with giving a presentation of a \v{C}ech cocycle in a way suited to our purpose.
Let $C$ be a complete integral curve and $\mathcal L$ be 
 an invertible sheaf on it.
Let $\{q_1, \dots, q_s\}$ be the set of singular points on $C$.
Let $\{p_1, \dots, p_e\}$ be a non-empty set of non-singular points on $C$.
Take an open covering $\{U_{1}, \dots, U_{m}\}$ of $C$.
We assume each of $p_i$ and $q_i$ is contained in a unique element of this covering, 
 and denote them as $U_{p_i}$ and $U_{q_i}$, respectively.
We also assume that 
 the normalization of each $U_j$ is a disc when it does not contain any $q_i$, whereas 
 the normalization of $U_{q_i}$ is 
 a disjoint union of $b_{q_i}$ discs, where $b_{q_i}$ is the number of branches of $C$ at $q_i$.
Here, a disc means an analytic subset which is analytically isomorphic to 
$D = \{z\in \Bbb C\;|\; |z|<1\}$.
We also assume $U_i\cap U_j$ is a disc or empty when $i\neq j$.




We associate a meromorphic section $\xi_j$ of $\mathcal L|_{U_j}$ with each $U_j$ 
 which can have a pole only at $\{p_1, \dots, p_e\}$.
In particular, if $U_j$ does not contain any of $\{p_1, \dots, p_e\}$, 
 $\xi_j$ is a regular section of $\mathcal L|_{U_j}$.
It follows that $\xi_{ij}=\xi_i-\xi_j$ is a section of 
 $\mathcal L|_{U_i\cap U_j}$, since $U_i\cap U_j$ does not contain any $p_k$.
Thus,
 the set of sections $\{\xi_i\}$ determines a \v{C}ech 1-cocycle 
 $\{\xi_{ij}\}$ with values in $\mathcal L$
 for the covering $\{U_i\}$ mentioned above.
Any class in $H^1(C, \mathcal L)$ can be represented 
 in this way (see Proposition \ref{prop:pair} below).

Let $\omega_{C}$ be the dualizing sheaf of $C$.
This sheaf is defined using appropriate meromorphic differential forms
 which have poles only at the (normalization of) singular points of $C$.
Namely, $\omega_C$ is given
 by the sheaf of Rosenlicht differentials \cite{R}. 
We do not need details of it in this paper, and we omit its definition.
See \cite{R} or other expositions, for example, \cite[Chapter VIII, Section 1]{AK} or \cite[Section II.6]{BHPV}. 
Let $\psi$ be a global section of $\mathcal L^{\vee}\otimes \omega_{C}$.
Then, $\{\xi_{ij}\}$ and $\psi$ make a natural pairing.
The value of this pairing is given as follows.
Namely, on $U_j$, the fiberwise pairing 
 between $\xi_j$ and $\psi$ gives a meromorphic 
 section $(\psi, \xi_j)$ of $\omega_C|_{U_j}$.
If $p_i\in U_{p_i}$, $(\psi, \xi_{p_i})$ may have a pole at $p_i$, and
 let $r_{p_i}$ be its residue.
The section $(\psi, \xi_{q_i})$ may also have a pole at 
 a singular point $q_i$ of $C$ due to the pole of $\omega_C$.
However, the contribution to the pairing from such a pole vanishes due to the defining property of 
 Rosenlicht differentials.

Then, the following is a special case of 
 \cite[Proposition 10]{N}.
\begin{prop}\label{prop:pair}
	\begin{enumerate}
		\item Any cohomology class in $H^1(C, \mathcal L)$ can be represented
		by some set $\{\xi_i\}$ of local meromorphic sections on open subsets $\{U_i\}$ as above.
		\item 	
		The pairing between $\{\xi_{ij}\}$ and $\psi$
		is given by
		\[
		\langle \psi, \{\xi_{ij}\}\rangle = \sum_{i=1}^e r_{p_i}.
		\]
		This gives the natural nondegenerate pairing between $H^1(C, \mathcal L)$
		and its dual space $H^0(C, \mathcal L^{\vee}\otimes \omega_{C})$.
	\end{enumerate}
\end{prop}

\subsection{Normal sheaf of a map}\label{subsec:localobst}
Let $\varphi\colon C\to X$ be a map from a complete integral curve to a smooth surface.
Let $\{p_1, \dots, p_e\}$ be the set of points where $\varphi$ is not a local embedding.
We assume each $p_i$ is a regular point of $C$.
We take a covering $\mathcal U = \{U_1, \dots, U_m\}$ of $C$ 
 as in Section  \ref{subsec:paring}.
In particular, for each point $p\in \{p_1, \dots, p_e\}$, there is a unique open subset in $\mathcal U$
 containing $p$.
We write it by $U_p$ as before.
Let $q\in C\setminus \{p_1, \dots, p_e\}$ be any point.
There is a neighborhood $U_q$ of $q$ in $C$ on which $\varphi$ is an isomorphism onto its image.
Let $V$ be a suitable open subset of $X$ such that $U_q$ is one of the connected components of 
 $\varphi^{-1}(V)$.
There is a usual normal sheaf of the image of $U_q$ defined by 
 $\mathcal N_{\varphi(U_q)} = \mathcal O_V(\varphi(U_q))|_{\varphi(U_q)}$.
We will regard this also as a sheaf on $U_q$ in an obvious way.

On the other hand, on a neighborhood $U_{p_i}$ of $p_i$, consider the sheaf $\mathcal N_{\varphi|_{U_{p_i}}}$
 defined by the exact sequence
\[
0\to \mathcal T_{U_{p_i}}\to \varphi_{U_{p_i}}^*\mathcal T_X\to \mathcal N_{\varphi|_{U_{p_i}}}\to 0.
\]
Here, $\mathcal T_{U_{p_i}}$ and $\mathcal T_X$ are the tangent sheaves.

The sheaves $\mathcal N_{\varphi(U_q)}$ and $\mathcal N_{\varphi|_{U_{p_i}}}$
 are naturally isomorphic on their intersection.
Namely, they are naturally identified with the pull back of the normal sheaf of the image 
 $\varphi(U_q\cap U_{p_i})$.
Thus, we obtain a global sheaf on $C$.
\begin{defn}\label{def:N_varphi}
We denote the sheaf on $C$ obtained in this way by $\mathcal N_{\varphi}$.
\end{defn}

The sheaf $\mathcal N_{{\varphi}}$ has torsion at singular points $\{p_1, \dots, p_e\}$ of ${\varphi}$.
In particular, there is an exact sequence of sheaves on $C$,
\[
0\to \mathcal H_{\varphi}\to \mathcal N_{{\varphi}}\to \bar{\mathcal N}_{{\varphi}}\to 0,
\]
 where $\mathcal H_{\varphi}$ is a torsion sheaf and $\bar{\mathcal N}_{{\varphi}}$ is locally free.
The sheaf $\bar{\mathcal N}_{{\varphi}}$ is also described on a neighborhood $U_{p_i}$ of $p_i$ as follows.
\begin{lem}\label{lem:equisingnormal}\cite[Section 3.4.3]{S}
There is an exact sequence
\[
0\to \mathcal T_{U_{p_i}}(Z_{p_i})\to \varphi|_{U_{p_i}}^*\mathcal T_X\to \bar{\mathcal N}_{{\varphi}}|_{U_{p_i}}\to 0.
\]
Here, $Z = (d\varphi)$ is the ramification divisor and $Z_{p_i}$ is its restriction to $p_i$.
 \qed
\end{lem}

The set $\{p_1, \dots, p_e\}$ is the support of $Z$.
Given a $k$-th order deformation ${\varphi}_k$ of $\varphi$ for a non-negative integer $k$, 
 the obstruction to deforming it one step further is represented by a cocycle
 defined by taking the difference of local deformations on $U_i$.
The obstruction class belongs to the cohomology group $H^1(C, \bar{\mathcal N}_{{\varphi}})$.

\subsection{Obstructions to deformations of singular curves on surfaces}\label{subsec:surfacecase}
By Lemma \ref{lem:equisingnormal}, we have the following.
\begin{lem}
We have an isomorphism
\[
\bar{\mathcal N}_{{\varphi}}\cong \varphi^*\omega_X^{-1}\otimes \omega_{C}(-Z),
\]
 of sheaves on $C$, 
 where $\omega_X$ is the canonical sheaf of $X$ and
 $\omega_C$ is the dualizing sheaf of the possibly singular reduced curve $C$.
\end{lem}
\proof
Outside the set $\{p_1, \dots, p_e\}$, this follows by the adjunction.
At $\{p_1, \dots, p_e\}$, this follows from Lemma \ref{lem:equisingnormal}.
These isomorphisms glue, since on the intersection
 of open subsets, $\bar{\mathcal N}_{\varphi}$ is naturally isomorphic to the
 normal sheaf of the image of $\varphi$.\qed\\

By the Serre duality, we have
\[
H^1(C, \bar{\mathcal N}_{{\varphi}})\cong H^0(C, \varphi^*\omega_X(Z))^{\vee}.
\]
Let $p$ be a singular point of $\varphi$.
Taking a suitable coordinate system
 $\{z, w\}$ around $\varphi(p)$ on $X$, 
 and a suitable local parameter $s$ on the open subset $U_p$ of $C$,
 the map $\varphi$ can be parameterized as
\[
(z, w) = (s^a, \; s^b + s^{b+1}g_0(s)), 
\]
 where $g_0(s)$ is an analytic function, 
 $a-1$ is the vanishing order of $d\varphi$ at $p$
 (in other words, the coefficient of $p$ in the ramification divisor $Z$),
 and $b$ is a positive integer larger than $a$.
See, for example, \cite[Chapter I, Corollary 3.8]{GLS}.
By a coordinate change on $X$, we can assume $b$ is not a multiple of $a$.
We assume this hereafter.


Its first order deformation is, up to a change of the parameter $s$, 
 of the following form:
\[
(z, w) = (s^a  + t\sum_{i=0}^{a-2}c_{a-i}s^i , \; s^b + s^{b+1}g_0(s) + tg_1(s)),
\]
 where $c_i$ is a complex number, and $g_1$ is an analytic function.
Note that since $U_p$ is nonsingular, its infinitesimal deformation is trivial.
Thus, if we have a $k$-th order deformation 
 $\varphi_k\colon C_k\to X\times \Spec\Bbb C[t]/t^{k+1}$ of $\varphi$, 
 the restriction of it to $U_p$ gives a map between locally ringed spaces
 from $U_p\times \Spec \Bbb C[t]/t^{k+1}$ to $X\times \Spec\Bbb C[t]/t^{k+1}$.
We call any function on $U_p\times \Spec \Bbb C[t]/t^{k+1}$ which reduces to the given parameter $s$
 on $U_p$ over $\Bbb C[t]/t$ a parameter on $U_p\times \Spec \Bbb C[t]/t^{k+1}$, and 
 we will denote it again by $s$.
Similarly, we will also consider a parameter on $U_p$ defined over $\Bbb C[[t]]$.

 
 
 
The part $t\sum_{i=0}^{a-2}c_{a-i}s^i$, which can be seen as an element of the $(a-1)$-dimensional 
 vector space $V_p:=\{c_a+c_{a-1}s+ \cdots + c_2s^{a-2}\,|\, c_i\in\Bbb C\}$, corresponds 
 to the torsion part of $\mathcal N_{\varphi}$ (precisely speaking, it is the sum of
 a torsion element and a non-torsion element, see Example \ref{ex:1} below).
The part $g_1$ corresponds to a section of the sheaf $\bar{\mathcal N}_{{\varphi}}$.
Similarly, given a $k$-th order deformation $\varphi_k$ of $\varphi$ whose restriction to
 $U_{p}\times\Spec\Bbb C[t]/t^{k+1}$ is parameterized as
\[
(z, w) = (s^a  + \sum_{j=1}^k\sum_{i=0}^{a-2}t^jc_{a-i, j}s^i , \; s^b + s^{b+1}g_0(s)+\sum_{j=1}^kt^jg_{j}(s)),
\]
 its local $(k+1)$-th order deformations are given by
\[
(z, w) = (s^a  + \sum_{j=1}^{k+1}\sum_{i=0}^{a-2}t^jc_{a-i, j}s^i , \; s^b + s^{b+1}g_0(s)+ \sum_{j=1}^{k+1}t^jg_{j}(s)),
\]
 where $c_{a-i, j}$ is a complex number and $g_j(s)$ is an analytic function.
Here, the parameter $s$ on $U_p\times \Spec \Bbb C[t]/t^{k+2}$ is chosen so that it reduces to a given parameter on 
 $U_p\times \Spec \Bbb C[t]/t^{k+1}$.

\begin{example}\label{ex:1}
In the case of the map $\varphi\colon \Bbb C\to \Bbb C^2$ given by $s\mapsto (s^2, s^3)$,
 whose image is the ordinary cusp, 
 the sheaf $\mathcal N_{\varphi}$
 is given by $\mathcal O_{\Bbb C}\langle\partial_z, \partial_w\rangle/(2s\partial_z + 3s^2\partial_w)$,
 where $\{z, w\}$ is the standard coordinate system on $\Bbb C^2$, and $\partial_z$, $\partial_w$ are
 the standard generators of the tangent sheaf of $\Bbb C^2$.
In this case, $t\sum_{i=0}^{a-2}c_{a-i}s^i = tc_2$ and it corresponds to 
 the section $c_2\partial_z$ of $\mathcal N_{\varphi}$.
It is the sum of a torsion element $c_2(\partial_z+\frac{3s}{2}\partial_w)$ and a non-torsion element
 $-\frac{3c_2s}{2}\partial_w$.
\end{example}

Let $k$ be a non-negative integer and 
 let $\varphi_k$ be a $k$-th order deformation of $\varphi$.
Let $\mathcal G_{k+1}$ be the group of automorphisms of $\mathcal O_{U_p}\times \Bbb C[t]/t^{k+2}$
 which are the identity over $\Bbb C[t]/t^{k+1}$,
 where $\mathcal O_{U_p}$ is the sheaf of functions on the neighborhood $U_p$ of $p$, 
 considered as an analytic locally ringed space.
The group $\mathcal G_{k+1}$ acts on the ringed space $U_p\times \Spec\Bbb C[t]/t^{k+2}$
 as automorphisms.
Consequently, it also acts on
 the set of $(k+1)$-th order local deformations of $\varphi$
 on $U_p$ which restricts to $\varphi_k$ over $\Bbb C[t]/t^{k+1}$.
Let $\varphi_{k+1}|_{U_p}, \varphi_{k+1}'|_{U_p}$ be such local deformations.
We call them equivalent if there is some $g\in \mathcal G_{k+1}$
 such that $\varphi_{k+1}|_{U_p} = \varphi_{k+1}'|_{U_p}\circ g$.
In the following statement, we use the same notation as above.
In particular, $V_p=\{c_a+c_{a-1}s+ \cdots + c_2s^{a-2}\,|\, c_i\in\Bbb C\}$.

\begin{prop}\label{prop:loccalnormal}
	Given a $k$-th order analytic deformation $\varphi_k$ of $\varphi$ on $U_p$, 
	 the set of equivalence classes of $(k+1)$-th order analytic deformations of $\varphi$ on $U_p$ 
	 which reduce to $\varphi_k$ 
	 is naturally isomorphic to 
	 the set of sections of the sheaf 
\[
\mathcal O_{U_p}\cdot \partial_w\oplus i_*V_p\cdot\partial_z
\]
 on $U_p$.
Here, $\{\partial_z, \partial_w\}$ is the pullback of the 
 natural basis of the tangent bundle of the coordinate neighborhood
 $\{z, w\}$ on $X$ by $\varphi$.
Also, the vector space $V_p$ is regarded as the constant sheaf on the point $p$, and
 $i\colon \{p\}\to U_p$ is the inclusion.
\end{prop} 
\proof
We define a map from the set of equivalence classes of deformations to the set of sections of  
 $\mathcal O_{U_p}\cdot \partial_w\oplus i_*V_p\cdot\partial_z$ by
\[
\varphi_{k+1}|_{U_p}\mapsto g_{k+1}(s)\cdot \partial_w + g_{1, k+1}(s)\cdot\partial_z,
\]
 where $g_{1, k+1}(s) = \sum_{i=0}^{a-2}c_{a-i, k+1}s^i$.
This is well-defined since if we apply an element of $\mathcal G_{k+1}$, it 
 changes the coordinate $s$ to another one of the form $s + t^{k+1}b_{k+1}(s)$, 
 where $b_{k+1}(s)$ is a holomorphic function.
Then, it breaks the form of $g_{1, k+1}(s)$.
Namely, it produces a non-zero coefficient of $s^{a'}$, $a'\geq a-1$.

Conversely, given a section $\alpha(s)\cdot \partial_w + \beta(s)\cdot\partial_z$
of $\mathcal O_{U_p}\cdot \partial_w\oplus i_*V_p\cdot\partial_z$,
 we define a deformation of $\varphi_k$ by taking 
$g_{k+1}(s) = \alpha(s)$ and $g_{1, k+1}(s) = \beta(s)$.
Then, we obtain the inverse mapping. \qed\\

The sheaf $\mathcal O_{U_p}\cdot \partial_w\oplus i_*V_p\cdot\partial_z$
 is naturally isomorphic to the sheaf $\mathcal N_{\varphi}$ restricted to $U_p$, though the direct sum decomposition 
 does not coincide with
 the more natural $\mathcal N_{\varphi}\cong \mathcal H_{\varphi}\oplus\bar{\mathcal N}_{\varphi}$
 appeared in the previous subsection, as in Example \ref{ex:1}.
However, the part $\mathcal O_{U_p}\cdot \partial_w$ is naturally isomorphic to $\bar{\mathcal N}_{\varphi}$.

%

\subsubsection{Explicit presentation of the obstruction cocycle in the first non-trivial case}\label{subsubsec:obstcocycle}
By Proposition \ref{prop:loccalnormal}, if a local deformation of $\varphi$ associated with
 the summand $i_*V_p\cdot\partial_z$ can be extended globally after modifying
 by sections of $\bar{\mathcal N}_{\varphi}$ if necessary, it gives a non-trivial deformation of $\varphi$.
Moreover, it deforms the singularity of $\varphi$ at $p$, since it lowers the multiplicity $a$ of 
 the singularity. 
So, we study the obstruction to the deformations
 associated with this part.
Take a covering $\mathcal U = \{U_1, \dots, U_m\}$ of $C$ as before.
Recall that for each singular point $p\in \{p_1, \dots, p_e\}$ of $\varphi$, there is a unique open subset $U_p$
 belonging to $\mathcal U$ which contains $p$.
Let $s$ be a parameter on $U_p$.
We also denote by $s$ a parameter on $U_p$ defined over $\Bbb C[[t]]$ which reduces to the given one over $\Bbb C[t]/t$.
Let $S$ be a local parameter over $\Bbb C[[t]]$ 
 on a punctured neighborhood $\mathring U_p = U_p\setminus\{p\}$ of $p$ in $C$ which satisfies
\[
S^a = s^a  + t\sum_{i=0}^{a-2}c_{a-i}s^i.
\]
Given a polynomial $\beta(s)$, since one of the $a$-th roots of $s^a+t\beta(s)$ is given by 
 $s(1 + \sum_{i=1}^{\infty}\prod_{j=0}^{i-1}(\frac{1}{a}-j)\frac{1}{i!}(t\frac{\beta(s)}{s^a})^i)$, 
 we can take
\[
S = s(1 + \sum_{i=1}^{\infty}\prod_{j=0}^{i-1}(\frac{1}{a}-j)\frac{1}{i!}(t\sum_{l=2}^{a}\frac{c_l}{s^l})^i).
\]
Consider a first order local deformation of $(z, w) = (s^a, s^b+s^{b+1}g_0(s))$ (now this is considered over $\Bbb C[t]/t$)
 of the form
\[
(z, w) = (S^a , \; S^b+S^{b+1}g_0(S)),
\]
 defined on $\mathring U_p$.
Evidently, we have $S^a = s^a  + t\sum_{i=0}^{a-2}c_{a-i}s^i$, which holds up to any order with respect to $t$.
One also sees that $S^b+S^{b+1}g_0(S)$, expanded in terms of $s$ as a function defined over $\Bbb C[t]/t^2$,
 can be extended to $p$.
In other words, the expansion does not include
 any term with negative powers of $s$.  
This means that the deformed curve is also extended to $U_p$, and is still locally defined by the same equation
 as that of $(z, w) = (s^a, s^b+s^{b+1}g_0(s))$, now considered over
 $\Bbb C[t]/t^2$.
On an open subset $U_j$ which does not contain a singular point of $\varphi$, 
 we can regard $\varphi$ itself as a local deformation of $\varphi$ over $\Bbb C[t]/t^2$.

The obstruction to the existence of a global deformation is given by the difference of local deformations
 on overlaps like $U_p\cap U_j$.
If we take local deformations as in the previous paragraph, the difference gives the zero section of 
 the sheaf $\bar{\mathcal N}_{\varphi}$ on each $U_p\cap U_j$.
Thus, there is no obstruction to the existence of a first order deformation.
Note that a deformation obtained in this way has the image 
 whose local defining equations are the same as those of $\varphi$, 
 but regarded as equations over $\Bbb C[t]/t^2$.
\begin{rem}\label{rem:image-map}
Note that although the image is the same, the map over $\Bbb C[t]/t^2$ may not be a trivial deformation
 of $\varphi$, since the difference between local deformations can give a nontrivial cocycle with values in 
 the tangent sheaf of $C$.
In this case, the domain curve of the deformed map is a nontrivial deformation of $C$, 
 and consequently, the map is also
 nontrivially deformed.
\end{rem}

The same holds until a singular term appears in the expansion of $S^b + S^{b+1}g_0(S)$ at some 
 $p\in \{p_1, \dots, p_e\}$.
Let $t^k$ be the minimal order where such a term appears.
Let $\varphi_{k-1}$ be the $(k-1)$-th order deformation of $\varphi$ obtained in the way described above.
Thus, the image of $\varphi_{k-1}$ is the same as that of $\varphi$.
In particular, on $U_p$, the map $\varphi_{k-1}$ has the parameterization $(z, w) = (S^a, S^b+S^{b+1}g_0(S))$.

Now, we take local deformations of $\varphi_{k-1}$ as follows.
On an open subset $U_i$, which does not contain a singular point $p$ at which a singular term 
 appears in $S^b+S^{b+1}g_0(S)$, 
 we take the trivial local deformation of $\varphi_{k-1}$
 as above,
 whose image is given by the same defining equation as the image of $\varphi$ restricted to $U_i$,
 but considered over $\Bbb C[t]/t^{k+1}$.
On the open subset $U_p$ which contains $p$, we have a parameterization
 of $\varphi_{k-1}$ given by $(z, w) = (S^a , S^b+S^{b+1}g_0(S))$ defined over $\Bbb C[t]/t^{k}$.
If we expand $S^b+S^{b+1}g_0(S)$ after substituting 
 $S = s(1 + \sum_{i=1}^{\infty}\prod_{j=0}^{i-1}(\frac{1}{a}-j)\frac{1}{i!}(t\sum_{l=2}^{a}\frac{c_l}{s^l})^i)$,
 it contains singular terms of the order $t^k$, and we take a local deformation on $U_p$ by simply
 discarding these singular terms.
We write it by $(z, w) = ({S}^a, \overline{S^b+S^{b+1}g_0(S)})$.

The difference between the local deformations gives the obstruction \v{C}ech 1-cocycle to 
 deforming $\varphi_{k-1}$.
In terms of the formulation in Section \ref{subsec:paring}, this is described as follows.
By taking a refinement of the open covering $\{U_i\}$ if necessary, we assume that 
 if each $U_i$ and $U_j$ ($i\neq j$) contains a singular point of $\varphi$,
 we have $U_i\cap U_j = \emptyset$.
\begin{prop}\label{prop:obstrep}
The obstruction cocycle associated with the local deformations of $\varphi_{k-1}$ described above
 is represented by the following set of local meromorphic sections of $\bar{\mathcal N}_{\varphi}$
 on each $U_i$.
Namely, give the zero section to all the open subsets $U_i$ which do not contain a
 singular point of the map $\varphi$.
On the open subset $U_p$ containing a singular point $p$, 
 attach a meromorphic section 
 $(\overline{S^b+S^{b+1}g_0(S)} - (S^b+S^{b+1}g_0(S)))\partial_w$
 (precisely speaking, the coefficient of $t^k$ of it) using the notation above.
\end{prop}
Note that if a singular term does not appear in the expansion of $S^b+S^{b+1}g_0(S)$ at some $q\in\{p_1, \dots, p_e\}$,
 the meromorphic section attached to $U_q$ is the zero section.

\proof
The coefficient of $t^k$ of the meromorphic function $S^b+S^{b+1}g_0(S)$, which we write by 
 $(S^b+S^{b+1}g_0(S))_{k}$, 
 corresponds to a holomorphic section $(S^b+S^{b+1}g_0(S))_{k}\partial_w$ of $\bar{\mathcal N}_{\varphi}$
 on the punctured disc $\mathring U_p$ under the correspondence 
 in Proposition \ref{prop:loccalnormal} (precisely speaking, its slightly modified version over $\mathring U_p$).
Note that this corresponds to a deformation of $\varphi_{k-1}$ on $\mathring U_p$
 whose image is defined by the same equation as the image of $\varphi$.
In particular, on the intersection $\mathring U_p\cap U_j$ with another open subset, 
 the images of the local deformation given by $(S^b+S^{b+1}g_0(S))_{k}\partial_w$ on $\mathring U_p$
 and that given by the zero section on $U_j$
 coincide, since they are both defined by the same equation.
It follows that the difference on $\mathring U_p\cap U_j$
 between these local deformations
 belongs to the tangent sheaf of $\varphi(C)$, which is zero in $\bar{\mathcal N}_{\varphi}$. 

Therefore, the section of $\bar{\mathcal N}_{\varphi}$ on the intersection $U_p\cap U_j$
 given by the difference of the local deformations
 corresponding to
 $(z, w) = ({S}^a, \overline{S^b+S^{b+1}g_0(S)})$ on $U_p$ and to the zero section on $U_j$
 equals to $(\overline{S^b+S^{b+1}g_0(S)} - (S^b+S^{b+1}g_0(S)))\partial_w$
 (precisely speaking, the coefficient of $t^{k}$ of 
 it).
This proves the claim.\qed\\

In general, such a section makes a non-trivial residue pairing with elements of 
 $H^0(C, \bar{\mathcal N}^{\vee}_{{\varphi}}\otimes\omega_C)$, which calculates a contribution to the obstruction
 to deforming $\varphi_{k-1}$, by Proposition \ref{prop:pair}.
At higher orders, the calculation of a representative of the obstruction class is more involved.
See Section \ref{subsec:calculation of o_{N+1}}.
In the rest of this paper, we study this contribution to the obstruction in more detail.

\subsubsection{Coefficients of singular terms}\label{subsubsec:regularobst}
Take a deformation of the form 
 $(z, w) = (S^a , S^b+S^{b+1}g_0(S))$
 of $(z, w) = (s^a, s^b+s^{b+1}g_0(s))$ 
 around a singular point $p$ of $\varphi$
 as before.
It is primarily defined only over $\mathring U_p$.
Here, $S$ satisfies $S^a = s^a + t\sum_{i=0}^{a-2}c_{a-i}s^i$.
Although we took $c_i$ to be complex numbers in the above argument, in general we can take them 
 to be elements of $\Bbb C[[t]]$.
For notational simplicity, we rewrite $S^a$ as $S^a = s^a + \sum_{i=0}^{a-2}c_{a-i}s^i$, 
 where $c_i\in t\Bbb C[[t]]$.
Moreover, it will be convenient to take {$c_i\in t^i\Bbb C[[t]]$} in view of certain homogeneity property, 
 see Definition \ref{def:homog}.
We assume this {hereafter}.

Recall that we can assume $b$ is not a multiple of $a$.
Then,
 we have 
\[
S^{b} = s^{b}(1+ \sum_{k=2}^a \frac{c_k}{s^k})^{\frac{b}{a}}
 = s^{b}(1 + \sum_{i=1}^{\infty}\prod_{l=0}^{i-1}(\frac{b}{a}-l)\frac{1}{i!}(\sum_{k=2}^{a}\frac{c_k}{s^k})^i).
\]
We write this in the form
\[
{S^{b} = s^{b}(1 + \sum_{i=1}^{\infty} f_{i}^{(b)}({\bf c})\frac{1}{s^i})},
\]
 where ${\bf c} = (c_2, \dots, c_a)$.
The coefficients $f_i^{(b)}$ are given as follows.
\begin{lem}\label{lem:coeff}
We have 
\[
f_{b+j}^{(b)}({\bf c}) = \sum_{\lambda\in\mathcal P(b+j; [2, a])} 
\begin{pmatrix}
\frac{b}{a}\\
\lambda(2) \; \cdots \; \lambda(a)
\end{pmatrix}
c_2^{\lambda(2)}\cdots c_a^{\lambda(a)}, 
\]
\begin{comment}
{(this is the coefficient of $s^{-j}$)}
\end{comment}
 here, $\mathcal P(b+j; [2, a])$ is the set of partitions of $b+j$ using only integers in $[2, a]$, 
 and $\lambda(h)$ is the multiplicity of the integer $h$ in the partition $\lambda$. 
Also, 
\[
\begin{pmatrix}
\alpha\\
\beta_1 \; \cdots \; \beta_k
\end{pmatrix}
 = \frac{\prod_{i=0}^{\beta_1+\cdots + \beta_k-1}(\alpha-i)}{\beta_1!\cdots \beta_k!}.
\]
\qed
\end{lem}

The coefficients of the expansions of the other terms in $S^b+S^{b+1}g_0(S)$ can be expressed similarly.
We write it as $S^{b'} = s^{b'}(1+\sum_{i=1}^{\infty}f_i^{(b')}(\bold c)\frac{1}{s^i})$, $b'\geq b$.

\subsection{Comparison of parameterizations}\label{subsec:comparison}
In later sections, we will construct deformations of $\varphi$.
In general, given an $N$-th order deformation $\varphi_N$ of $\varphi$ for some positive integer $N$, 
 there is a non-trivial obstruction to deforming $\varphi_N$ one step further.
It means that the map $\varphi_N$ itself does not deform any more.
To overcome this problem, we will change the values of $c_i$ at each singular point of $\varphi$
 at the orders lower than $t^{N+1}$
 to eliminate the obstruction associated with $\varphi_N$.
This results in changing the map $\varphi_N$ in lower orders.
More precisely, we construct a new map
 $\bar{\varphi}_N$ which is equal to $\varphi_N$ only up to some order $t^{N'}$, $N'<N$,
 but the obstruction to deforming $\bar{\varphi}_N$ vanishes.
Thus, there is a map $\bar{\varphi}_{N+1}$ extending $\bar{\varphi}_N$.
We will do this in a way that $N'$ goes to $\infty$ when $N$ goes to $\infty$.
Thus, eventually we will obtain a formal deformation of $\varphi$.
Then, by applying an appropriate algebraization theorem \cite{Ar}, we will have
 an actual deformation.

In this argument, knowing the relation between the parameterizations 
 on a neighborhood of a singular point of $\varphi$ associated with different values of $c_i$ is important.
In this section, we will study this issue.

At a singular point $p$ of $\varphi$, write 
\[
S^{b} + S^{b+1}g_0(S) = \sum_{l=-\infty}^{\infty} \sigma_{-l}^{}({\bf c})s^l, 
\]
 where $\sigma_{-l}^{}$ is a series, which is the sum of polynomials 
 of the form $f_{b'-l}^{(b')}$ in Lemma \ref{lem:coeff},
 with $b'\geq b$.
Here, $S = s(1 + \sum_{i=1}^{\infty}\prod_{j=0}^{i-1}(\frac{1}{a}-j)\frac{1}{i!}(\sum_{k=2}^{a}\frac{c_k}{s^k})^i)$
 as in Section \ref{subsubsec:obstcocycle}.
As in the previous subsection, we take $c_i\in t^i\Bbb C[[t]]$.
Assume that we have an $N$-th order deformation $\varphi_N$ of $\varphi$.
Let $\{z, w\}$ be a local coordinate system on $X$ around $\varphi(p)$.
We can take it so that the pull back of it
 by $\varphi_N$ is given in the form 
\[
(z, w) = (s^{a} + c_2^{}(N)s^{a-2} + \cdots + c_{a}^{}(N), 
 \sum_{l=0}^{\infty}\sigma_{-l}^{}(c_2^{}(N), \dots, c_{a}^{}(N))s^l + h_N(s, t)),
 \;\; \text{mod $t^{N+1}$},
\]
 where $c_i(N)\in t^i\Bbb C[[t]]$ and $h_N(s, t)$ is a holomorphic function.
From now on, we will write $(c_2^{}(N), \dots, c_{a}^{}(N)) = \bold c(N)$
 for notational simplicity.

We will perturb $\bold c(N)$
 to $\bold c(N+1)$ in the form
\begin{equation}\label{eq:c_i}
c_i^{}(N+1) = c_i^{}(N) + \delta_i,
\end{equation}
 where $\delta_i\in t^{i+1}\Bbb C[[t]]$. 
We compare the parameterizations
\[
(z, w) = (s^{a} + c_2^{}(N)s^{a-2} + \cdots + c_{a}^{}(N), 
 \sum_{l=0}^{\infty}\sigma_{-l}^{}(\bold c(N))s^l)
\]
 and
\[
(z, w) = (s^{a} + c_2^{}(N+1)s^{a-2} + \cdots + c_{a}^{}(N+1), 
 \sum_{l=0}^{\infty}\sigma_{-l}^{}(\bold c(N+1))s^l).
\]
Note that we have dropped the part $h_N(s, t)$.
A calculation needed for this part will be provided later in this subsection.

We first note the following.
\begin{lem}\label{lem:coordchange}
On a {punctured} neighborhood of $p$, there is a change of the parameter $s$ 
 defined over $\Bbb C[[t]]$ which 
 transforms $s^{a} + c_2^{}(N+1)s^{a-2} + \cdots + c_{a}^{}(N+1)$
 into $s^{a} + c_2^{}(N)s^{a-2} + \cdots + c_{a}^{}(N)$.
\end{lem}
\proof
Put
\[
s(N+1) = s - \frac{1}{a}\sum_{i=2}^{a}\frac{\delta_i'}{s^{i-1}} +
  \sum_{i = a+1}^{\infty}\frac{\varepsilon_i}{s^{i-1}},
\]
 where 
 $\delta_i'$ and $\varepsilon_i$ are unknown series in $\Bbb C[[t]]$.
The equation 
\[
s(N+1)^{a} + c_2^{}(N+1)s(N+1)^{a-2} + \cdots + c_{a}^{}(N+1)
 = s^{a} + c_2^{}(N)s^{a-2} + \cdots + c_{a}^{}(N)
\]
 can be solved order by order with respect to $s$ so that the unknown variables
 $\delta_i'$ and $\varepsilon_i$ are uniquely determined.
Explicitly, in the above equation, the condition that the coefficients of 
 $s^{a-i}$ on the left and the right hand sides coincide for $2\leq i\leq a$,
 is equivalent to the
 equation of the form
\[
\delta_i' = F_i(\delta_2, \dots, \delta_i, \delta_2', \dots, \delta_{i-1}', c_2(N), \dots, c_a(N)),
\]
 where $F_i$ is a polynomial.
This fixes $\delta_i'$ uniquely, since by induction we can assume $\delta_2', \dots, \delta_{i-1}'$ are
 already fixed.
Similarly, for $l<0$, the condition that the coefficient of $s^l$ in the above equation vanishes
 is equivalent to the equation
\[
\varepsilon_{a-l} = F_{a-l}(\delta_2, \dots, \delta_a, \delta_2', \dots, \delta_a',
  \varepsilon_{a+1}, \dots, \varepsilon_{a-l-1}, c_2(N), \dots, c_a(N)), 
\]
 for some polynomial $F_{a-l}$.
Again, this determines $\varepsilon_{a-l}$ uniquely.\qed

\begin{defn}
For an element $\delta\in \Bbb C[[t]]$, let $ord(\delta)$ be the maximal integer $k$ such that $\delta$ is divisible
 by $t^k$.
\end{defn}

\begin{lem}\label{lem:d'}
The coefficient $\delta_i'$ in the proof of Lemma \ref{lem:coordchange}
 is given by 
\[
\delta_i' = \delta_i - \sum_{j=2}^{i-2}\frac{a-j}{a}\delta_{i-j}c_j(N) + O(\delta^2),
\]
 for $4\leq i\leq a$, where $O(\delta^2)$ is the sum of terms which are quadratic or more with respect to 
 $\delta_2, \dots, \delta_a$.
Also, we have $\delta_2' = \delta_2$ and $\delta_3' = \delta_3$.
Moreover, we have 
\[
ord(\delta_i')
 \geq i + \min_{j\in \{2, 3,\dots, i-2, i\}}\{ord(\delta_j)-j\}.
\]
\end{lem}
\proof
This follows from direct calculation.
For the last claim, we note that if a monomial 
 $\prod_{j=2}^a \delta_{j}^{p_j}\prod_{k=2}^a (\delta_{k}')^{q_k}\prod_{l=2}^a c_{l}(N)^{r_l}$
 is contained in the part $O(\delta^2)$, we have
\[
\sum_{j=2}^a jp_j + \sum_{k=2}^a kq_k + \sum_{l=2}^a lr_l = i,
\]
 and at least two of $p_j, q_k$ are not zero.
Then, the claim follows by induction
 using the fact $ord(c_l(N)) \geq l$.\qed\\
 
Note that since we assume {$ord(\delta_{i})\geq i+1$}, we have {$ord(\delta_i')\geq i+1$}. 
Similarly, we have the following.
\begin{lem}\label{lem:e_i}
For the coefficient $\varepsilon_i$ in the proof of Lemma \ref{lem:coordchange},
 we have
\[
ord(\varepsilon_i)\geq i + \min_{2\leq j \leq a}\{ord(\delta_j) - j\}.
\]
\end{lem}
\proof
As in the proof of Lemma \ref{lem:coordchange}, $\varepsilon_{a-l}$ can be expressed as a polynomial $F_{a-l}$ of 
 $\delta_i, \delta_i', c_i(N)$ and $\varepsilon_{a+1}, \dots, \varepsilon_{a-l-1}$.
Also, if $\prod_{i=2}^a \delta_{i}^{p_i}\prod_{j=2}^a (\delta_{j}')^{q_j}\prod_{k=2}^a c_{k}(N)^{r_k}
 \prod_{m=a+1}^{a-l-1} \varepsilon_{m}^{s_m}$
 is a monomial in $F_{a-l}$, we have
\[
\sum_{i=2}^a ip_i + \sum_{j=2}^a jq_j + \sum_{k=2}^a kr_k + \sum_{m=a+1}^{a-l-1}ms_m = a-l.
\]
We can inductively assume 
 $ord(\varepsilon_i)\geq i + \min_{2\leq j\leq a}\{ord(\delta_j) - j\}$ for $i< a-l$.
Also, at least one of $p_i, q_j, s_m$ is not zero.
From these observations and Lemma \ref{lem:d'}, we obtain the claim. \qed\\
 
By definition of $c(N+1)$, we have the following.
\begin{prop}
We have
\[
\sum_{l=-\infty}^{\infty} \sigma_{-l}^{}(\bold c(N+1))s(N+1)^l
 = \sum_{l=-\infty}^{\infty} \sigma_{-l}^{}(\bold c(N))s^l,
\]
 which holds over $\Bbb C[[t]]$.
Here, $s(N+1)$ is given in Lemma \ref{lem:coordchange}.
\end{prop}
\proof
The parameterization 
\[
(z, w) =  (s^{a} + c_2^{}(N+1)s^{a-2} + \cdots + c_{a}^{}(N+1), 
 \sum_{l=-\infty}^{\infty}\sigma_{-l}^{}(\bold c(N+1))s^l)
\]
 is obtained by substituting 
 $S = s(1 + \sum_{i=1}^{\infty}\prod_{j=0}^{i-1}(\frac{1}{a}-j)\frac{1}{i!}(\sum_{k=2}^{a}\frac{c_k^{}(N+1)}{s^k})^i)$
 to $(z, w) = (S^{a}, S^{b}+S^{b+1}g_0(S))$.
Substituting $s(N+1)$ to $s$, we have
\[
(z, w) =  (s^{a} + c_2^{}(N)s^{a-2} + \cdots + c_{a}^{}(N), 
 \sum_{l=-\infty}^{\infty}\sigma_{-l}^{}(\bold c(N+1))s(N+1)^l)
\]
 according to the definition of $s(N+1)$.
However, looking at $z$, this must be identical to the substitution of 
 $S = s(1 + \sum_{i=1}^{\infty}\prod_{j=0}^{i-1}(\frac{1}{a}-j)\frac{1}{i!}(\sum_{k=2}^{a}\frac{c_k^{}(N)}{s^k})^i)$
  to $(z, w) = (S^{a}, S^{b}+S^{b+1}g_0(S))$.
Then, by comparing $w$, we obtain the claimed identity.\qed\\

Thus, we have the following. 
\begin{cor}\label{cor:pm}
The equality
\[
 \sum_{l=0}^{\infty}\sigma_{-l}^{}(\bold c(N+1))s(N+1)^l
  - \sum_{l= 0}^{\infty} \sigma_{-l}^{}(\bold c(N))s^l\\
  = \sum_{l=-\infty}^{-1} \sigma_{-l}^{}(\bold c(N))s^l
     - \sum_{l=-\infty}^{-1}\sigma_{-l}^{}(\bold c(N+1))s(N+1)^l
\]
 holds.\qed
\end{cor}

We also have the following by direct calculation.
\begin{lem}\label{lem:perturb}
The equality
\[
\begin{array}{l}
\sum_{l=-\infty}^{-1}\sigma_{-l}^{}(\bold c(N+1))s(N+1)^l \\
 = \sum_{l=-\infty}^{-1}\sigma_{-l}^{}(\bold c(N+1))
  (s - \frac{1}{a}\sum_{i=2}^{a}\frac{\delta_i'}{s^{i-1}} + \sum_{i = a+1}^{\infty}\frac{\varepsilon_i}{s^{i-1}})^l\\
   = \sum_{l=-\infty}^{-1}(\sigma_{-l}^{}(\bold c(N+1))
      - \sum_{i=2}^a\frac{l+i}{a}\delta_i'\bar{\sigma}_{-l-i}
      (\bold c(N+1))
      +\sum_{i=a+1}^{\infty}(l+i)\varepsilon_i\bar{\sigma}_{-l-i}^{}(\bold c(N+1)) + \nu_l)s^l
\end{array}
\]
 holds.
Here, 
 $\nu_l$ is the sum of terms which are quadratic or more with respect to $\delta_i$ and $\varepsilon_i$.
Also, $\bar{\sigma}_m = \sigma_m$ for $m>0$ and $0$ otherwise. \qed.
\end{lem}

Finally, we give a calculation for the part $h_N(s, t)$.
Let $g(s, t)$ be any holomorphic function.
Substituting $s(N+1)$ to $s$, we obtain a meromorphic function 
 $g(s(N+1), t)$.
Let $g(s(N+1), t)_{reg}$ be its regular part and 
 $g(s(N+1), t)_{sing} = g(s(N+1), t) - g(s(N+1), t)_{reg}$ be its singular part
 with respect to the expansion using $s$.
\begin{lem}\label{lem:h_N}
There is a holomorphic function $\bar h_N(s, t)$ such that 
\[
\bar h_N(s(N+1), t)_{reg} = h_N(s, t), \;\; \text{mod $t^{N+2}$}.
\]
\end{lem}
\proof
Substituting $s(N+1)$ to $s$ in $h_N(s, t)$, the difference 
\[
H_1(s, t)=h_N(s(N+1), t)_{reg} - h_N(s, t)
\]
 can be divided by $t^{m}$, where $m$ is the minimal order of $\{\delta_i', \varepsilon_i\}$.
This follows from the definition of $s(N+1)$, 
 see also Lemmas \ref{lem:d'} and \ref{lem:e_i}.
Applying the similar process to $h_N - H_1$, 
 we see that the difference
\[
\begin{array}{ll}
H_2(s, t) & = (h_N - H_1)(s(N+1), t)_{reg} - h_N(s, t)\\
 & = (h_N(s(N+1), t)_{reg}-h_N(s, t)) - H_1(s(N+1), t)_{reg}\\
 & = H_1(s, t) - H_1(s(N+1), t)_{reg}
\end{array}
\]
 can be divided by $t^{2m}$.
Similarly, $H_3(s, t) = (h_N-H_1-H_2)(s(N+1), t)_{reg}-h_N(s, t)$  
 can be divided by $t^{3m}$.
Repeating this, if we put $\bar h_N(s, t) = (h_N-H_1-\cdots - H_k)(s, t)$ so that $(k+1)m>N+1$, 
 it satisfies the requirement of the claim. \qed




\subsection{The dual space of obstructions}
Recall that the obstruction to deforming the map $\varphi$ lies in 
 $H^1(C, \bar{\mathcal N}_{{\varphi}})\cong H^0(C, \varphi^*\omega_X(Z))^{\vee}$.
There is a natural map 
\[
i\colon H^0(X, \omega_X)\to H^0(C, \varphi^*\omega_X(Z))
\]
 given by the pullback.
The argument in \cite{N} (see also \cite{B, KS}) shows the following.
\begin{prop}\label{prop:obst}
Let $\varphi_N$ be an $N$-th order deformation of $\varphi$ and 
 let $o_N\in H^1(C, \bar{\mathcal N}_{{\varphi}})$ be the obstruction class to deforming it one step further.
Then, under the pairing between 
 $H^1(C, \bar{\mathcal N}_{{\varphi}})$ and $H^0(C, \varphi^*\omega_X(Z))$, 
 elements in $\mathrm{Im}\, i$ pair with $o_N$ trivially.
In particular, to identify the obstruction class $o_N$, it suffices to 
 consider the pairing between it and elements in $\mathrm{Coker}\, i$.
\end{prop}
\proof
We have a natural inclusion $\bar{\mathcal N}_{\varphi}\to \varphi^*\mathcal N_{\varphi(C)}$.
Therefore, we have a natural map 
\[
\varphi_*\bar{\mathcal N}_{\varphi}\to 
 \varphi_*\varphi^*\mathcal N_{\varphi(C)}\cong \mathcal N_{\varphi(C)}\otimes \varphi_*\mathcal O_C.
\]
By the Leray's spectral sequence, we have an isomorphism
 $H^1(C, \bar{\mathcal N}_{\varphi})\cong H^1(\varphi(C), \varphi_*\bar{\mathcal N}_{\varphi})$.
Thus, we have a natural map 
\[
\pi\colon H^1(C, \bar{\mathcal N}_{\varphi})\to H^1(\varphi(C), \varphi_*\varphi^*{\mathcal N}_{\varphi(C)})
 \cong H^1(\varphi(C), {\mathcal N}_{\varphi(C)}).
\]
The latter isomorphism is due to the natural exact sequence
\[
0\to {\mathcal N}_{\varphi(C)}\to \varphi_*\varphi^*{\mathcal N}_{\varphi(C)} \to \mathcal Q\to 0,
\]
 where $\mathcal Q$ is a torsion sheaf.
Under this map, the obstruction class $o_N$ gives the obstruction class 
 $\bar{o}_N$ to deforming the image $\varphi(C)$.
The dual of it gives a map $H^1(\varphi(C), {\mathcal N}_{\varphi(C)})^{\vee}\to H^1(C, \bar{\mathcal N}_{\varphi})^{\vee} $.
By the Serre duality and the adjunction, we have a natural isomorphism 
\[
H^1(\varphi(C), \mathcal N_{\varphi(C)})^{\vee} \cong H^0(\varphi(C), \omega_X|_{\varphi(C)}).
\]
Then, the composition of natural maps
\[
H^0(X, \omega_X)\to H^0(\varphi(C), \omega_X|_{\varphi(C)})\cong H^1(\varphi(C), \mathcal N_{\varphi(C)})^{\vee}
 \to H^1(C, \bar{\mathcal N}_{\varphi})^{\vee}\cong H^0(C, \varphi^*\omega_X(Z))
\]
 is the map $i$.

The argument in \cite{N} shows that in the natural pairing between 
 $H^1(\varphi(C), \mathcal N_{\varphi(C)})$ and $H^0(\varphi(C), \omega_X|_{\varphi(C)})$, 
 the class of $H^1(\varphi(C), \mathcal N_{\varphi(C)})$ which is the obstruction to 
 deforming $\varphi(C)$ pairs trivially with those classes of 
 $H^0(\varphi(C), \omega_X|_{\varphi(C)})$ coming from $H^0(X, \omega_X)$.
In particular, the classes in $H^0(X, \omega_X)$ annihilates the class $\bar{o}_N$.
It follows that the image of the map $i$ annihilates the class $o_N$.\qed

\begin{defn}\label{def:sr}
We call a map $\varphi$ \emph{semiregular} if the map 
 $i\colon H^0(X, \omega_X)\to H^0(C, \varphi^*\omega_X(Z))$
 induces a surjection onto the subspace $H^0(C, \varphi^*\omega_X)$.
\end{defn}

Classically, the semiregularity was defined for subvarieties \cite{B, KS, S1, S2}.
Namely, in the case of curves on surfaces, a curve $C\subset X$ is semiregular 
 if the embedding $i\colon C\to X$ is semiregular in the above sense.

\begin{example}
If the surface $X$ is Fano or Calabi-Yau, any map $\varphi\colon C\to X$ is semiregular.
\end{example}

In general, we have an exact sequence on the closed subvariety $\varphi(C)$ of $X$
\[
0\to \omega_X|_{\varphi(C)}\to \varphi_*\varphi^*\omega_X\to \mathcal Q\to 0,
\]
 where $\mathcal Q$ is a torsion sheaf defined by this sequence.
Taking the cohomology, we have
\[
\begin{array}{ll}
0\to H^0(\varphi(C), \omega_X|_{\varphi(C)}) &\to H^0(\varphi(C), \varphi_*\varphi^*\omega_X)
 \to H^0(\mathcal Q)\\
  & \to H^1(\varphi(C), \omega_X|_{\varphi(C)})\to H^1(\varphi(C), \varphi_*\varphi^*\omega_X)
 \to 0.
\end{array}
\]
Then, map $\varphi$ is semiregular in the above sense if
\begin{itemize}
\item the curve $\varphi(C)$ is semiregular in the classical sense (that is, the inclusion $\varphi(C)\to X$ is 
 semiregular), and
\item the map $H^0(\varphi(C), \omega_X|_{\varphi(C)}) \to H^0(\varphi(C), \varphi_*\varphi^*\omega_X)$
 is surjective.
\end{itemize}
Taking the dual, we have
\[
0\to {\rm Hom}_{\mathcal O_{\varphi(C)}}(\varphi_*\varphi^*\omega_X, \omega_{\varphi(C)})\to 
 {\rm Hom}_{\mathcal O_{\varphi(C)}}(\omega_X|_{\varphi(C)}, \omega_{\varphi(C)})
 \to  H^0(\mathcal Q)^{\vee}, 
\]
 where $\omega_{\varphi(C)}$ is the dualizing sheaf of $\varphi(C)$.
The map $H^0(\varphi(C), \omega_X|_{\varphi(C)}) \to H^0(\varphi(C), \varphi_*\varphi^*\omega_X)$
 is surjective if and only if the map 
 ${\rm Hom}_{\mathcal O_{\varphi(C)}}(\omega_X|_{\varphi(C)}, \omega_{\varphi(C)})
  \to  H^0(\mathcal Q)^{\vee}$
 is surjective.
Note that we have
\[
{\rm Hom}_{\mathcal O_{\varphi(C)}}(\omega_X|_{\varphi(C)}, \omega_{\varphi(C)})
 \cong {\rm Hom}_{\mathcal O_{\varphi(C)}}(\mathcal O_{\varphi(C)}, \omega_X^{\vee}|_{\varphi(C)}\otimes \omega_{\varphi(C)})
 \cong H^0(\varphi(C), \mathcal N_{\varphi(C)}).
\]
The pairing between $H^0(\mathcal Q)$ and the image of $H^0(\varphi(C), \mathcal N_{\varphi(C)})$
 in $H^0(\mathcal Q)^{\vee}$ is given as follows.
Namely, an element of $H^0(\mathcal Q)$ is represented by a germ $\xi$ of the sheaf $\varphi_*\varphi^*\omega_X$
 (modulo the germs of the sheaf $\omega_X|_{\varphi(C)}$), which is a section of $\omega_X$
 whose coefficient belongs to $\varphi_*\mathcal O_C$.
Given a section $\alpha$ of $H^0(\varphi(C), \mathcal N_{\varphi(C)})$, since we have
 $\omega_X|_{\varphi(C)}\otimes \mathcal N_{\varphi(C)}\cong \omega_{\varphi(C)}$, 
 $\xi$ pairs naturally with $\alpha$ to give a germ of a section of $\omega_{\varphi(C)}$
 with coefficient in $\varphi_*\mathcal O_C$.
Pulling it back to $C$ gives a germ of a meromorphic 1-form on $C$, and its residue is the value of the pairing.

Thus, if $\varphi(C)$ is sufficiently ample and the normal sheaf $\mathcal N_{\varphi(C)}$ has
 plenty of global sections, the map 
 ${\rm Hom}_{\mathcal O_{\varphi(C)}}(\omega_X|_{\varphi(C)}, \omega_{\varphi(C)})
 \to  H^0(\mathcal Q)^{\vee}$ will be surjective.
Also, if $\varphi(C)$ is sufficiently positive, it is semiregular in the classical sense.
Thus, if $\varphi(C)$ is sufficiently positive (compared to the number of singular points), 
 the map $\varphi$ is semiregular in the sense of Definition \ref{def:sr}.


\section{Deformation of singular curves}
Let $\varphi\colon C\to X$ be a map from a complete integral curve to a smooth surface
 as before.
Let $\{p_1,\, \dots, p_e\}$ be the set of points on $C$ where $\varphi$ is singular.
Namely, we assume $C$ is non-singular at $p_i$ and $d\varphi = 0$ there, as in Section \ref{subsec:localobst}.
At each $p_j$, the image of $\varphi$ is parameterized as
 $(z, w) = (s^{a_j}, s^{b_j}+s^{b_j+1}g_0(s))$ for some integers $1<a_j<b_j$.
Let us write $Z = (d\varphi)= \sum_{j=1}^e (a_j-1)p_j$.
From each singular point $p\in \{p_1\, \dots, p_e\}$,
 there is a contribution to the obstruction controlled by the functions
 $\sigma_{-l}(\bold c)$ in the notation of the previous section.
We write the parameter $c_i$ at $p_j$ by $c_i^{(j)}$ and $\bold c$ by $\bold c^{(j)}$ for clarity.
Also, we write the functions $\sigma_{-l}$ defined on a neighborhood of $p_j$ 
 by $\sigma_{-l}^{(j)}$.
We assume that each $c_i^{(j)}$ belongs to $t^{d_ji}\Bbb C[[t]]$, where
 $d_j$ is a positive integer (to be fixed later, see Definition \ref{def:Mb}), and write $c_i^{(j)} =t^{d_ji}\bar c_i^{(j)}$. 
Assume we have constructed an $N$-th order deformation $\varphi_N$ of $\varphi$
 for some non-negative integer $N$. 

Recall that the function $\sigma_{-l}^{(j)}(\bold c^{(j)})$ is a sum of polynomials of the form
 $f_{b'-l}^{(b')}$, $b'\geq b_j$, $l< b'$, in the notation of Section \ref{subsubsec:regularobst}.
Also, recall that we assume $b_j$ is not a multiple of $a_j$.
Among these functions, the ones with $l<0$ are relevant to the calculation of the obstruction.
Substituting $c_i^{(j)} = t^{d_ji}\bar c_i^{(j)}$ to $f_{b_j-l}^{(b_j)}$, we have
\[
f_{b_j-l}^{(b_j)}(\bold c^{(j)}) 
 = t^{d_j(b_j-l)}f_{b_j-l}^{(b_j)}(\bar{\bold c}^{(j)}).
\]
Properties of the functions $f_{b_j-l}^{(b_j)}$ are crucial to the study of obstructions.

\subsection{Functions $F_{-n}^{(j)}$}\label{subsec:F}
We introduce some functions related to $f_{b_j-l}^{(b_j)}$ for later purposes, see Definition \ref{def:star}.
Recall that around a singular point $p_j$ of $\varphi$, we introduced a function 
 $S = s_j(1 + \sum_{i=1}^{\infty}\prod_{l=0}^{i-1}(\frac{1}{a_j}-l)\frac{1}{i!}(\sum_{k=2}^{a_j}\frac{c_k}{s_j^k})^i)$
 by solving $S^{a_j} = s_j^{a_j}  + \sum_{i=0}^{a_j-2}c_{a_j-i}s_j^i$.
We write it as 
\[
S = s_j(1+\sum_{i=-\infty}^{-1}\gamma_{-i}^{(j)}s_j^i).
\]
Note that we have 
\[
\gamma_{-i}^{(j)} = f^{(1)}_{-i}(\bold c^{(j)}) = 
 \sum_{\lambda\in\mathcal P(-i; [2, a_j])} 
\begin{pmatrix}
\frac{1}{a_j} \\
\lambda(2) \; \cdots \; \lambda(a_j)
\end{pmatrix}
c_2^{\lambda(2)}\cdots c_{a_j}^{\lambda(a_j)}
\]
 in the notation of Section \ref{subsubsec:regularobst}.
We can solve this and express $s_j$ in terms of $S$,
\[
s_j = S(1+\sum_{i=-\infty}^{-1}\theta_{-i}^{(j)}S^i).
\]
Furthermore, for a negative integer $l$, we write
\[
s_j^l = S^l\sum_{i = -\infty}^0\Theta_{-i}^{(j; l)}S^i.
\]
Using this notation, we introduce the following functions.
\begin{defn}\label{def:F}
For a positive integer $n$, we define the function $F_{-n}^{(j)}$ by
\[
F_{-n}^{(j)} = \sum_{i = -n}^{-1} \Theta_{i+n}^{(j; i)} f_{b_j-i}^{(b_j)}.
\]
\end{defn}
The function $\Theta_{-i}^{(j;l)}$ is given by 
\[
\Theta_{-i}^{(j;l)} = 
 \sum_{\lambda\in \mathcal P(-i;[2, \infty))}
 \begin{pmatrix}
 l\\
 \lambda(2) \; \lambda(3) \; \cdots \; 
 \end{pmatrix}
 \prod_{k=2}^{\infty}(\theta_{k}^{(j)})^{\lambda(k)}.
\]
On the other hand, $\theta_k$ is determined by the condition 
\[
\begin{array}{ll}
S &= S(1+\sum_{i=-\infty}^{-1}\theta_{-i}^{(j)}S^i)(1+\sum_{i=-\infty}^{-1}\gamma_{-i}^{(j)}
 (S(1+\sum_{n=-\infty}^{-1}\theta_{-n}^{(j)}S^n))^i)\\
 & = S + \sum_{i=-\infty}^{-1}\theta_{-i}^{(j)}S^{i+1}
   + \sum_{i=-\infty}^{-1}\gamma_{-i}^{(j)}S^{i+1}(1+\sum_{n=-\infty}^{-1}\theta_{-n}^{(j)}S^n)^{i+1}.
\end{array}
\]
Thus, we have
\begin{equation}\label{eq:theta}
\sum_{i=-\infty}^{-1}\theta_{-i}^{(j)}S^{i+1}
  =-\sum_{i=-\infty}^{-1}\gamma_{-i}^{(j)}S^{i+1}(1+\sum_{n=-\infty}^{-1}\theta_{-n}^{(j)}S^n)^{i+1}.
%
\end{equation}

Comparing the coefficients of $S^{i+1}$ on the left and right hand sides, we see
\[
\theta_{-i}^{(j)}
 = -\gamma_{-i}^{(j)} 
 - \sum_{k = i+2}^{-2} \gamma_{-k}^{(j)}(\sum_{n = i-k}^{-1} \frac{(k+1)\cdot k\cdot\cdots\cdot (k+n+2)}{(-n)!})
 \sum_{\lambda\in\mathcal P_{-n}(k-i;[2, \infty))}
 \begin{pmatrix}
  -n\\
  \lambda(2) \; \lambda(3) \; \cdots \; 
  \end{pmatrix}
  \prod_{l=2}^{\infty}(\theta_{l}^{(j)})^{\lambda(l)},
\]
 where $\mathcal P_{-n}(k-i;[2, \infty))$ is the set of partitions of $k-i$ of length $-n$, 
 using integers larger than one.

Note that $\gamma_{-i}^{(j)}$ is weighted homogeneous of degree $-i$ in the sense of 
 Definition \ref{def:homog} below.
Also, we have $\gamma_1^{(j)} = 0$.
We can recursively solve Eq.(\ref{eq:theta}), 
 and it is easy to see that we can write
\[
\theta_{-i}^{(j)} = -\gamma_{-i}^{(j)} + O(\gamma^2), 
\]
 where $O(\gamma^2)$ is the sum of monomials of $\gamma_{-k}^{(j)}$ which is quadratic or more, 
 and $\theta_{-i}^{(j)}$ is also weighted homogeneous of degree $-i$.
It follows that $\Theta_{-i}^{(j;l)}$ is also weighted homogeneous of degree $-i$.
Therefore, the function $F_{-n}^{(j)}$ is  weighted homogeneous of degree $b_j+n$.

\subsection{The condition $({\rm T})$}\label{subsec:TD}
In this subsection, we introduce the condition $({\rm T})$, which ensures 
 a certain transversality property of the set of polynomials $\{f_{b+i}^{(b)}\}$ appeared in
 Lemma \ref{lem:coeff}.
Our final goal is to prove the existence of deformations of $\varphi$ when
 these conditions
 are met at the singular points $p_1, \dots, p_e$.

\begin{defn}\label{def:perturb}
Let $2<a<b$ be integers where $b$ is not a multiple of $a$.
We say that the polynomials $f_{b+1}^{(b)}, \dots, f_{b+a-1}^{(b)}$ satisfy the condition (${\rm T}$)
 at a point $\tilde{\bold c} = (\tilde c_2, \dots, \tilde c_a)\in\Bbb C^{a-1}$
 if the hypersurfaces in $\Bbb C^{a-1}$ defined by 
\[
\bar f_{b+j}^{(b)}(\bold c) = f_{b+j}^{(b)}(\tilde{\bold c}),\;\; j\in \{1, \dots, a-1\}
\]
 have a transversal intersection at $\tilde{\bold c}$.
Here, 
\[
 \bar f_{b+j}^{(b)}(\bold c) 
 = 
  f_{b+j}^{(b)}(\bold c)
   + \sum_{k=2}^{j-1}\frac{(j-k)(c_{k}-\tilde c_{k})}{a}
       f_{b+j-k}^{(b)}(\tilde{\bold c})
  - \sum_{k=2}^{j-1}\frac{j-k}{a}\sum_{l=2}^{k-2}
    \frac{a-l}{a}(c_{k-l}-\tilde c_{k-l})
    \tilde c_lf_{b+j-k}^{(b)}(\tilde{\bold c}).
\]
When $a = 2$, we say $f_{b+1}{(b)}$ satisfies the condition $\rm{(T)}$ at any $\tilde c\in \Bbb C^{\times}$
 by definition.
 \end{defn}
The additional terms $\sum_{k=2}^{j-1}\frac{(j-k)(c_{k}-\tilde c_{k})}{a}
       f_{b+j-k}^{(b)}(\tilde{\bold c})
  - \sum_{k=2}^{j-1}\frac{j-k}{a}\sum_{l=2}^{k-2}
    \frac{a-l}{a}(c_{k-l}-\tilde c_{k-l})
    \tilde c_lf_{b+j-k}^{(b)}(\tilde{\bold c})$
 reflect
 the calculation in Section \ref{subsec:comparison}, see Lemma \ref{lem:leading}.

Note that $f_{b+j}^{(b)}$ is the sum of the monomials of the 
 lowest 
 degree of $\sigma_j^{(b)}$ in view of the following Definition \ref{def:homog}.

\begin{defn}\label{def:homog}
We call a holomorphic function $f$ of ${\bold c}\in\Bbb C^{a-1}$ satisfying 
 $f(\alpha^2 c_2, \dots, \alpha^a c_a) = \alpha^df({\bold c})$
 weighted homogeneous of the degree $d$, where $\alpha$ is any constant
 and $d$ is a non-negative integer.
In particular, the function $f^{(b)}_{b+j}$ is weighted homogeneous of the degree $b+j$.
\end{defn}

In Proposition \ref{prop:perturb} below, the variables $c_2, \dots, c_a$ of $f_{b+j}^{(b)}$ will take values
 in $\Bbb C[[t]]$ as before.
More precisely, we will take $c_i\in t^{di}\Bbb C[[t]]$ for a fixed positive integer $d$.
Accordingly, 
 in the above definition of $\bar f_{b+j}^{(b)}$, we may replace the constant $\tilde c_i$ by 
 some $c_i(-\infty) \in t^{di}\Bbb C[[t]]$ satisfying $c_i(-\infty) = t^{di}\tilde c_i$ mod $t^{di+1}$.
In this case, $\bar f_{b+j}^{(b)}$ becomes
\begin{equation}\label{eq:fbar}
\begin{array}{l}
\bar f_{b+j}^{(b)}(\bold c) 
 = f_{b+j}^{(b)}(\bold c) + \sum_{k=2}^{j-1}\frac{(j-k)(c_{k}-c_{k}(-\infty))}{a}
        f_{b+j-k}^{(b)}({\bold c}(-\infty))\\
  \hspace{1in}   - \sum_{k=2}^{j-1}\frac{j-k}{a}\sum_{l=2}^{k-2}
     \frac{a-l}{a}(c_{k-l}-c_{k-l}(-\infty))
     c_l(-\infty)f_{b+j-k}^{(b)}({\bold c}(-\infty)),
\end{array}
\end{equation}
 where $\bold c(-\infty) = (c_1(-\infty), \dots, c_{a-1}(-\infty))$.
In the construction of deformations of $\varphi$, we need to eliminate obstructions, which essentially involves
 solving a system of polynomial equations.
These equations have solutions, primarily thanks to 
 the transversality property guaranteed by the condition $({\rm T})$.
However, the actual equations contain additional higher order terms, though these terms do not 
 significantly affect the
 properties of the equations.
Proposition \ref{prop:perturb} is formulated in such a way that it can be applied to cases
 with these extra terms.
To simplify the description, we introduce the following notation.

\begin{defn}\label{def:h-}
Fix ${\bold c(-\infty)} = (c_2(-\infty), \dots, c_a(-\infty))\in \Bbb C[[t]]$.
We symbolically write 
 $o_{b+j}(\bold c) = o_{b+j}(c_2, \dots, c_a)\in \Bbb C[c_2, \dots, c_a][[t]]$ 
 for any series which can be expressed as a sum of terms of the following form 
\begin{enumerate}
\item $\alpha{t^lc_2^{l_2}\cdots c_a^{l_a}}$, ${l+d\sum_{p=2}^apl_p \geq d(b+j)+1}$, or 
\item $\alpha{t^lc_2^{l_2}\cdots c_a^{l_a}(c_k-c_k(-\infty))(c_l-c_l(-\infty))}$,
 ${l+d\sum_{p=2}^apl_p \geq d(b+j)-d(k+l)}$,
\end{enumerate}
 where $\alpha$ is a complex number, $l, l_2, \dots, l_a$ are non-negative integers.
The explicit form of $o_{b+j}(\bold c)$ may vary depending on the equations in which they appear. 
\end{defn}

In particular, a polynomial weighted homogeneous of the degree larger than 
 $b+j$ can be a summand of $o_{b+j}$. 
The following calculation is crucial for solving the relevant polynomial equations. 
 
\begin{prop}\label{prop:perturb}
Fix a positive integer $k$.
Assume $f_{b+1}^{(b)}, \dots, f_{b+a-1}^{(b)}$ 
 satisfy the condition $\rm{(T)}$ at $\tilde{\bold c} \in \Bbb C^{a-1}$.
Assume the system of equations 
\[
\bar f_{b+j}^{(b)}(\bold c) = t^{d(b+j)}f_{b+j}^{(b)}(\tilde{\bold c})
                                                  + o_{b+j}(\bold c),
                                                 \;\; \text{\rm{mod} $t^{d(b+j)+k}$},
                                                 \;\; j\in \{1, \dots, a-1\}                                                
\]
 has a solution $\bold c(k-1) = (c_2(k-1), \dots, c_a(k-1))\in (\Bbb C[[t]])^{a-1}$ satisfying 
 $c_i(k-1)=t^{di}\tilde c_i$ mod $t^{di+1}$, $i = 1, \dots, a-1$.
Here, the functions $\bar f_{b+j}^{(b)}(\bold c)$ are defined by Eq.(\ref{eq:fbar}) with respect to 
 some fixed $(c_2(-\infty), \dots, c_a(-\infty))$
 satisfying $c_i(-\infty) = t^{di}\tilde c_i$ mod $t^{di+1}$.
Then, the same system of equations with $k$ replaced by $k+1$ has a solution $\bold c(k)$
 which extends the given solution in the sense that  
 $c_i(k)-c_i(k-1) = 0$ mod $t^{di + k}$
 holds. 
\end{prop}
\proof 
In the case $a = 2$, it is easy to see that we have
 $\bar f_{b+1}^{(b)} = f_{b+1}^{(b)} = c_2^{\frac{b+1}2}$.
In this case, the proof of the claim is easy and will be omitted.

So, assume we have $a>2$.
First, we will perturb $c_2(k-1), \dots, c_a(k-1)$ so that the equation 
\[
\bar f_{b+1}^{(b)}(\bold c) = t^{d(b+1)}f_{b+1}^{(b)}(\tilde{\bold c}) + o_{b+1}(\bold c),
 \;\; \text{mod $t^{d(b+1)+k+1}$}
\]
 holds.
Let $\tilde h_{d(b+1)+k}$ be the coefficient of $t^{d(b+1)+k}$ in $o_{b+1}(\bold c(k-1))-\bar f_{b+1}^{(b)}(\bold c(k-1))$.
By the condition $\rm{(T)}$, we can take a complex vector 
\[
\tilde{\bold c}_{1} = (c_{2, 1}\dots, c_{a, 1})\in 
      \cap_{j\in \{2, \dots, a-1\}}\ker d\bar f_{b+j}^{(b)}(\tilde{\bold c})\subset\Bbb C^{a-1}
\]
 such that 
\[
\bar f_{b+1}^{(b)}(\tilde{\bold c} + \varepsilon \tilde{\bold c}_{1}) 
 = f_{b+1}^{(b)}(\tilde{\bold c}) + \varepsilon \tilde h_{d(b+1)+k} + O(\varepsilon^2)
\] 
 holds for any small positive real number $\varepsilon$.
Then, we have 
\[
\begin{array}{l}
\bar f_{b+1}^{(b)}(t^{2d}(\bar c_2(k-1) + t^{k}c_{2, 1}), \dots, t^{ad}(\bar c_a(k-1) + t^{k}c_{a, 1}))\\
 =t^{d(b+1)}f_{b+1}^{(b)}(\tilde{\bold c}) +
  o_{b+1}(t^{2d}(\bar c_2(k-1) + t^{k}c_{2, 1}), \dots, t^{ad}(\bar c_a(k-1) + t^{k}c_{a, 1}))
  \;\; \text{mod $t^{d(b+1)+k+1}$},
\end{array}
\]
 as required.
Here, we write $c_l(k-1) = t^{dl}\bar c_l(k-1)$.
On the other hand, the equalities
\[
\begin{array}{l}
\bar f_{b+j}^{(b)}(t^{2d}(\bar c_2(k-1) + t^{k}c_{2, 1}), \dots, t^{ad}(\bar c_a(k-1) + t^{k}c_{a, 1}))\\
 = t^{d(b+j)}f_{b+j}^{(b)}(\tilde{\bold c}) +
   o_{b+j}(t^{2d}(\bar c_2(k-1) + t^{k}c_{2, 1}), \dots, t^{ad}(\bar c_a(k-1) + t^{k}c_{a, 1})),
 \;\; \text{mod $t^{d(b+j)+k}$, \;\;
 $j\geq 2$},
\end{array}
\]
 still hold by the homogeneity of $f_{b+j}^{(b)}$ and definition of $o_{b+j}$.

To make the equations 
\[
\bar f_{b+j}^{(b)}(\bold c) = 
 t^{d(b+j)}f_{b+j}^{(b)}(\tilde{\bold c})  + o_{b+j}(\bold c),
  \;\; \text{mod $t^{d(b+j)+k+1}$}, \;\; j\geq 2
\]
 hold, 
 we add appropriate vectors 
 $t^k\tilde{\bold c}_j$,
 where $\tilde{\bold c}_j\in \cap_{l\in \{1, \dots, a-1\}\setminus \{j\}}\ker d\bar f_{b+l}^{(b)}(\tilde{\bold c})$, to 
 ${\bold c}(k-1) + t^k\tilde{\bold c}_{1}$ for all $j \in \{2, \dots, a-1\}$, as in the case of $j = 1$ above.
By the condition $\tilde{\bold c}_j\in \cap_{l\in \{1, \dots, a-1\}\setminus \{j\}}\ker d\bar f_{b+l}^{(b)}$,
 adding $t^k\tilde{\bold c}_j$ changes $\bar f_{b+l}^{(b)}$, $l\in \{1, \dots, a-1\}\setminus \{j\}$, 
 only by quadratic or higher terms with respect to $t^{k}\tilde{\bold c}_j$.
It follows that adding $t^k\tilde{\bold c}_j$ changes $\bar f_{b+l}^{(b)}$, $l\in \{1, \dots, a-1\}\setminus \{j\}$, 
 and $o_{b+l}$
 only at the orders higher than $t^{d(b+l)+k}$.
This shows that $\bold c(k-1) + \sum_{j=1}^{a-1}t^{k}\tilde{\bold c}_j$
 solves the given system of equations for $k$ replaced by $k+1$.  \qed

\subsection{Basis of the dual space of obstructions}\label{subsec:obbasis}

Take an element $\eta$ of $H^0(C, \varphi^*\omega_X(Z))$.
Fix a local coordinate system $\{z_j, w_j\}$ around $\varphi(p_j)$ on $X$.
Also, fix a local parameter $s_j$ on $C$ around $p_j$
 as in Section \ref{subsec:surfacecase}.
%
%

By Proposition \ref{prop:pair},
 an obstruction cocycle to deforming $\varphi$ 
 can be represented by a set of
 local meromorphic sections of the sheaf
 $\bar{\mathcal N}_{{\varphi}}$ on a suitable covering of $C$. 
We take a covering $\mathcal U$ as in Section \ref{subsec:localobst}.
In particular, for each $p_j$, there is a unique open subset $U_{p_j}$ containing it.

The natural pairing between such a representative of the obstruction class and a
 section {$\eta\in H^0(C, \varphi^*\omega_X(Z))$} is also given by Proposition \ref{prop:pair}.
Explicitly, recall that the sheaf $\bar{\mathcal N}_{{\varphi}}$
 is isomorphic to $\varphi^*\omega_X^{-1}\otimes\omega_{C}(-Z)$.
Therefore, $\varphi^*\omega_X(Z)$ is isomorphic to $\bar{\mathcal N}_{{\varphi}}^{\vee}\otimes \omega_C$.
Write $\eta$ in the form 
\[
\eta = \varphi^*(dz_j\wedge dw_j)\widetilde{\eta},
\]
 on a neighborhood of $p_j$,
 where $\widetilde{\eta}$ is a local section of $\mathcal O_{C}(Z)$.
Using the notation in Section \ref{subsec:surfacecase}, 
 the fiberwise pairing between $\bar{\mathcal N}_{{\varphi}}$
 and $\bar{\mathcal N}_{{\varphi}}^{\vee}\otimes \omega_{C}$ 
 is explicitly given by
\[
\varphi^*(dz_j\wedge dw_j)\widetilde{\eta}\otimes \xi\partial_{w_j}
 \mapsto \xi\widetilde{\eta}\varphi^*dz_j 
  = a_j\xi\widetilde{\eta}{s_j^{a_j-1}}ds_j,
\]
 where $\xi$ is a meromorphic function on a neighborhood of $p_j$
 representing the obstruction class. 
We write this by $(\eta, \xi\partial_{w_j})$.
\begin{defn}\label{def:order}
For an element $\eta\in H^0(C, \varphi^*\omega_X(Z))$, 
 we define a subset $psupp(\eta)\subset  \{1, \dots, e\}\times \Bbb Z_{> 0}$ by the property that  
 $(j, m)\in \{1, \dots, e\}\times \Bbb Z_{> 0}$
 belongs to $psupp(\eta)$
 if and only if 
\[
Res_{p_j}(\eta, \frac{1}{s_j^{a_j-m}}\partial_{w_j}) \neq 0.
\] 
This is equivalent to the condition that if we expand $\widetilde{\eta}$
 in terms of $s_j$, its coefficient of $\frac{1}{s_j^m}$ is non-zero.
\end{defn}

\begin{rem}
The parameter $s$ on the domain curve $C$ was chosen so that the map $\varphi$ is locally represented in the 
 form $(z, w) = (s^a, \; s^b + s^{b+1}g_0(s))$ as in Section \ref{subsec:surfacecase}.
Here $b>a$ and $a$ does not divide $b$.
When we fix local coordinates $z, w$ on the target $X$, the choice of $s$ with this property is unique.
When we choose another coordinates $z', w'$ on $X$, with the point $z'=w'=0$ corresponding to the image of 
 a singularity of $\varphi$, we will need to reparameterize the curve $C$ to represent the map $\varphi$
 in the above form.
However, if we write the new parameter as $s'$, it is related to the original by
\[
s' = \alpha s + O(s^a),
\]
 where $\alpha$ is a nonzero constant.
This replacement does not affect the condition 
 $Res_{p}(\eta, \frac{1}{s^{a-m}}\partial_{w}) \neq 0$,
 so the definition of $psupp(\eta)$ does not depend on the choice of the coordinates,
 as long as we choose them so that $(z, w) = (s^a, \; s^b + s^{b+1}g_0(s))$ holds, as we always do in this paper.
\end{rem}

The residue of the meromorphic 1-form $(\eta, \xi\partial_{w_j})$ is the local contribution at the point $p_j$
 to the obstruction paired with $\eta$.
In our case, the coefficient of $s_j^{-m}$ of the section $\xi$ has the form
 $F_{-m}^{(j)}$ as defined in Definition \ref{def:F},
 with some modification in higher order terms, see Section \ref{subsec:calculation of o_{N+1}}.

We prove our main theorem (Theorem \ref{thm:main} below) by explicitly constructing a deformation of $\varphi$.
First, we construct a basis of the space $H^0(C, \varphi^*\omega_X(Z))$
 suitable for our calculations.
We introduce the integers $M$ and $d_j$, $j = 1, \dots, e$, as follows.
\begin{defn}\label{def:Mb}
Let $M$ be the least common multiple of $b_1+1, \dots, b_e+1$.
Also, at each $p_j$, we take $d_j = \frac{M}{b_j+1}$.
Note that with this definition, we have
\[
f_{b_j+1}^{(b_j)}(t^{2d_j}\tilde c_2, \dots, t^{a_jd_j}\tilde c_{a_j})
 = t^Mf_{b_j+1}^{(b_j)}(\tilde{\bold c}).
\]
 \end{defn} 
\begin{defn}\label{def:eZorder}
We introduce a total order to the set $\{1, \dots, e\}\times \Bbb Z_{> 0}$
 by the rule that $(j, m)>(j', m')$ if and only if
\begin{enumerate}
\item $d_j(b_j + a_j - m) < d_{j'}(b_{j'}+a_{j'} - m')$, or
\item $d_j(b_j + a_j - m) = d_{j'}(b_{j'}+a_{j'} - m')$ and $j>j'$.
\end{enumerate}
\end{defn}
In particular, if $j = j'$, we have $(j, m)>(j, m')$ if and only if $m>m'$.

\begin{defn}\label{def:Pord}
For $\eta\in H^0(C, \varphi^*\omega_X(Z))$,
 let $P(\eta) = (j(\eta), m(\eta))$ be the maximal element of $psupp(\eta)$
 with respect to the above order.
Using this notation, let us define $ord(\eta) \in \Bbb Z$ by
\[
ord(\eta) = d_{j(\eta)}(b_{j(\eta)}+a_{j(\eta)} - m(\eta)).
\] 
Note that we have $ord(\eta) = \min_{(j, m)\in psupp(\eta)}\{d_j(b_j+a_j-m)\}$.
We set $ord(\eta) = \infty$ if $\eta$ belongs to $H^0(C, \varphi^*\omega_X)$.
\end{defn}

For a positive integer $N$, let us define the subspace $V_N$ of $H^0(C, \varphi^*\omega_X(Z))$
 by
\[
V_N = \{\eta\in H^0(C, \varphi^*\omega_X(Z))\;|\;
 ord(\eta)\geq N\}.
\]
Also, define $V_{\infty} = H^0(C, \varphi^*\omega_X)$.
They {compose} a sequence of subspaces
\[
V_{\infty}\subset \cdots \subset V_{N+1}\subset V_N\subset \cdots 
\]
 of $H^0(C, \varphi^*\omega_X(Z))$.
Let 
\[
V_{\infty}\subset V_{i_k}\subset \cdots V_{i_{2}}\subset V_{i_1} = H^0(C, \varphi^*\omega_X(Z))
\]
 be the maximal strictly increasing subsequence.
That is, we have $V_{\infty}=V_{i_k+1}\neq V_{i_k}$ and 
 $V_{i_{j+1}} = V_{i_{j+1}-1} = \cdots = V_{i_{j}+1} \neq V_{i_{j}}$, 
 for $j = 1, \dots, k-1$.

 
We have a refinement of the above sequence
\[
V_{i_{j+1}} \subset V_{i_{j}, 1}\subset V_{i_{j}, 2}\subset \cdots\subset V_{i_{j}, e} = V_{i_{j}},
\]
 where
\[
V_{i_{j}, n} 
 = \{\eta\in H^0(C, \varphi^*\omega_X(Z))\;|\; \text{$ord(\eta) \geq i_{j}$, and $j(\eta)\leq n$ if $ord(\eta) = i_j$}\},
\]
 using the notation $P(\eta) = (j(\eta), m(\eta))$.
Let 
\[
V_{i_{j+1}} \subset V_{i_{j}, n_1}\subset V_{i_{j}, n_2}\subset \cdots\subset V_{i_{j}, n_{u_j}} = V_{i_{j}}
\]
 be the subsequence such that 
\[
\dim V_{i_{j}, n_{r+1}} = \dim V_{i_{j}, n_r} + 1,\;\; r = 0, \dots, u_j-1,
\]
 where we define $V_{i_{j}, n_0} = V_{i_{j+1}}$, and $u_j = \dim V_{i_{j}}-\dim V_{i_{j+1}}$.

For $r = 1, \dots, u_j$, 
 let $\eta_{r}^{(i_{j})}\in V_{i_{j}, n_r}\setminus V_{i_{j}, n_{r-1}}$ be any vector. 
Let $W_{i_j}$ be the subspace of $H^0(C, \varphi^*\omega_X(Z))$
 spanned by $\{\eta_1^{(i_j)}, \dots, \eta_{u_j}^{(i_j)}\}$.
Then, the following is obvious.
\begin{lem}
There is a direct sum decomposition
\[
H^0(C, \varphi^*\omega_X(Z)) = H^0(C, \varphi^*\omega_X)\oplus W_{i_1}\oplus\cdots\oplus W_{i_k}. 
\]
Also, the set $\{\eta_1^{(i_j)}, \dots, \eta_{u_j}^{(i_j)}\}$ is a basis of $W_{i_j}$. \qed
\end{lem}

By Proposition \ref{prop:obst}, we have the following.
\begin{cor}\label{cor:dualsection}
Assume that we have an $N$-th order deformation $\varphi_N$ of $\varphi$
 for some non-negative integer $N$.
Then, to prove the existence of deformations of $\varphi_N$ over $\Bbb C[t]/t^{N+2}$, 
 it suffices to show that the pairings between the obstruction class and elements in 
 $\{\eta_1^{(i_j)}, \dots, \eta_{u_j}^{(i_j)}\}$, $j = 1, \dots, k$,
 vanish.\qed
\end{cor}

\begin{defn}\label{def:obstsections}
We write $\mathcal I = \cup_{j=1}^k\{\eta_1^{(i_j)}, \dots, \eta_{u_j}^{(i_j)}\}$.
\end{defn}

\subsection{Main theorem}\label{subsec:mainthm}
We use the same notation as in the previous subsection.
Also, let $\{(l_1, m_1), \dots, (l_{v}, m_{v})\}$ be the subset of $\{1, \dots, e\}\times \Bbb Z_{>0}$
 consisting of the elements satisfying
 $d_{l_q}(b_{l_q} + a_{l_q} - m_{l_q}) = i_{j}$, $q = 1, \dots, v$.
Let us take $\eta\in \{\eta_1^{(i_j)}, \dots, \eta_{u_j}^{(i_j)}\}$.
We write
\[
\{(l_1', m_1'), \dots, (l_w', m_w')\} = psupp(\eta)\cap \{(l_1, m_1), \dots, (l_{v}, m_{v})\}.
\]
Recall that we introduced a function $F_{-n}^{(j)}$ in Section \ref{subsec:F}.

\begin{defn}\label{def:star}
For $\eta\in \{\eta_1^{(i_j)}, \dots, \eta_{u_j}^{(i_j)}\}$, 
 define the equation $(\star_{\eta})$ on $\prod_{j=1}^e \Bbb C^{a_j-1}$ by
\begin{equation*}\label{eq:leading}
(\star_{\eta})\;\;\;\; \sum_{q=1}^w
  Res_{p_{l_q'}}(\eta, F_{-(a_{l_q'}-m_q')}^{(l_q')}(\tilde{\bold c}^{({l_q'})})s^{-(a_{l_q'}-m_q')}\partial_{w_{l_q'}}) = 0,
\end{equation*}
 where $\tilde{\bold c}^{({l_q'})}\in \Bbb C^{a_{l_q'}-1}$.
We will be dealing with those constants of the form 
 $\bold c^{(j)} = (c_2^{(j)}, \dots, c_{a_j}^{(j)})\in (\Bbb C[[t]])^{a_j-1}$, where 
 $c_i^{(j)} = t^{d_ji}\tilde c_i^{(j)} + o(t^{d_ji+1})$, $\tilde c_i^{(j)}\in \Bbb C$.
We say that $\{\bold c^{(j)}\}_{j = 1, \dots, e}$ satisfies the equation $(\star_{\eta})$ if
 $\{\tilde c_i^{(j)}\}_{j = 1, \dots, e; i = 2, \dots, a_j}\in \prod_{j=1}^e \Bbb C^{a_j-1}$ does.
\end{defn}

The following is the main theorem of this paper.
\begin{thm}\label{thm:main}
Assume $\varphi$ is semiregular in the sense of Definition \ref{def:sr}.
If there is a point $\tilde{\bold c}^{(j)}\in\Bbb C^{a_j-1}$ at each $p_j\in \{p_1\, \dots, p_e\}$ satisfying 
 the condition $\rm{(T)}$ for $j = 1, \dots, e$, 
 and also the equations $(\star_{\eta})$ for all $\eta\in\mathcal I$, 
 then there is a deformation of $\varphi$ which deforms the singularity of $\varphi$ at each $p_j$ non-trivially.
\end{thm}
The equations $(\star_{\eta})$ can usually be replaced by much simpler conditions, 
 see Section \ref{sec:examples}.
Note that if $\bold{\tilde c}$ satisfies ${\rm (T)}$, it is not a zero vector, 
 because if $\tilde{\bold c} = 0$, we have $\bar f_{b+j}^{(b)} = f_{b+j}^{(b)}$, and 
 these polynomials are singular at $\tilde{\bold c} = 0$.

\section{Proof of the main theorem}\label{sec:4}

\subsection{Deformation up to the order $t^{M-1}$}\label{subsec:lowerdeformation}
Now, we begin the proof of Theorem \ref{thm:main}.
It is easy to show that the map $\varphi$ has a deformation up to the order $t^{M-1}$, as we will see shortly.
For each $j\in \{1, \dots, e\}$, take 
 $\tilde{\bold c}^{(j)}\in \Bbb C^{a_j-1}$ at which
 the condition of Theorem \ref{thm:main} is satisfied.

Let $M$ be the integer introduced in Definition \ref{def:Mb}.
Recall that up to the order $t^{M-1}$, there is no singular term in the expansion of 
 $S^{b_j} + S^{b_j+1}g_0(S)$ at each $p_j$, see Section \ref{subsubsec:obstcocycle}.
This means that the expression $(z_j, w_j) = (S^{a_j}, S^{b_j} + S^{b_j+1}g_{0}(S))$
 still makes sense,
 and defines a curve whose image is the same as that of $\varphi$.
In particular, we can take local deformations
 such that the difference of them gives a section of the tangent sheaf of $C$, 
 which is zero on the sheaf $\bar{\mathcal N}_{\varphi}$ where the obstruction takes value.
Thus, there is no obstruction to deforming $\varphi$.

\subsection{Leading terms of the obstruction}\label{subsec:calculation of o_{N+1}}
To construct higher order deformations, we need to know and control the obstruction classes.
As we discussed in Section \ref{subsubsec:obstcocycle}, 
 at low orders where the expression $S^b+S^{b+1}g_0(S)$ does not contain singular terms, 
 the obstruction trivially vanishes.
Then, at the order when $S^b+S^{b+1}g_0(S)$
 gives a singular term for the first time, 
 the obstruction is contributed only from the singular points,
 and is calculated as in the proof of Proposition \ref{prop:obstrep}.
At higher orders, there are additional contributions to the obstruction.
In general, finding a representative of the obstruction class using meromorphic sections on 
 open subsets, as discussed in Section \ref{subsec:paring} will be very hard.
The point is that we can still find it for the leading order terms.
In this subsection, we will calculate them.

We fix an open covering $\{U_i\}$ of the domain $C$ of the map $\varphi$ as before.
Let $\{p_1, \dots, p_e\}$ be the set of singular points of the map $\varphi$.
Recall that the curve $C$ is regular at $p_i$, while $C$ may have singular points elsewhere.
Assume we have an $N$-th order deformation $\varphi_N\colon C_N\to X$ of $\varphi$
 for some non-negative integer $N$.
The obstruction class is calculated by 
 the difference of local $(N+1)$-th order deformations.
As we saw in Section \ref{subsec:paring}, such a class can be represented by a
 set of local meromorphic sections $\{\xi_i\}$ 
 of $\bar{\mathcal N}_{\varphi}$ associated with the covering $\{U_i\}$.
We write by $C_N\setminus \{p_1, \dots, p_e\}$ the locally ringed space, which is the restriction of the structure of 
 a locally ringed space on $C_N$ to the topological space underlying $C\setminus\{p_1, \dots, p_e\}$. 
We first show that the restriction of $\varphi_N$ to $C_N\setminus\{p_1, \dots, p_e\}$ is an immersion.
To see this, it suffices to prove the following.
 \begin{lem}\label{lem:isodeform}
Let $C_{N, i}^{\circ}$, $i = 1, 2$, be flat deformations of $C\setminus\{p_1, \dots, p_e\}$ over $\Bbb C[t]/t^{N+1}$. 
Assume that there is a map
 $\tau_N^{\circ}\colon C_{N, 1}^{\circ}\to C_{N, 2}^{\circ}$
  over $\Spec\Bbb C[t]/t^{N+1}$
  which reduces to $id_{C\setminus\{p_1, \dots, p_e\}}$ over $\Bbb C[t]/t$.
 Then, $\tau_N^{\circ}$ is an isomorphism.
 \end{lem}
 \proof
When $N = 0$, the claim is trivial.
So, we assume $N$ is a positive integer.
It suffices to show that for each $q\in C\setminus\{p_1, \dots, p_e\}$, there is an open neighborhood of it
 such that the restriction of $\tau_N^{\circ}$ to it is an isomorphism. 
Note that $q$ may be a singular point of the curve $C$.
Consider an affine open neighborhood $U_q = \Spec R_q$
  of $q$ in $C\setminus\{p_1, \dots, p_e\}$, where $R_q$ is some ring. 
Let $U_{q, N, i} = \Spec \mathcal R_{q, N, i}$, 
 $i = 1, 2$, be the restriction of the structure of locally ringed spaces on $C_{N, i}^{\circ}$
 to the topological space underlying $U_q$.
 The restriction of $\tau_N^{\circ}$ to $U_{q, N, 1}$ gives a map 
 \[
\mathcal R_{q, N, 2}\to \mathcal R_{q, N, 1}.
 \]
We have a commutative diagram 
 \[
 \xymatrix{
 0\ar[r] & R_q\ar[d]^{id}\ar[r] & \mathcal R_{q, N, 2}\ar[r]\ar[d] &  \mathcal R_{q, N-1, 2}\ar[r]\ar[d]  & 0\\
 0\ar[r] & R_q\ar[r] & \mathcal R_{q, N, 1}\ar[r] & \mathcal R_{q, N-1, 1}\ar[r] & 0
 }
 \]
By assumption, 
 when $N = 1$, the map $\mathcal R_{q, 0, 2}\to \mathcal R_{q, 0, 1}$ is the identity map of $R_q$.
Therefore, the map $\mathcal R_{q, 1, 2}\to \mathcal R_{q, 1, 1}$ is also an isomorphism.
 By induction, it follows that the map $\mathcal R_{q, N, 2}\to \mathcal R_{q, N, 1}$ is an isomorphism, too.\qed\\
 
\begin{cor}\label{cor:immersion}
The restriction of the map $\varphi_N$ to $C_N\setminus\{p_1, \dots, p_e\}$ is an immersion.
\end{cor}
\proof
For any point $q\in C\setminus\{p_1, \dots, p_e\}$, there is an open subset 
 $U_q\subset C\setminus\{p_1, \dots, p_e\}$ such that the restriction of $\varphi$ to $U_q$
 is an embedding.
Apply Lemma \ref{lem:isodeform} to the case where $C_{N, 1}^{\circ}$ is the restriction of 
 the structure of a locally ringed space on $C_N$ to $U_q$ (which we write by $U_q$ again), and 
 $C_{N, 2}^{\circ}$ is the image (in the sense of Section \ref{subsec:notation}) of the 
 restriction of $\varphi_N$ to it.
Then, the restriction of $\varphi_N$ to $U_q$ is an isomorphism by Lemma \ref{lem:isodeform}.
The claim follows from this. \qed\\
  
\begin{cor}\label{cor:deformprod}
In the range $N<M$, the curve $C_N\setminus\{p_1, \dots, p_e\}$ is isomorphic to the product 
 $C\setminus\{p_1, \dots, p_e\}\times \Spec\Bbb C[t]/t^{N+1}$.
\end{cor} 
\proof
This follows from the fact that the image of $\varphi_N$ is the same as that of $\varphi$
 for $N<M$.
In particular, 
 we have a map $\psi\colon C_N\setminus\{p_1, \dots, p_e\}\to 
  C\setminus\{p_1, \dots, p_e\}\times\Spec\Bbb C[t]/t^{N+1}$
  which restricts to the identity map over $\Bbb C[t]/t$.
Then, the claim follows from Lemma \ref{lem:isodeform}.
 \qed

\begin{rem}\label{rem:convention}
If $U$ is an open subset of $C$, we will often use the same letter to denote the 
 locally ringed space which is the restriction of $C_N$ to $U$,
 if any confusion would not happen. 
\end{rem}
We observe the following generalization of Proposition \ref{prop:obstrep}.
\begin{prop}\label{prop:2M-1}
In the range where $N<2M-1$ holds, a representative $\{\xi_i\}$ 
 of the obstruction class to deforming $\varphi_N$ can be taken so that $\xi_i = 0$
 unless $U_i = U_{p_j}$ for some $p_j$, where $U_{p_j}$ is the unique open subset in the covering $\{U_i\}$
 containing $p_j$.
\end{prop}
\proof

Let $U_1$ and $U_2$ be open subsets in $C\setminus\{p_1, \dots, p_e\}$.
Let $F_1 = 0$ and $F_2 = 0$ be defining equations, defined over $\Bbb C[t]/t^{N+1}$,
 of the images of
 $\varphi_N|_{U_1}$ and $\varphi_N|_{U_2}$ on suitable open subsets $W_1$ and $W_2$ of $X$,
 respectively.
By taking these open subsets sufficiently small, we can assume
\[
\varphi(U_1\cup U_2)\cap (W_1\cap W_2) = \varphi(U_1\cap U_2)
\]
 holds.
Recall that the image of $\varphi_N$ is the same as that of $\varphi$ over $\Bbb C[t]/t^M$.
That is, we can assume $F_1=0$ and $F_2=0$
 reduce to the defining equations of the images of $\varphi|_{U_1}$ and $\varphi|_{U_2}$
 over $\Bbb C[t]/t^{M}$, respectively. 
In particular, $F_1$ and $F_2$ do not contain terms of the order lower than $M$ with respect to $t$.
On the intersection $U_1\cap U_2$, we have
\[
F_1 = g_{12} F_2 \;\; \text{mod $t^{N+1}$},
\]
 where $g_{12}$ is an invertible holomorphic function on an open subset of $X$.
The function $g_{12}$ does not contain terms of the order lower than $M$ with respect to $t$, either.

One obtains local deformations of $\varphi_N$ on $U_1$ and $U_2$ by 
 considering the equations $F_1 = 0$ and $F_2 = 0$ as defined over $\Bbb C[t]/t^{N+2}$.
Note that, since $\varphi$ is an immersion on $C\setminus\{p_1, \dots, p_e\}$ by Corollary \ref{cor:immersion}, 
 this automatically fixes a local deformation of the domain curve $C_N$ over $U_1$ and $U_2$.
By the above observation, the equality
\[
F_1 = g_{12} F_2 \;\; \text{mod $t^{N+2}$},
\]
 still holds in the range $N<2M-1$.
This implies that there is no contribution to the obstruction cocycle from the difference
 of these local deformations at these orders.



Now, take any local deformation of $\varphi_N$ on the open subset $U_{p_j}$ containing $p_j$.
Let $U_i$ and $U_j$ be open subsets of $C$ that intersect $U_{p_j}$, but do not contain $p_j$.
The difference between local deformations of $\varphi_N$ on $U_{p_j}$ and on such open subsets as $U_i$ and $U_j$
 contributions to the obstruction cocycle.
Recall that the difference between local deformations on $U_{p_j}$ and on $U_i$ 
 gives a section of 
  $\bar{\mathcal N}_{\varphi}|_{U_{p_j}\cap U_i}$.
Then, by the above observation, the restriction of such a section to $U_{p_j}\cap U_i\cap U_j$
 coincides
 with that associated with the open subsets $U_{p_j}$ and $U_j$, in the range $N<2M-1$.
It follows that the difference between local deformations on $U_{p_j}$ and $U_i$
 can be extended to a meromorphic section of $\bar{\mathcal N}_{\varphi}$ on $U_{p_j}$.
On the other hand, assign zero sections to the other open subsets of $\{U_i\}$.
Then, by construction, this represents the obstruction class in the sense of Proposition \ref{prop:pair}.\qed\\ 

Now, we assume $N\geq M-1$.
Let $\xi_j$ be the meromorphic section of $\bar{\mathcal N}_{\varphi}$ on $U_{p_j}$ constructed in 
 the above proposition.
For later purposes, we need to calculate $\{\xi_j\}$ more precisely.
Namely, we will show that under the condition $(\star_{\eta})$, $\eta\in\mathcal I$, in Theorem \ref{thm:main},
 the obstruction contributed by $\{\xi_j\}$ 
 can be absorbed in the term $o_{b+j}(\bold c)$ in the system of equations of Proposition \ref{prop:perturb},
 which we will need to solve in later sections to construct deformations of $\varphi_N$. 
The calculation is rather subtle, since when we compare local deformations at higher orders, 
 full non-linearity of various coordinate changes comes into play. 
We can achieve this by using the special nature of the deformation $\varphi_N$.
Namely, the fact that 
 it has the same image as $\varphi$ up to the order $t^{M-1}$.

We use the same notation as in the proof of Proposition \ref{prop:2M-1}.
The map $\varphi_N$ has a parameterization 
 on $U_{p_j}$ of the form
\begin{equation}\label{eq:0}
\begin{array}{ll}
(z_j, w_j) &=  (s_j^{a_j} + c_2^{(j)}(N)s_j^{a_j-2} + \cdots + c_{a_j}^{(j)}(N), 
 \sum_{l=0}^{\infty}\sigma_{-l}^{(j)}(\bold c^{(j)}(N))s_j^l + t^MH_N(s_j, t))\\
 & =  (S^{a_j}, S^{b_j} + S^{b_j+1}g_0(S) 
  -\sum_{l=-\infty}^{-1}\sigma_{-l}^{(j)}(\bold c^{(j)}(N))s_j^l + t^MH_N(s_j, t)),
  \;\; \text{mod $t^{N+1}$},
\end{array}
\end{equation}
 as in Section \ref{subsec:comparison}.
Here, $S = s_j(1 + \sum_{i=1}^{\infty}
    \prod_{l=0}^{i-1}(\frac{1}{a_j}-l)\frac{1}{i!}(\sum_{k=2}^{a_j}\frac{c_k^{(j)}(N)}{s_j^k})^i)$
    and $H_N$ is a holomorphic function.
Also, $\bold c^{(j)}(N) = (c_2^{(j)}(N), \dots, c_{a_j}^{(j)}(N))$, where $c_i^{(j)}(N)\in t^{d_ji}\Bbb C[[t]]$.
{Considering this as an expression over $\Bbb C[t]/t^{N+2}$}, we obtain a local deformation of $\varphi_N$
 on $U_{p_j}$.
Precisely speaking, we first regard
  $-\sum_{l=-\infty}^{-1}\sigma_{-l}^{(j)}(\bold c^{(j)}(N))s_j^l$ as an expression over $\Bbb C[[t]]$, 
  and reduce it over $\Bbb C[t]/t^{N+2}$.
Then, the term $S^{b_j} + S^{b_j+1}g_0(S) -\sum_{l=-\infty}^{-1}\sigma_{-l}^{(j)}(\bold c^{(j)}(N))s_j^l$
 does not contain a pole up to the order $N+1$ with respect to $t$ after substituting 
 $S = s_j(1 + \sum_{i=1}^{\infty}
     \prod_{l=0}^{i-1}(\frac{1}{a_j}-l)\frac{1}{i!}(\sum_{k=2}^{a_j}\frac{c_k^{(j)}(N)}{s_j^k})^i)$.
On the other hand, 
 we regard $t^MH_N(s_j, t)$ simply as an expression over $\Bbb C[t]/t^{N+2}$.
In other words, 
 it does not contain a term of the order $N+1$ with respect to $t$.

By taking a refinement of the open covering $\{U_i\}$ if necessary, we assume that 
 if $U_i$ contains a singular point of the domain curve $C$,  
 it does not intersect any $U_{p_j}$.
Note that if $U_i$ does not contain a singular point of $C$,
 the restriction of the structure of a locally ringed space on $C_N$ to 
 $U_i$ is isomorphic to the product $U_i\times \Spec\Bbb C[t]/t^{N+1}$.
The same holds for $U_{p_j}$, since it does not contain a singular point of $C$, either.
Then, we have the following.

\begin{lem}\label{lem:atleastM}
On an open subset $U_i$ which intersects $U_{p_j}$ but does not contain $p_j$, we have a local parameterization
 of $\varphi_N$ of the form
\[
(z_{i}, w_{i}) = (f_{i}(s_{i}) + t^Mf_{i, 1}(s_{i}, t), g_{i}(s_{i}) + t^Mg_{i, 1}(s_{i}, t)),
 \;\; \text{mod $t^{N+1}$}.
\]
Here, 
 $f_i, g_i$,  $f_{i, 1}$ and $g_{i, 1}$ are holomorphic functions, 
 $s_i$ is a suitable local parameter on $U_i\times \Spec\Bbb C[t]/t^{N+1}$ (that is, a function which reduces to 
 a parameter on $U_i$ over $\Bbb C[t]/t$),
 and $\{z_i, w_i\}$ is a suitable local coordinate system on $X$.
\end{lem}
\proof
By Corollaries \ref{cor:immersion} and \ref{cor:deformprod},
 the reduction over $\Bbb C[t]/t^M$ of the map $\varphi_N|_{C_N\setminus\{p_1, \dots, p_e\}}$
 is an immersion.
Also, the image of $\varphi_N|_{C_N\setminus\{p_1, \dots, p_e\}}$ is the same as that of 
 $\varphi|_{C\setminus\{p_1, \dots, p_e\}}$
 over $\Bbb C[t]/t^M$.
Thus, by pulling back a fixed parameter on $\varphi(U_i)$, we obtain a
 parameter $s_i$ on $U_i\times \Spec\Bbb C[t]/t^{M}$.
By extending this arbitrarily, we obtain a parameter on $U_i\times \Spec\Bbb C[t]/t^{N+1}$.
Using this parameter, it is clear that the map $\varphi_N$ has a parameterization of the given form.\qed\\

{Regarding this parameterization as defined over $\Bbb C[t]/t^{N+2}$}, 
 we have a local deformation of $\varphi_N$ on $U_i$.
In the range $N < 2M-1$, this is a type of local deformations taken in the proof of Proposition \ref{prop:2M-1}.
Thus, we can use it for the calculation of $\xi_j$.

We now have local deformations of $\varphi_N$ on $U_{p_j}$ and on open subsets $U_i$ which 
 intersect $U_{p_j}$.
Since the obstruction is represented by the difference between them, we need to study their properties.
First, let us clarify the relation between the coordinate functions used in these expressions.
On the curve $C_N$, we have a coordinate change between $s_i$ and $s_j$ defined over $\Bbb C[t]/t^{N+1}$.
When we compare local deformations of $\varphi_N$ on $U_{p_j}$ and $U_i$, 
 the relation between $s_i$ and $s_j$ is given simply by regarding this coordinate change
 over $\Bbb C[t]/t^{N+1}$ as a relation over $\Bbb C[t]/t^{N+2}$, so that $s_i$ is described as a function of 
 $s_j$ which does not contain terms of the order $t^{N+1}$
 (in fact, any other choice which reduces to the given relation over $\Bbb C[t]/t^{N+1}$ will suffice).
Recall that the function $S$ is related to $s_j$ by 
 $S = s_j(1 + \sum_{i=1}^{\infty}
     \prod_{l=0}^{i-1}(\frac{1}{a_j}-l)\frac{1}{i!}(\sum_{k=2}^{a_j}\frac{c_k^{(j)}(N)}{s_j^k})^i)$.
By solving this, we can express $s_j$ in terms of $S$ (see Section \ref{subsec:F}).
Substituting it into $s_j$, we obtain an expression of $s_i$ in terms of $S$, 
 defined over $\Bbb C[t]/t^{N+2}$.

On the other hand, on the target space $X$,
 the functions
 $z_j$ and $w_j$ are expressed as holomorphic functions of $z_i$ and $w_i$ on the intersection of local charts.
Note that 
\[
(z_{i}, w_{i}) = (f_{i}(s_{i}), g_{i}(s_{i})),\;\; \text{mod $t^M$}
\]
 gives a parameterization of the image of the map $\varphi|_{U_i}$.
Also, recall that 
\[
(z_j, w_j) = (S^a, S^b + S^{b+1}g_0(S))
\]
 gives a parameterization of the image of the map $\varphi$
 restricted to $U_{p_j}\setminus\{p_j\}$, the punctured neighborhood of $p_j$.
This holds at any order with respect to $t$.

From these observations, we see the following.
The point of this lemma is that although the relation between $s_i$ and $s_j$ may depend on 
 $t$ at the orders lower than $t^M$ since the domain curve $C$ may deform in general, 
 the relation between $s_i$ and $S$ does not.
\begin{lem}
The relation between $s_i$ and $S$ is given by 
\[
s_i = a(S) + t^Mb(S, t),\;\; \text{mod $t^{N+1}$},
\] 
 where $a(S)$ and $b(S, t)$ are holomorphic functions.
\end{lem}
\proof
By Corollary \ref{cor:immersion}, the map $\varphi_N$ restricted to $C\setminus \{p_1, \dots, p_e\}$
 is an immersion.
Then, 
 by taking $U_i$ small enough, we can assume that the image
 $\varphi_N(U_i)$ is isomorphic to $U_i$.
Note that we use the convention of Remark \ref{rem:convention}.

Then, $s_i$ gives a local coordinate on the image $\varphi_N(U_i)$.
On the other hand, the function $S$ on $C_N$, when reduced over $\Bbb C[t]/t^M$, 
 is the pull back of a fixed parameter on the punctured neighborhood of $\varphi(p_j)$
 on $\varphi(C)$ by the map $\varphi_N$ reduced over $\Bbb C[t]/t^M$.
Since the reduction of 
 $\varphi_N$ over $\Bbb C[t]/t^{M}$,
 when restricted to $C\setminus\{p_1, \dots, p_e\}$, does not depend on $t$ by Corollaries
 \ref{cor:immersion} and \ref{cor:deformprod}, the coordinate change between $s_i$ and $S$ does not contain 
 terms of the order lower than $M$ with respect to $t$.
This proves the claim.\qed\\
 

Now, let us calculate the section $\xi_j$ associated with the open subset $U_{p_j}$
 in Proposition \ref{prop:2M-1}. 
First, let us compute the coordinate transformation of  
\[
(z_{i}, w_{i}) = (f_{i}(s_{i}) + t^Mf_{i, 1}(s_{i}, t), g_{i}(s_{i}) + t^Mg_{i, 1}(s_{i}, t)),
 \;\; \text{mod $t^{N+2}$}, 
\] 
 into $\{z_j, w_j\}$ in terms of $S$.
By the observation so far, it will be of the form
\[
(z_j, w_j) = (S^a + t^MF(S, t), S^b + S^{b+1}g_0(S) + t^MG(S, t)) \;\; \text{mod $t^{N+2}$},
\]
 where $F, G$ are holomorphic functions.
\begin{lem}\label{lem:FGN+1}
In the range $N<2M-1$, the functions $t^MF(S, t)$ and $t^MG(S, t)$ do not contain terms  
 of the order $N+1$ with respect to $t$.
\end{lem}
\proof
Note that, by construction, the expression 
 $(z_{i}, w_{i}) = (f_{i}(s_{i}) + t^Mf_{i, 1}(s_{i}, t), g_{i}(s_{i}) + t^Mg_{i, 1}(s_{i}, t))$,
 mod $t^{N+2}$,
 does not contain terms of the order $N+1$ with respect to $t$.
From this and the relation $s_i = a(S) + t^Mb(S, t)$, and the fact that 
 the coordinate change between $\{z_i, w_i\}$ and $\{z_j, w_j\}$ does not depend on $t$, 
 the claim follows.\qed\\
 
Also, note that the above expression for $(z_j, w_j)$
 coincides with 
 $(z_j, w_j) = (S^a, S^b + S^{b+1}g_0(S) 
   -\sum_{l=-\infty}^{-1}\sigma_{-l}^{(j)}(\bold c^{(j)}(N))s_j^l + t^MH_N(s_j, t))$
 over $\Bbb C[t]/t^{N+1}$, see Eq.(\ref{eq:0}).
Combined with the above lemma, we have the following.
\begin{lem}\label{lem:F}
We have
\[
t^MF(S, t) = 0,\;\; \text{mod $t^{N+2}$},
\]
 in the range $N<2M-1$.\qed
\end{lem}
 
It follows that the local deformations of $\varphi_N$ on $U_{p_j}$ and $U_i$ have the same $z_j$-part.
Therefore, the difference between these local deformations is given by the coefficient of $t^{N+1}$ of 
 the difference of the $w_j$-part, namely, 
\begin{equation}\label{eq:A}
t^MG(S, t) - (-\sum_{l=-\infty}^{-1}\sigma_{-l}^{(j)}(\bold c^{(j)}(N))s_j^l + t^MH_N(s_j, t))
\end{equation}
 after substituting 
 $S = s_j(1 + \sum_{i=1}^{\infty}
    \prod_{l=0}^{i-1}(\frac{1}{a_j}-l)\frac{1}{i!}(\sum_{k=2}^{a_j}\frac{c_k^{(j)}(N)}{s_j^k})^i)$.
Note that it has no term lower than the order $t^{N+1}$ by the remark before Lemma \ref{lem:F}.
Recall that $t^MG(S, t)$ does not contain terms of the order $N+1$ with respect to $t$ before
 substituting $S = s_j(1 + \sum_{i=1}^{\infty}
     \prod_{l=0}^{i-1}(\frac{1}{a_j}-l)\frac{1}{i!}(\sum_{k=2}^{a_j}\frac{c_k^{(j)}(N)}{s_j^k})^i)$,
 while $-\sum_{l=-\infty}^{-1}\sigma_{-l}^{(j)}(\bold c^{(j)}(N))s_j^l + t^MH_N(s_j, t)$ does.
We will evaluate (\ref{eq:A}) using these observations.
Since the calculation is a little tricky, we outline the strategy:
\begin{itemize}
\item Although our ultimate goal is to evaluate (\ref{eq:A}) after substituting 
 $S = s_j(1 + \sum_{i=1}^{\infty}
    \prod_{l=0}^{i-1}(\frac{1}{a_j}-l)\frac{1}{i!}(\sum_{k=2}^{a_j}\frac{c_k^{(j)}(N)}{s_j^k})^i)$
    to $G(S, t)$, we begin by comparing $G(S, t)$ with the result of substituting 
    $s_j = S(1 + \sum_{i = -\infty}^{-1}\theta_{-i}^{(j)}S^i)$ (see below) to 
    $-\sum_{l=-\infty}^{-1}\sigma_{-l}^{(j)}(\bold c^{(j)}(N))s_j^l + t^MH_N(s_j, t)$ (Lemma \ref{lem:Gex}).
\item Then, by studying the result of substituting $S = s_j(1 + \sum_{i=1}^{\infty}
    \prod_{l=0}^{i-1}(\frac{1}{a_j}-l)\frac{1}{i!}(\sum_{k=2}^{a_j}\frac{c_k^{(j)}(N)}{s_j^k})^i)$
    back to the result of Lemma \ref{lem:Gex} carefully (Lemmas \ref{lem:Gorder}, \ref{lem:49}
    and Proposition \ref{prop:G}),
    we obtain the desired result (Corollary \ref{cor:obstexp}). 
\end{itemize}

Now, we begin the study of the term $t^MG(S, t)$ defined over $\Bbb C[t]/t^{N+2}$.
We recall some notations from  Section \ref{subsec:F}.
We can solve
 $S = s_j(1 + \sum_{i=1}^{\infty}
    \prod_{l=0}^{i-1}(\frac{1}{a_j}-l)\frac{1}{i!}(\sum_{k=2}^{a_j}\frac{c_k^{(j)}(N)}{s_j^k})^i)$,
 and express $s_j$ as a Laurent series of $S$.
Note that under the condition $c_k^{(j)}(N)\in t^{d_jk}\Bbb C[[t]]$, if we write
\[
S = s_j(1 + \sum_{i = -\infty}^{-1}\gamma_{-i}^{(j)}s_j^i),
\]
 we have $\gamma_{-i}^{(j)}\in t^{-d_ji}\Bbb C[[t]]$.
Also, note that $\gamma_{-i}^{(j)} = f_{-i}^{(1)}(\bold c^{(j)}(N))$ in the notation of 
 Section \ref{subsubsec:regularobst}.
From this, it is not difficult to see that if we write
\[
s_j = S(1 + \sum_{i = -\infty}^{-1}\theta_{-i}^{(j)}S^i),
\]
 we also have $\theta_{-i}^{(j)}\in t^{-d_ji}\Bbb C[[t]]$.

The following is easy to see.
\begin{lem}\label{lem:Gex}
The term $t^MG(S, t)$ is obtained by substituting $s_j = S(1 + \sum_{i = -\infty}^{-1}\theta_{-i}^{(j)}S^i)$
 to $-\sum_{l=-\infty}^{-1}\sigma_{-l}^{(j)}(\bold c^{(j)}(N))s_j^l + t^MH_N(s_j, t)$
 and discarding all the terms whose order is higher than $N$ with respect to $t$. 
\end{lem}
\proof
Since we have 
 $t^MG(S, t) - (-\sum_{l=-\infty}^{-1}\sigma_{-l}^{(j)}(\bold c^{(j)}(N))s_j^l + t^MH_N(s_j, t)) = 0$,
 mod $t^{N+1}$, 
 the claim is obvious over $\Bbb C[t]/t^{N+1}$.
Then, the claim follows from Lemma \ref{lem:FGN+1}. \qed\\

To compute Eq.(\ref{eq:A}), 
 we need to calculate the term of the order $N+1$ with respect to $t$ when we substitute 
 $S = s_j(1 + \sum_{i = -\infty}^{-1}\gamma_{-i}^{(j)}s_j^i)$ to $t^MG(S, t)$.
By substituting $s_j = S(1 + \sum_{i = -\infty}^{-1}\theta_{-i}^{(j)}S^i)$ to $s_j^l$, $l<0$, 
 we have
\[
s_j^l = S^l(1 + \sum_{i = -\infty}^{-1}\theta_{-i}^{(j)}S^i)^l.
\]
As in Section \ref{subsec:F}, we write this as
\[
s_j^l = S^l\sum_{m=-\infty}^{0}\Theta_{-m}^{(j;l)}S^m,
\]
 where $\Theta_{-m}^{(j;l)}\in t^{-d_jm}\Bbb C[[t]]$.
We write $N-(M-1) = d_jn_j + r$, where $0\leq r< d_j$.
Since we have $M = d_j(b_j+1)$, we can write it as
\[
N +1= d_j(b_j+n_j+1)+r.
\]
For $l<0$, 
 let 
\[
(\sigma_{-l}^{(j)}(\bold c^{(j)}(N))S^l\sum_{m=-\infty}^{0}\Theta_{-m}^{(j;l)}S^m)^{\leq N}
\]
 be the sum of terms of 
 $\sigma_{-l}^{(j)}(\bold c^{(j)}(N))s_j^l = \sigma_{-l}^{(j)}(\bold c^{(j)}(N))S^l\sum_{m=-\infty}^{0}\Theta_{-m}^{(j;l)}S^m$
 whose order with respect to $t$ is at most $N$, in the expression using $S$.
\begin{lem}\label{lem:Gorder}
Any term of $(\sigma_{-l}^{(j)}(\bold c^{(j)}(N))S^l\sum_{m=-\infty}^{0}\Theta_{-m}^{(j;l)}S^m)^{\leq N}$ has
 the order at least $-n_j-1$ with respect to $S$. 
When $N+1$ is a multiple of $d_j$, the bound is given by $-n_j$.
In particular, if we have $-n_j-1>l$, 
 $(\sigma_{-l}^{(j)}(\bold c^{(j)}(N))S^l\sum_{m=-\infty}^{0}\Theta_{-m}^{(j;l)}S^m)^{\leq N}$ is 
 equal to zero.
\end{lem}
\proof 
A term of $S^l\sum_{m=-\infty}^{0}\Theta_{-m}^{(j;l)}S^m$ which is of the order $n\leq l$ with respect to $S$
 has the coefficient in $t^{d_j(l-n)}\Bbb C[[t]]$.
Therefore, if a term of the order $n$ with respect to $S$ is contained in 
 $(\sigma_{-l}^{(j)}(\bold c^{(j)}(N))S^l\sum_{m=-\infty}^{0}\Theta_{-m}^{(j;l)}S^m)^{\leq N}$,
 we have 
\[
d_j(b_j-l)+d_j(l-n)\leq N = d_j(b_j+n_j+1)+r-1.
\]
It follows that the inequality 
\[
n\geq -n_j-1 - \frac{r-1}{d_j}
\]
 holds.
The claim follows from this.\qed

\begin{lem}\label{lem:49}
Assume we have $-n_j-1\leq l$.
If we substitute $S = s_j(1 + \sum_{i = -\infty}^{-1}\gamma_{-i}^{(j)}s_j^i)$ to
 $(\sigma_{-l}^{(j)}(\bold c^{(j)}(N))S^l\sum_{m=-\infty}^{0}\Theta_{-m}^{(j;l)}S^m)^{\leq N}$, 
 then we have 
\[
\sigma_{-l}^{(j)}(\bold c^{(j)}(N))s_j^l + \sum_{m = -n_j-1}^l o_{b_j-m}s_j^m,\;\; \text{mod $t^{N+2}$},
\]
 when
 $N+1$ is not a multiple of $d_j$, and 
\[
\sigma_{-l}^{(j)}(\bold c^{(j)}(N))s_j^l 
 -\sigma_{-l}^{(j)}(\bold c^{(j)}(N))\Theta_{l+n_j+1}^{(j;l)}s_j^{-n_j-1}+
   \sum_{m = -n_j}^l o_{b_j-m}s_j^m,\;\; \text{mod $t^{N+2}$},
\] 
 when $N+1 = d_j(b_j+n_j+1)$.
Note that we have $\Theta_0^{(j;l)} = 1$.
Here, $o_{b_j-m} = o_{b_j-m}(\bold c^{(j)}(N))$ is the notation in Definition \ref{def:h-}.
\end{lem}
\proof
Let $J_l(s_j, t)$ be the series obtained by 
 substituting $S = s_j(1 + \sum_{i = -\infty}^{-1}\gamma_{-i}^{(j)}s_j^i)$ to
 $(\sigma_{-l}^{(j)}(\bold c^{(j)}(N))S^l\sum_{m=-\infty}^{0}\Theta_{-m}^{(j;l)}S^m)^{\leq N}$.
Then, $J_l(s_j, t)$ is the sum of $\sigma_{-l}^{(j)}(\bold c^{(j)}(N))s_j^l$ and terms of the order at least $N+1$
 with respect to $t$.
We write it as 
\[
J_l(s_j, t)=\sigma_{-l}^{(j)}(\bold c^{(j)}(N))s_j^l+\rho(s_j, t).
\]
 
Assume $N+1$ is not a multiple of $d_j$.
Then, we have $J_l(s_j, t) = 0$, mod $t^{N+2}$, if $N+1 < d_j(b_j-l)$.
If we have $N+1 > d_j(b_j-l)$, 
 we can write $J_l(s_j, t)$ in the form
\[
J_l(s_j, t) = \sigma_{-l}^{(j)}(\bold c^{(j)}(N))s_j^l + \sum_{m = -n_j-1}^l o_{b_j-m}s_j^m, \;\; \text{mod $t^{N+2}$}.
\]
Namely, since we have $d_j(b_j+n_j+1) < N+1$, the coefficient of $s_j^m$, $m\geq -n_j-1$,
 is absorbed in $o_{b_j-m}$ when it is of the order $N+1$ with respect to $t$.

Assume $N+1$ is a multiple of $d_j$, that is, $N+1 = d_j(b_j+n_j+1)$.
{If $-n_j-1 = l$}, so that $N+1 = d_j(b_j-l)$, it is clear that we have
 $(\sigma_{-l}^{(j)}(\bold c^{(j)}(N))S^l\sum_{m=-\infty}^{0}\Theta_{-m}^{(j;l)}S^m)^{\leq N} = 0$.
{Assume $-n_j-1<l$}.
Let us write by 
\[
(\sigma_{-l}^{(j)}(\bold c^{(j)}(N))S^l(\sum_{m=-\infty}^{0}\Theta_{-m}^{(j;l)}S^m))^{= N+1}
\]
 the sum of the terms of $(\sigma_{-l}^{(j)}(\bold c^{(j)}(N))S^l(\sum_{m=-\infty}^{0}\Theta_{-m}^{(j;l)}S^m))^{\leq N+1}$
 whose coefficient is of the order ${N+1}$ with respect to $t$.
In this case, if we substitute $S = s_j(1 + \sum_{i = -\infty}^{-1}\gamma_{-i}^{(j)}s_j^i)$ to
 $(\sigma_{-l}^{(j)}(\bold c^{(j)}(N))S^l(\sum_{m=-\infty}^{0}\Theta_{-m}^{(j;l)}S^m))^{\leq N+1}$, 
 the result is $\sigma_{-l}^{(j)}(\bold c^{(j)}(N))s_j^l$, mod $t^{N+2}$.
It follows that the term $\rho(s_j, t)$ in $J_l(s_j, t)$ and 
 the term $(\sigma_{-l}^{(j)}(\bold c^{(j)}(N))S^l(\sum_{m=-\infty}^{0}\Theta_{-m}^{(j;l)}S^m))^{= N+1}$
 cancel, mod $t^{N+2}$, after substituting 
 $S = s_j(1 + \sum_{i = -\infty}^{-1}\gamma_{-i}^{(j)}s_j^i)$ to the latter.
On the other hand, the term 
 $(\sigma_{-l}^{(j)}(\bold c^{(j)}(N))S^l(\sum_{m=-\infty}^{0}\Theta_{-m}^{(j;l)}S^m))^{= N+1}$ is written in the form
\[
\sigma_{-l}^{(j)}(\bold c^{(j)}(N))\Theta_{l+n_j+1}^{(j:l)}S^{-n_j-1}
 + \sum_{m = -n_j}^lo_{b_j-m}S^m, \;\; \text{mod $t^{N+2}$}.
\]
The claim follows from this.\qed\\

Now, we can show the following.
\begin{prop}\label{prop:G}
If $N+1= d_j(b_j+n_j+1)+r$ is not a multiple of $d_j$, 
 after substituting $S = s_j(1 + \sum_{i = -\infty}^{-1}\gamma_{-i}^{(j)}s_j^i)$,
 we can write $t^MG(S, t)$ 
 in the form
\[
-\sum_{m = -n_j-1}^{-1}\sigma_{-m}^{(j)}(\bold c^{(j)}(N))s_j^m +\sum_{m=-\infty}^{-1} o_{b_j-m}s_j^{m} + t^Mh(s_j, t),
\]
 over $\Bbb C[t]/t^{N+2}$.
Here, $h(s_j, t)$ is a holomorphic function which does not contain a negative power of $s_j$.
If $N+1 = d_j(b_j+n_j+1)$, 
 we can write it
 in the form
\[
-\sum_{m = -n_j}^{-1}\sigma_{-m}^{(j)}(\bold c^{(j)}(N))s_j^m +
  \sum_{m=-n_j}^{-1}
  \sigma_{-m}^{(j)}(\bold c^{(j)}(N))\Theta_{m+n_j+1}^{(j;m)}s_j^{-n_j-1}+ \sum_{m=-\infty}^{-1} o_{b_j-m}s_j^{m}
  + t^Mh(s_j, t),
\]
 over $\Bbb C[t]/t^{N+2}$.
\end{prop}
\proof
Recall that $G(S, t)$ is given by 
\[
(-\sum_{l=-\infty}^{-1}\sigma_{-l}^{(j)}(\bold c^{(j)}(N))S^l\sum_{m=-\infty}^{0}\Theta_{-m}^{(j;l)}S^m
 + t^MH_N(S\sum_{m=-\infty}^{0}\Theta_{-m}^{(j;1)}S^m, t))^{\leq N},
\]
 by Lemma \ref{lem:Gex}.
Note that $H_N(s_j, t)$ is a holomorphic function, that is, it does not
 contain a singular term with respect to $s_j$.
Therefore, for $l<0$, the coefficient of $S^l$ in 
 $t^MH_N(S\sum_{m=-\infty}^{0}\Theta_{-m}^{(j;1)}S^m, t)$
 belongs to $t^{M+d_j(-l+1)}\Bbb C[[t]]$.
Since we have $M+d_j(-l+1) = d_j(b_j-l+2) > d_j(b_j-l)$, 
 such a term can be absorbed in $o_{b_j-l}S^l$.
Thus, we can write 
\[
(t^MH_N(S\sum_{m=-\infty}^{0}\Theta_{-m}^{(j;1)}S^m, t))^{\leq N} = \sum_{m=-\infty}^{-1} o_{b_j-m}S^{m}+ 
 t^M\tilde h(S, t),
\]
 for any $N$.
Here, $\tilde h$ does not contain a negative power of $S$.
By the same argument, it is easy to see that after substituting $S = s_j(1 + \sum_{i = -\infty}^{-1}\gamma_{-i}^{(j)}s_j^i)$,
 this is written in the form 
 $\sum_{m=-\infty}^{-1} o_{b_j-m}s_j^{m} + t^Mh(s_j, t)$.

On the other hand, by Lemmas \ref{lem:Gorder} and \ref{lem:49}, we have
\[
(-\sum_{l=-\infty}^{-1}\sigma_{-l}^{(j)}(\bold c^{(j)}(N))S^l\sum_{m=-\infty}^{0}\Theta_{-m}^{(j;l)}S^m)^{\leq N}
 = -\sum_{m = -n_j-1}^{-1}\sigma_{-m}^{(j)}(\bold c^{(j)}(N))s_j^m + \sum_{m=-\infty}^{-1} o_{b_j-m}s_j^{m},\;\; 
  \text{mod $t^{N+2}$},
\] 
 when $N+1$ is not a multiple of $d_j$, and
\[
\begin{array}{l}
(-\sum_{l=-\infty}^{-1}\sigma_{-l}^{(j)}(\bold c^{(j)}(N))S^l\sum_{m=-\infty}^{0}\Theta_{-m}^{(j;l)}S^m)^{\leq N}\\
 = -\sum_{m = -n_j}^{-1}\sigma_{-m}^{(j)}(\bold c^{(j)}(N))s_j^m +
  \sum_{m=-n_j}^{-1}\sigma_{-m}^{(j)}(\bold c^{(j)}(N))\Theta_{m+n_j+1}^{(j;m)}s_j^{-n_j-1}
  + \sum_{m=-\infty}^{-1} o_{b_j-m}s_j^{m},\;\; 
    \text{mod $t^{N+2}$},
\end{array}
\]
 when $N+1 = d_j(b_j+n_j+1)$, 
 after substituting $S = s_j(1 + \sum_{i = -\infty}^{-1}\gamma_{-i}^{(j)}s_j^i)$.
Note that in the latter case, {$\sigma_{n_j+1}^{(j)}(\bold c^{(j)}(N))s_j^{-n_j-1}$ is excluded}.
The claim follows from the observation so far.\qed\\

Recall that we are calculating the difference
 $t^MG(S, t) - (-\sum_{l=-\infty}^{-1}\sigma_{-l}^{(j)}(\bold c^{(j)}(N))s_j^l + t^MH_N(s_j, t))$,
 after substituting $S = s_j(1 + \sum_{i = -\infty}^{-1}\gamma_{-i}^{(j)}s_j^i)$.
By the above proposition, we have the following.
\begin{cor}\label{cor:obstexp}
If $N+1= d_j(b_j+n_j+1)+r$ is not a multiple of $d_j$, we can write
\[
t^MG(S, t) - (-\sum_{l=-\infty}^{-1}\sigma_{-l}^{(j)}(\bold c^{(j)}(N))s_j^l + t^MH_N(s_j, t))
 = \sum_{m=-\infty}^{-1} o_{b_j-m}s_j^{m} + t^Mk(s_j, t),\;\; \text{mod $t^{N+2}$}.
\]
If $N+1= d_j(b_j+n_j+1)$, we can write
\[
\begin{array}{l}
t^MG(S, t) - (-\sum_{l=-\infty}^{-1}\sigma_{-l}^{(j)}(\bold c^{(j)}(N))s_j^l + t^MH_N(s_j, t))\\
 = 
  \sum_{m=-n_j-1}^{-1}\sigma_{-m}^{(j)}(\bold c^{(j)}(N))\Theta_{m+n_j+1}^{(j;m)}s_j^{-n_j-1}
   + \sum_{m=-\infty}^{-1} o_{b_j-m}s_j^{m}
   + t^Mk(s_j, t), \;\; \text{mod $t^{N+2}$}.
\end{array}
\]
Here, $k(s_j, t)$ is a holomorphic function which does not contain a negative power of $s_j$.
\qed
\end{cor} 

By Proposition \ref{prop:2M-1}, the functions in Corollary \ref{cor:obstexp}
 give a representative $\{\xi_j\}$ of the obstruction class.

Recall that the obstruction is calculated by the pairing between the local meromorphic section $\xi_j$
 on the open subset $U_{p_j}$ and elements in 
 $\mathcal I$ of Definition \ref{def:obstsections}
 as in Proposition \ref{prop:pair}, in the range $N<2M-1$.
Now, let us consider these pairings.
Note that the fiberwise pairing between elements of $\mathcal I$ and a meromorphic section 
 of $\bar{\mathcal N}_{\varphi}|_{U_{p_j}}$ gives a meromorphic 1-form on 
 $U_{p_j}$. 
 
Let us recall some notations introduced in Sections \ref{subsec:obbasis}
 and \ref{subsec:mainthm}.
Recall that we introduced the notation $P(\eta) = (j(\eta), m(\eta))$ and
 $ord(\eta) = d_{j(\eta)}(b_{j(\eta)}+a_{j(\eta)}-m(\eta))$ in Definition \ref{def:Pord}.
We write $d_{j(\eta)}(b_{j(\eta)}+a_{j(\eta)}-m(\eta)) = i_n$ for some $n\in \{1, \dots, k\}$.
Let $\{(l_1, m_1), \dots, (l_{v}, m_{v})\}$ be the subset of $\{1, \dots, e\}\times \Bbb Z_{>0}$
 consisting of the elements satisfying
 $d_{l_q}(b_{l_q} + a_{l_q} - m_{q}) = i_{n}$, $q = 1, \dots, v$.
We write $psupp(\eta)\cap \{(l_1, m_1), \dots, (l_v, m_v)\} = \{(l_1', m_1'), \dots, (l_w', m_w')\}$.
Recall that when we have a deformation $\varphi_N$ of $\varphi$ over $\Bbb C[t]/t^{N+1}$, 
 we have constants $\bold c^{(j)}(N) = (c_2^{(j)}(N), \dots, c_{a_j}^{(j)}(N))$
 at each singular point $p_j$ of $\varphi$.
The term $c_i^{(j)}(N)$ is of the form $c_i^{(j)}(N) = t^{d_ji}\bar c_i^{(j)}(N)$, 
 where $\bar c_i^{(j)}(N)\in \tilde c_i^{(j)} + t\Bbb C[[t]]$.
We can consider the condition $(\star_{\eta})$ for the constants $(\tilde c_2^{(j)}, \dots, \tilde c_{a_j}^{(j)})$, 
 $j = 1, \dots, e$, see Definition \ref{def:star}.

\begin{cor}
Assume that the condition $(\star_{\eta})$ holds for any $\eta\in \mathcal I$.
Then, for any $\eta\in \mathcal I$ and $N<2M-1$,
 the pairing between $\eta$ and the obstruction cocycle 
  to deforming $\varphi_N$ in Proposition \ref{prop:2M-1}
  is of the form $o_{b_{j(\eta)}+a_{j(\eta)}-m(\eta)}(\bold c^{(j(\eta))}(N))$. 

\end{cor}
\proof 
According to the definition of $psupp(\eta)$ and $P(\eta)=(j(\eta),m(\eta))$, 
 when we examine the summands of $\sum_{m=-\infty}^{-1}o_{b_{j}-m}(\bold c^{(j)}(N))s_{j}^m$
 in Corollary \ref{cor:obstexp}
 at each $j = 1, \dots, e$,
 only the terms $o_{b_j-m}(\bold c^{(j)}(N))s_j^m$, 
 that satisfy the condition $d_{j}(b_j-m)\geq d_{j(\eta)}(b_{j(\eta)} + a_{j(\eta)} - m(\eta))$, pair with $\eta$
 in a non-trivial manner.
Consequently, the contribution from the part $\sum_{m=-\infty}^{-1}o_{b_{j}-m}(\bold c^{(j)}(N))s_{j}^m$
 of $\xi_j$
 to the pairing between $\eta$ and the obstruction cocycle takes the form 
 $o_{b_{j(\eta)}+a_{j(\eta)}-m(\eta)}(\bold c^{(j(\eta))}(N))$.
Note that the part $t^Mk(s_{j}, t)$ pairs with $\eta$ trivially, for any $j = 1, \dots, e$.
This observation confirms the claim when $N+1$ is not a multiple of $d_{j(\eta)}$.

Now, let us assume that $N+1$ is a multiple of $d_{j(\eta)}$.
If the equality $N+1 = d_{j(\eta)}(b_{j(\eta)} + a_{j(\eta)} -m(\eta))$
 does not hold, it is easy to see that the pairing between $\eta$ and the obstruction class takes the form 
 $o_{b_{j(\eta)}+a_{j(\eta)}-m(\eta)}(\bold c^{(j(\eta))}(N))$ by the construction of $\mathcal I$.
So, let us assume that the equality $N+1 = d_{j(\eta)}(b_{j(\eta)} + a_{j(\eta)} -m(\eta))$ holds.
We use the notation before this corollary.
Thus, we write $d_{j(\eta)}(b_{j(\eta)} + a_{j(\eta)} -m(\eta)) = i_n$ for some $n$,
 and define the pairs $\{(l_1', m_1'), \dots, (l_w', m_w')\}$ as above.
Then, the pairing between $\eta$ and the obstruction cocycle is the sum of terms of the form 
 $o_{b_{j(\eta)}+a_{j(\eta)}-m(\eta)}(\bold c^{(j(\eta))}(N))$, and
 the terms contributed from 
\[
\sum_{m=-(a_{l_r'}-m_{r}')}^{-1}f_{b_{l_r'}-m}^{(b_{l_r'})}(\bold c^{(l_r')}(N))
 \Theta_{m+a_{l_r'}-m_{r}'}^{(l_r';m)}s_{l_r'}^{-(a_{l_r'}-m_r')}
 = F_{-(a_{l_r'}-m_r')}^{(l_r')}(\bold c^{(l_r')}(N))s_{l_r'}^{-(a_{l_r'}-m_r')},\;\; 
 r = 1, \dots, w.
\]
Here, recall that $f_{b_{l_r'}-m}^{(b_{l_r'})}$ is the leading term of $\sigma_{-m}^{({l_r'})}$.
This contribution is nothing but the one that appeared in the condition $(\star_{\eta})$. 
Thus, under the condition $(\star_{\eta})$, the pairing between $\eta$ and the obstruction cocycle is
 of the form $o_{b_{j(\eta)}+a_{j(\eta)}-m(\eta)}(\bold c^{(j(\eta))}(N))$. \qed\\

In the range $N\geq 2M-1$, the obstruction class may not be written in the form 
 as in Proposition \ref{prop:2M-1}.
However,  
 the pairing between any section 
 $\eta\in \mathcal I$
 and the obstruction cocycle is of the order $t^{N+1}$.
Since $N+1\geq 2M > \max_{j=1, \dots, e}\{d_j(b_j+a_j-1)\}$, 
 the pairing is of the form $o_{b_{j(\eta)}+a_{j(\eta)}-m(\eta)}(\bold c^{(j)}(N))$.
Thus, we can conclude the following.
\begin{prop}\label{prop:obsteval}
For any $N\geq M-1$, assume we have constructed a deformation $\varphi_N$ of $\varphi$
 over $\Bbb C[t]/t^{N+1}$.
Also, assume that the condition $(\star_{\eta})$ holds for any $\eta\in \mathcal I$.
Let $\xi_{N+1}$ be an obstruction cocycle to deforming $\varphi_N$ obtained by 
 the difference of local deformations of $\varphi_N$.
Then, we can take $\xi_{N+1}$ so that the pairing between $\eta\in \mathcal I$
 and $\xi_{N+1}$ is of the form $o_{b_{j(\eta)}+a_{j(\eta)}-m(\eta)}(\bold c^{(j(\eta))}(N))$, 
 where $\bold c^{(j)}(N) = (c_2^{(j)}(N), \dots, c_{a_j}^{(j)}(N))$ is given by
 Eq.(\ref{eq:0}) above.\qed 
\end{prop}
Note that $c_i^{(j)}(N) \in t^{d_ji}\Bbb C[[t]]$, and 
 the pairing between $\eta$ and $\xi_{N+1}$ gives an element of 
 $t^{N+1}\Bbb C[t]/t^{N+2}$.
The coefficient of $t^{N+1}$ is the value of the pairing between $\eta$ and
 the obstruction cohomology class $o_{N+1}\in H^1(C, \bar{\mathcal N}_{\varphi})$
 represented by $\xi_{N+1}$.

\subsection{Deformation at higher orders}\label{subsec:step3}

Assume that we have constructed a deformation $\varphi_{N}$ of $\varphi$ up to the order $t^N$
 for some $N\geq M-1$.
In particular, we have fixed the constants $\bold c^{(j)}(N) = (c_2^{(j)}(N), \dots, c_{a_j}^{(j)}(N))$,
 for each $j = 1, \dots, e$.
Here, $\{p_1, \dots, p_e\}$ is the set of singular points of $\varphi$.
We write $c_i^{(j)}(N) = t^{d_ji}\bar c_i^{(j)}$, $\bar c_i^{(j)}\in \tilde c_i^{(j)}+t\Bbb C[[t]]$, 
 and we assume that $(\tilde c_2^{(j)}, \dots, \tilde c_{a_j}^{(j)})$ satisfies the condition of Theorem \ref{thm:main}.
 
Let us consider the deformation at the order $t^{N+1}$.
At this order, there is an obstruction cocycle of the order $t^{N+1}$
 defined by the difference of local deformations. 
Our goal is to make its cohomology class zero by making further modification to $c_{i}^{(j)}(N)$, and
 to construct a deformation of $\varphi$ of the order $t^{N+1}$.
In this process, we have two primary tasks.
The first task is to describe the system of equations whose solution corresponds to 
 the value of $c_2^{(j)}, \dots, c_{a_j}^{(j)}$, at which the obstruction at the order $t^{N+1}$ vanishes.
We need to check that this system is of the form described in Proposition \ref{prop:perturb} in order to 
 obtain a solution.
This will be done in Proposition \ref{prop:cal}.
Second, we need to verify that the solution actually corresponds to some curve.
Namely, to eliminate the obstruction at the order $t^{N+1}$, we have to change the values of
 $c_2^{(j)}, \dots, c_{a_j}^{(j)}$, at lower orders.
Thus, it is unclear that using the new values of $c_2^{(j)}, \dots, c_{a_j}^{(j)}$, we can deform
 $\varphi$ even up to the order $t^N$.
We will show that with the new values of $c_2^{(j)}, \dots, c_{a_j}^{(j)}$,
 we can construct a new deformation $\bar{\varphi}_N$ of $\varphi$ of
 the order $t^N$, and the obstruction to deforming it vanishes.
This is the content of Section \ref{subsec:main}.

\subsubsection{Local and virtual local deformations}\label{subsec:virtdeform}
Locally on a neighborhood $U_{p_j}$ of $p_j$, the map $\varphi_N$ has a parameterization of the form
\begin{equation}\label{eq:8}
(z_j, w_j) = (s_j^{a_j} + c_2^{(j)}(N)s_j^{a_j-2} + \cdots + c_a^{(j)}(N), 
 \sum_{l=0}^{\infty}\sigma_{-l}^{(j)}(\bold c^{(j)}(N))s_j^l + t^MH_j(s_j, t)),\;\; \text{mod $t^{N+1}$},
\end{equation}
 where 
 $H_j(s_j, t)$ is a holomorphic function on a neighborhood of $p_j$.
 
To calculate the obstruction to deforming $\varphi_N$ to the next order, we
 considered specific local deformations of $\varphi_N$ in the previous subsection.
Namely, regarding Eq.(\ref{eq:8}) as an expression over $\Bbb C[t]/t^{N+2}$,
 it gives a local deformation of $\varphi_N$ around $p_j$.
Away from the singular points of $\varphi_N$, we 
 take local deformations as in the proof of Proposition \ref{prop:2M-1} in the range $N<2M-1$.
In the range $N\geq 2M-1$, we take any local deformation.
Let $o_{N+1}\in H^1(C, \bar{\mathcal N}_{\varphi})$ be the obstruction class associated with these local deformations.

We will compare this with a curve given by the parameterization 
\begin{equation}\label{eq:vir}
(z_j, w_j) = (s_j^{a_j} + c_2^{(j)}(N+1)s_j^{a_j-2} + \cdots + c_{a_j}^{(j)}(N+1), 
  \sum_{l=0}^{\infty}\sigma_{-l}^{(j)}(\bold c^{(j)}(N+1))s_j^l + t^M\bar H_N(s_j, t)),
 \;\; \text{mod $t^{N+2}$},
\end{equation}
 around each $p_j$,
 and find values of $c_i^{(j)}(N+1)$ in a way that the map $\varphi$ can be deformed up to 
 the order $t^{N+1}$.
Here, we take 
\[
c_i^{(j)}(N+1) = c_i^{(j)}(N) + \delta_i,
\]
 to be a modification of $c_i^{(j)}(N)$, where $\delta_i\in t^{d_ji+1}\Bbb C[[t]]$ is to be determined.
By Lemma \ref{lem:coordchange}, for any choice of $c_i^{(j)}(N+1)$, 
 there is a change of parameters from $s_j$ to 
 $s_j(N+1)$ on a {punctured} neighborhood of $p_j$ such that 
\[
s_j(N+1)^{a_j} + c_2^{(j)}(N+1)s_j(N+1)^{a_j-2} + \cdots + c_{a_j}^{(j)}(N+1)
 = s_j^{a_j} + c_2^{(j)}(N)s_j^{a_j-2} + \cdots + c_{a_j}^{(j)}(N).
\]
This holds over $\Bbb C[[t]]$.
Moreover, $\bar H_N(s_j, t)$ is a holomorphic function determined by $H_N(s_j, t)$ and $c_i^{(j)}(N+1)$
 by Lemma \ref{lem:h_N}.
In particular, we have {$\bar H_N(s_j(N+1), t)_{reg} = H_N(s_j, t)$}, mod $t^{N+2} $. 
Note that the parameterization Eq.(\ref{eq:vir}) may {not} coincide with the restriction
 of $\varphi_N$ over $\Bbb C[t]/t^{N+1}$.
That is, it is not a local deformation of $\varphi_N$ in general,
 {contrary to the one associated with Eq.(\ref{eq:8}) above.}
Therefore, we call the curve given by the parameterization Eq.(\ref{eq:vir})
 a \emph{virtual local deformation}.

\subsubsection{Comparison between local and virtual local deformations}\label{subsec:virtuallocal}
Substituting $s_j(N+1)$ to $s_j$ on the right hand side of Eq.(\ref{eq:vir}),
 we obtain a parameterization
\[
\begin{array}{l}
(z_j, w_j) \\
 = (s_j^{a_j} + c_2^{(j)}(N)s_j^{a_j-2} + \cdots + c_{a_j}^{(j)}(N), \\
 \hspace{1in}\sum_{l=0}^{\infty}\sigma_{-l}^{(j)}(\bold c^{(j)}(N+1))s_j(N+1)^l
 + t^M\bar H_N(s_j(N+1), t)),\\
  = (s_j^{a_j} + c_2^{(j)}(N)s_j^{a_j-2} + \cdots + c_{a_j}^{(j)}(N), \\
  \hspace{1in} \sum_{l=0}^{\infty}\sigma_{-l}^{(j)}(\bold c^{(j)}(N+1))s_j(N+1)^l 
  + t^M(\bar H_N(s_j(N+1), t)_{sing} + H_N(s_j, t))),
\end{array}
\]
 {mod $t^{N+2}$},
 by Lemma \ref{lem:h_N}.
Now, the difference between local and virtual local deformations is given by 
 the difference of the coordinate $w_j$:
\begin{equation}\label{eq:diamond}
\begin{array}{ll}
\sum_{l=0}^{\infty}\sigma_{-l}^{(j)}(\bold c^{(j)}(N))s_j^l + t^MH_N(s_j, t)\\
 \hspace{.5in} - (\sum_{l=0}^{\infty}\sigma_{-l}^{(j)}(\bold c^{(j)}(N+1))s_j(N+1)^l 
   + t^M(\bar H_N(s_j(N+1), t)_{sing} + H_N(s_j, t)))\\
   =
   - \sum_{l=-\infty}^{-1} \sigma_{-l}^{(j)}(\bold c^{(j)}(N))s_j^l
       + \sum_{l=-\infty}^{-1}\sigma_{-l}^{(j)}(\bold c^{(j)}(N+1))s_j(N+1)^l 
               - t^M\bar H_N(s_j(N+1), t)_{sing},\\
\end{array}
\end{equation}
 mod $t^{N+2}$,
 by  Corollary \ref{cor:pm}.
Note that for fixed $\bold c^{(j)}(N)$, the parameter $s_j(N+1)$ is determined by $c_i^{(j)}(N+1)$.
Thus, the coefficient of $s_j^l$ in 
 $\sum_{l=-\infty}^{-1}\sigma_{-l}^{(j)}(\bold c^{(j)}(N+1))s_j(N+1)^l 
   - t^M\bar H_N(s_j(N+1), t)_{sing}$ is a function of $\bold c^{(j)}(N+1)$, and we write it as 
   $\tilde{\sigma}_{-l}^{(j)}$, that is, 
\[
\sum_{l=-\infty}^{-1}\sigma_{-l}^{(j)}(\bold c^{(j)}(N+1))s_j(N+1)^l 
  - t^M\bar H_N(s_j(N+1), t)_{sing}
 = {\sum_{l=-\infty}^{-1}\tilde{\sigma}_{-l}^{(j)}(\bold c^{(j)}(N+1))s_j^l} 
\]
 over $\Bbb C[t]/t^{N+2}$.
We note the following equality.
\begin{lem}\label{lem:ori}
When $c_i^{(j)}(N) = c_i^{(j)}(N+1)$, we have
\[
\tilde{\sigma}_{-l}^{(j)}(\bold c^{(j)}(N))
 = {\sigma}_{-l}^{(j)}(\bold c^{(j)}(N)).
\]
This holds over $\Bbb C[[t]]$.
\end{lem}
\proof
Recall that when $c_i^{(j)}(N)$ is fixed, 
 both $s_j(N+1)$ and $\bar H_N$ are determined by $c_i^{(j)}(N+1) = c_i^{(j)}(N) + \delta_i$, 
 and when $\delta_i = 0$, we have $s_j(N+1) = s_j$ and $\bar H_N = H_N$.
The claim follows from this.\qed\\

In the argument below, the notation $o_{b_j-l}$
 symbolically refers to functions of the form described in Definition \ref{def:h-}, 
 and their explicit values may vary in different equations.
Recall that $o_{b_j-l}$ is an element of $\Bbb C[c_2^{(j)}, \dots, c_{a_j}^{(j)}][[t]]$.
In the following argument, $c_i^{(j)}$ is regarded as a variable, while $c_i^{(j)}(N)$ is
 a fixed element in $\Bbb C[[t]]$ which plays the role of $c_i(-\infty)$ in Definition \ref{def:h-}.
Based on Lemma \ref{lem:perturb}, we observe the following.
\begin{lem}\label{lem:leading}
In the range $-(a_j-1)\leq l\leq -1$, the term ${\tilde{\sigma}_{-l}^{(j)}(\bold c^{(j)})}
 = \tilde{\sigma}_{-l}^{(j)}(c_2^{(j)}, \dots, c_{a_j}^{(j)})$
 can be written as
\[
{\tilde{\sigma}_{-l}^{(j)}(\bold c^{(j)}) 
 = \bar f_{b_j-l}^{(b_j)}(\bold c^{(j)})
       + o_{b_j-l}(\bold c^{(j)})}
\]
 in the notation of Proposition \ref{prop:perturb}.
Here
 $c_i^{(j)}$ is a variable which takes values in $c_i^{(j)}(N) + t^{d_ji+1}\Bbb C[[t]]$, 
 and 
\[
\begin{array}{l}
{\bar f_{b_j+\alpha}^{(b_j)}(\bold c^{(j)})} 
 = 
 f_{b_j+\alpha}^{(b_j)}(\bold c^{(j)})
  + \sum_{k=2}^{\alpha-1}\frac{(\alpha-k)(c_{k}^{(j)}-c_{k}^{(j)}(N))}{a_j}
      f_{b_j+\alpha-k}^{(b_j)}(\bold c^{(j)}(N))\\
 \hspace{1.5in} - \sum_{k=2}^{\alpha-1}\frac{\alpha-k}{a_j}\sum_{l=2}^{k-2}
   \frac{a_j-l}{a_j}(c_{k-l}^{(j)}-c_{k-l}^{(j)}(N))
   c_l^{(j)}(N)f_{b_j+\alpha-k}^{(b_j)}(\bold c^{(j)}(N)),
 \end{array}
\]
 as in Proposition \ref{prop:perturb}.
\end{lem}
\proof 
By definition, the term $\tilde{\sigma}_{-l}^{(j)}(\bold c^{(j)})$
 has contributions from 
 $\sum_{l=-\infty}^{-1}\sigma_{-l}^{(j)}(\bold c^{(j)})s_j(N+1)^l$
  and from $t^M\bar H_N(s_j(N+1), t)_{sing}$.
First, we study the contribution from 
 $\sum_{l=-\infty}^{-1}\sigma_{-l}^{(j)}(\bold c^{(j)})s_j(N+1)^l$.
 
By Lemma \ref{lem:perturb}, we have
\[
\displaystyle\sum_{l=-\infty}^{-1}\sigma_{-l}^{(j)}(\bold c^{(j)})s_j(N+1)^l \\
\displaystyle = \sum_{l=-\infty}^{-1}(\sigma_{-l}^{(j)}(\bold c^{(j)})-
     \sum_{k=2}^{a_j}\frac{l+k}{a_j}\delta_k'\bar{\sigma}_{-l-k}^{(j)}(\bold c^{(j)})
         + \sum_{k=a_j+1}^{\infty}(l+k)\varepsilon_k\bar{\sigma}_{-l-k}^{(j)}(\bold c^{(j)}) + \nu_{l})s_j^l,
\]
 using the notation there.
Here, we write $c_i^{(j)} - c_i^{(j)}(N) = \delta_i$.
The constants $\delta_i'$ and $\varepsilon_i$
 are determined by $\delta_i$ as in the proof of Lemma \ref{lem:coordchange},
 and these constants in turn determine $s_j(N+1)$.
Moreover, $\nu_{l}$ is the sum of terms which depend on $\delta_i$, 
 $\delta_i'$ and $\varepsilon_i$ quadratically or more,
 and it is easy to see that $\nu_{l}$ can 
 be written in the form $o_{b_j-l}(\bold c^{(j)})$.
Note that in the sum
 $\sum_{k=a_j+1}^{\infty}(l+k)\varepsilon_k\bar{\sigma}_{-l-k}^{(j)}(\bold c^{(j)}(N))$,
 we have $-l-k<0$ for $-(a_j-1)\leq l\leq -1$.
Thus, $\bar{\sigma}_{-l-k}^{(j)} = 0$ by definition and we do not need to deal with this sum.

Note that we have
\[
-\sum_{k=2}^{a_j}\frac{l+k}{a_j}\delta_k'\bar{\sigma}_{-l-k}^{(j)}(\bold c^{(j)})
 = -\sum_{k=2}^{-l-1}\frac{l+k}{a_j}\delta_k'\bar{\sigma}_{-l-k}^{(j)}(\bold c^{(j)}),
\]
 by definition of $\bar{\sigma}_{-l}$.
The term $\frac{l+k}{a_j}\delta_k'\bar{\sigma}_{-l-k}(\bold c^{(j)})$
  can be written in the form 
\[
\begin{array}{l}
\frac{l+k}{a_j}\delta_k'\bar{\sigma}_{-l-k}(\bold c^{(j)})\\
 = \frac{l+k}{a_j}(\delta_{k}-\sum_{m=2}^{k-2}
    \frac{a_j-m}{a_j}\delta_{k-m}c_{m}^{(j)}(N)+ O(\delta^2))\bar{\sigma}_{-l-k}^{(j)}(\bold c^{(j)}) \\
  =  (\frac{(l+k)(c_{k}^{(j)}-c_{k}^{(j)}(N))}{a_j}
   - \frac{l+k}{a_j}\sum_{m=2}^{k-2}
     \frac{a_j-m}{a_j}(c_{k-m}^{(j)}-c_{k-m}^{(j)}(N))
     c_{m}^{(j)}(N)
     + O(\delta^2))\bar{\sigma}_{-l-k}^{(j)}(\bold c^{(j)}).
\end{array}
\]
Here, $O(\delta^2)$ is the sum of terms which depends on $\delta_i$ quadratically or more, 
 and it is easy to see that $O(\delta^2)\bar{\sigma}_{-l-k}^{(j)}(\bold c^{(j)})$
 can be written in the form $o_{b_j-l}(\bold c^{(j)})$.

Also, we have
\[
\bar{\sigma}_{-l-k}(\bold c^{(j)}) = f_{b_j-l-k}^{(b_j)}(\bold c^{(j)}(N)) + o_{b_j-l-k}^{(b_j)}(\bold c^{(j)}(N)),
\]
 for $-l-k>0$.
Thus, we have
\[
\begin{array}{l}
\frac{l+k}{a_j}\delta_k'\bar{\sigma}_{-l-k}(\bold c^{(j)})\\
  =  (\frac{(l+k)(c_{k}^{(j)}-c_{k}^{(j)}(N))}{a_j}
   - \frac{l+k}{a_j}\sum_{m=2}^{k-2}
     \frac{a_j-m}{a_j}(c_{k-m}^{(j)}-c_{k-m}^{(j)}(N))
     c_{m}^{(j)}(N))f_{b_j-l-k}^{(b_j)}(\bold c^{(j)}(N)) + o_{b_j-l}^{(b_j)}(\bold c^{(j)}(N)).
\end{array}
\]

From these observations, we can write
\[
\sum_{l=-\infty}^{-1}(\sigma_{-l}^{(j)}(\bold c^{(j)})-
     \sum_{k=2}^{-l-1}\frac{l+k}{a_j}\delta_k'\bar{\sigma}_{-l-k}^{(j)}(\bold c^{(j)})
         + \nu_l)\\
  = {\bar f_{b_j-l}^{(b_j)}(\bold c^{(j)})}
   + o_{b_j-l}(\bold c^{(j)}).
\] 

Finally, the coefficient of $s_j^l$ in $t^M\bar H_N(s_j(N+1), t)_{sing}$
 has the order at least $M + d_j(-l+1)$ with respect to $t$, and it
 can be absorbed into $o_{b_i-l}$, since $M > d_jb_j$.
This proves the claim.\qed

\subsubsection{Fixing a virtual local deformation}
Let us recall some notations.
At $p_j\in \{p_1, \dots, p_e\}$, we write 
\[
S^{b_j} + S^{b_j+1}g_{0}(S) = \sum_{l=-\infty}^{\infty} \sigma_{-l}^{(j)}(\bold c^{(j)})s_j^l, 
\]
 as before.
Here, $S = 
 s_j(1 + \sum_{l=1}^{\infty}\prod_{i=0}^{l-1}(\frac{1}{a_j}-i)\frac{1}{l!}(\sum_{k=2}^{a_j}\frac{c_k^{(j)}}{s_j^k})^l)$
 and $s_j$ is a local coordinate on $C$ around $p_j$.
In particular, for $-(a_j-1)\leq l\leq -1$, $\sigma_{-l}^{(j)}$ has the form 
 $f_{b_j-l}^{(b_j)} + o_{b_j-l}$ in the notation of Proposition \ref{prop:perturb}.

Recall that we have a subset $\mathcal I$
 of $H^0(C, \varphi^*\omega_X(Z))$, see Definition \ref{def:obstsections}.
The obstruction class associated with the map $\varphi_N$ pairs with these sections.
The map $\varphi_N$ determines the constants $\bold c^{(j)}(N) = (c_2^{(j)}(N), \dots, c_{a_j}^{(j)}(N))$, 
 $c_i^{(j)}(N)\in t^{d_ji}\Bbb C[[t]]$ for each $j = 1, \dots, e$.
Of course, terms of $c_i^{(j)}(N)$ of sufficiently high order with respect to $t$
 does not affect $\varphi_N$. 
So, the map $\varphi_N$ does not determine $\bold c^{(j)}(N)$
 uniquely, but it does not matter in the following argument.
We need to change the values of $\bold c^{(j)}(N)$
 to $\bold c^{(j)}(N+1)$ to cancel the obstruction.


In this subsection, we prove the following.
Recall that $o_{N+1}\in H^1(C, \bar{\mathcal N}_{\varphi})$ denotes the obstruction class to deforming
 $\varphi_N$.
\begin{prop}\label{prop:cal}
Given $\bold c^{(j)}(N)$, $j = 1, \dots, e$, satisfying $(\star_{\eta})$ for any $\eta$ (see Definition \ref{def:star}),
 there is a set of series $\bold c^{(j)}(N+1)\in \Bbb C[[t]]^{a_j-1}$
 which satisfy the following
 conditions.
\begin{enumerate}
\item The equality
\[
c_i^{(j)}(N+1)-c_i^{(j)}(N) = 0\;\;
\begin{cases} \text{mod $t^{d_ji+N+1-d_j(b_j+a_j-1)}$, when $N \geq d_j(b_j+a_j-1)$},\\
 \text{mod $t^{d_ji+1}$, otherwise,}
\end{cases}
\]
 holds for $j = 1, \dots, e$, $i = 2, \dots, a_j$.
\item The equality
\[
\tilde{\sigma}_{-l}^{(j)}(\bold c^{(j)}(N+1))
 = \sigma_{-l}^{(j)}(\bold c^{(j)}(N)),\;\; \text{mod $t^{N+1}$},
\]
 holds for $l<0$, $j = 1, \dots, e$.
\item The equality
\[
 \sum_{j=1}^e Res_{p_j}(\eta, \sum_{l=-\infty}^{-1}(\sigma_{-l}^{(j)}(\bold c^{(j)}(N)) 
  - \tilde{\sigma}_{-l}^{(j)}(\bold c^{(j)}(N+1)))s_j^l\partial_{w_j})\\
 = t^{N+1}(\eta, o_{N+1}),\;\; \text{mod $t^{N+2}$},
\]
 holds for any $\eta\in \mathcal I$.
Here, $(\eta, o_{N+1})$ is the pairing between 
 $H^0(C, \varphi^*\omega_X(Z))$ and 
 $H^1(C, \bar{\mathcal N}_{\varphi}) \cong H^0(C, \varphi^*\omega_X(Z))^{\vee}$.
\end{enumerate}
\end{prop}
\proof
We will construct such $c_i^{(j)}(N+1)$ by a bootstrapping type argument based on Proposition \ref{prop:perturb}.
To use Proposition \ref{prop:perturb}, we consider the equation 
\[
(\ast_{n, \eta})\;\;\;\;
 \sum_{j=1}^e Res_{p_j}(\eta, \sum_{l=-\infty}^{-1}(\sigma_{-l}^{(j)}(\bold c^{(j)}(N)) 
  - \tilde{\sigma}_{-l}^{(j)}(\bold c^{(j)}))s_j^l\partial_{w_j})\\
   = t^{N+1}(\eta, o_{N+1}),\;\; \text{mod $t^{n+1}$},
\]
 for any non-negative integer $n$, not just for $n = N+1$.
Here, $\bold c^{(j)} = (c_2^{(j)}, \dots, c_{a_j}^{(j)})$ are variables such that $c_i^{(j)}$ 
 takes values in $c_i^{(j)}(N) + t^{d_ji+1}\Bbb C[[t]]$.
Recall that $P(\eta) = (j(\eta), m(\eta))$ is the largest element in $psupp(\eta)$
 with respect to the order on the set $\{1, \dots, e\}\times\Bbb Z_{>0}$
 introduced in Definition \ref{def:eZorder}.
In the following lemma, we regard
 $\bold c^{(j(\eta))}$ as variables and $\bold c^{(j)}$, $j\neq j(\eta)$,
 as constants in $(\ast_{n, \eta})$.
Also, recall that $c_i^{(j)}(N)$ is of the form 
 $c_i^{(j)}(N) = t^{d_ji}\tilde c_i^{(j)} + t^{d_ji+1}\Bbb C[[t]]$, 
 where $\tilde c_i^{(j)}\in \Bbb C$.
\begin{lem}\label{lem:N}
The equation $(\ast_{n, \eta})$ can be written in the form 
\begin{equation}\label{eq:*}
\bar f_{b_{j(\eta)}+a_{j(\eta)}-m(\eta)}^{(b_{j(\eta)})}(\bold c^{(j(\eta))})\\
  = t^{d_{j(\eta)}(b_{j(\eta)}+a_{j(\eta)}-m(\eta))}f_{b_{j(\eta)}+a_{j(\eta)}-m(\eta)}^{(b_{j(\eta)})}(\tilde{\bold c}^{(j(\eta))}) 
   + o_{b_{j(\eta)}+a_{j(\eta)} - m(\eta)}(\bold c^{(j(\eta))}),\;\; \text{mod $t^{n+1}$}.
\end{equation}
\end{lem}
\proof
Consider the summand 
 $Res_{p_j}(\eta, \sum_{l=-\infty}^{-1}(\sigma_{-l}^{(j)}(\bold c^{(j)}(N)) 
   - \tilde{\sigma}_{-l}^{(j)}(\bold c^{(j)}))s_j^l\partial_{w_j})$, $j\neq j(\eta)$, 
   in $(\ast_{n, \eta})$.
Since $c_i^{(j)}$ takes values in $c_i^{(j)}(N) + t^{d_ji+1}\Bbb C[[t]]$, 
 by Lemma \ref{lem:ori}, it is easy to see that 
 $Res_{p_j}(\eta, \sum_{l=-\infty}^{-1}(\sigma_{-l}^{(j)}(\bold c^{(j)}(N)) 
    - \tilde{\sigma}_{-l}^{(j)}(\bold c^{(j)}))s_j^l\partial_{w_j})$ 
  is contained in $t^{d_j(b_j+a_j-k_j)+1}\Bbb C[[t]]$,
  where $k_j = \max\{k\;|\; (j, k)\in psupp(\eta)\}$ if $\{k\;|\; (j, k)\in psupp(\eta)\}$ is non-empty and 
   $k_j=0$ otherwise.
Here, we have $d_j(b_j+a_j-k_j) \geq d_{j(\eta)}(b_{j(\eta)} + a_{j(\eta)} - m(\eta))$
 by the construction of $\mathcal I$.
Thus, $Res_{p_j}(\eta, \sum_{l=-\infty}^{-1}(\sigma_{-l}^{(j)}(\bold c^{(j)}(N)) 
    - \tilde{\sigma}_{-l}^{(j)}(\bold c^{(j)}))s_j^l\partial_{w_j})$ 
  can be written in the form $o_{b_{j(\eta)}+a_{j(\eta)} - m(\eta)}(\bold c^{(j(\eta))})$.

On the other hand, the summand 
\[
Res_{p_{j(\eta)}}(\eta, \sum_{l=-\infty}^{-1}(\sigma_{-l}^{(j(\eta))}(\bold c^{(j(\eta))}(N)) 
  - \tilde{\sigma}_{-l}^{(j(\eta))}(\bold c^{(j(\eta))}))s_{j(\eta)}^l\partial_{w_{j(\eta)}})
\]
  can be written in the form 
 $-\bar f_{b_{j(\eta)}+a_{j(\eta)}-m(\eta)}^{(b_{j(\eta)})}(\bold c^{(j(\eta))})
   + f_{b_{j(\eta)}+a_{j(\eta)}-m(\eta)}^{(b_{j(\eta)})}(\bold c^{(j(\eta))}(N))
    + o_{b_{j(\eta)}+a_{j(\eta)} - m(\eta)}(\bold c^{(j(\eta))})$
 up to a multiplicative constant
 by a same argument as in Lemmas \ref{lem:leading}.
Note that the term $f_{b_{j(\eta)}+a_{j(\eta)}-m(\eta)}^{(b_{j(\eta)})}(\bold c^{(j(\eta))}(N))$ 
 comes from $\sigma_{-l}^{(j(\eta))}(\bold c^{(j(\eta))}(N))$, $l = -a_{j(\eta)}+m(\eta)$, 
 and it can be written in the form 
 $t^{d_{j(\eta)}(b_{j(\eta)}+a_{j(\eta)}-m(\eta))}f_{b_{j(\eta)}+a_{j(\eta)}-m(\eta)}^{(b_{j(\eta)})}(\tilde{\bold c}^{(j(\eta))}) 
    + o_{b_{j(\eta)}+a_{j(\eta)} - m(\eta)}(\bold c^{(j(\eta))})$.

Finally, the term $t^{N+1}(\eta, o_{N+1})$
 can be written in the form $o_{b_{j(\eta)}+a_{j(\eta)} - m(\eta)}(\bold c^{(j(\eta))})$
 by Proposition \ref{prop:obsteval}.
\qed\\


Now, we define a system of equations of the form in Proposition \ref{prop:perturb}
 for each $p_j$.
For those $(j, k)$, $j=1, \dots, e$, $1 \leq k\leq a_j-1$, such that there is some $\eta\in \mathcal I$ 
 where $P(\eta) = (j(\eta), m(\eta)) = (j, k)$, 
 we assign the equation given in Lemma \ref{lem:N} above.

For other $(j, k)$, we assign the equation 
\begin{equation}\label{eq:7}
\tilde{\sigma}_{-(a_j-k)}^{(j)}(\bold c^{(j)}) = \sigma_{-(a_j-k)}^{(j)}(\bold c^{(j)}(N)), 
\end{equation}
 mod $t^{n+1}$.
This can be written in the form
\[
\bar f_{b_j+a_j-k}^{(b_j)}(\bold c^{(j)}) = t^{d(b_j+a_j-k)}f_{b_j+a_j-k}^{(b_j)}(\tilde{\bold c}^{(j)})
                                                  + o_{b_j+a_j-k}(\bold c^{(j)}),\;\; \text{\rm{mod} $t^{n+1}$},\;\;  k\in\{1, \dots, a_j-1\}.
\]

\begin{defn}
We write the equation (\ref{eq:*}) by $(S_{j(\eta), m(\eta)})_n$, 
 when it is considered over $\Bbb C[t]/t^{n+1}$.
Similarly, we write the equation (\ref{eq:7}) by $(S_{j, k})_n$, 
 when it is considered over $\Bbb C[t]/t^{n+1}$.
\end{defn}

The equations (\ref{eq:*}) and (\ref{eq:7})
 give a system of equations of the form described in Proposition \ref{prop:perturb}
 for each $j = 1, \dots, e$. 
For example, for the case $j = 1$, 
 let $\eta_1^{(1)}, \dots, \eta_{q_1}^{(1)}$ be the subset of $\mathcal I$
 consisting of those $\eta$ satisfying $j(\eta) = 1$, where $P(\eta) = (j(\eta), m(\eta))$.
Then, the system of equations consists of 
\[
(S_{1, m(\eta_i^{(1)})})_{d_1(b_1+a_1-m(\eta_i^{(1)})) + \beta},\;\; i = 1, \dots, q_1,
\] 
 coming from (\ref{eq:*}) and
\[
(S_{1, k})_{d_1(b_1+a_1-k)) + \beta},\;\; k\in \{1, \dots, a_1-1\}\setminus \{m(\eta_1^{(1)}), \dots, m(\eta_{q_1}^{(1)})\},
\] 
 coming from (\ref{eq:7}).
Here, $\beta$ is an integer.
Note that if we have $i\neq j$, then $m(\eta_i^{(1)})\neq m(\eta_j^{(1)})$.
%
Similarly, let $(S_{j, k})_{d_j(b_j+a_j-k)) + \beta}$, 
 $k=1, \dots, a_j-1$, be the system of equations for the other $j = 2, \dots, e$.

As we noted above, in general, the equation $(S_{j, k})_n$ associated with the point $p_j$
 depends on the variables $\bold c^{(i)} = (c_2^{(i)}, \dots, c_{a_i}^{(i)})$,
 $i\neq j$, attached to the other singular points of $\varphi$.
Here, each $c_l^{(i)}$ takes values in $t^{d_il}\tilde c_l^{(i)} + t^{d_il+1}\Bbb C[[t]]$.
Thus, we need to be careful when we apply Proposition \ref{prop:perturb}
 to the system of equations $(S_{j, k})_n$, $k=1, \dots, a_j-1$,
 for each $j$, so that the solutions provided for each $j$ are compatible with 
 the system of equations $(S_{j', k})_n$, $k=1, \dots, a_{j'}-1$, associated with the other singular points.
This point is addressed in Lemma \ref{lem:beta} below.
In fact, Lemma \ref{lem:beta} asserts that if we modify the values of $c_l^{(i)}$
 as demonstrated in the proof of Proposition \ref{prop:perturb}, 
 the equations $(S_{j, k})_n$, $j\neq i$,
 remain independent of how they are modified.
This is important for the bootstrapping argument used below to work.

\begin{lem}\label{lem:beta}
Let $l\in \{1, \dots, e\}$ and
 let $\beta$ be a positive integer.
Let $n\geq M-1$ be an integer.
If we modify $c_i^{(j)}$, $j> l$, by adding a term $\delta_i\in t^{d_ji+\beta}\Bbb C[[t]]$, 
 the terms in the system of equations $(S_{l, k})_n$, $k=1, \dots, a_l-1$, are modified 
 only at the orders higher than $t^{d_l(b_l+a_l-k)+\beta}$.
\end{lem}
\proof 
If the equation $(S_{l, k})_n$ does not correspond to some $(\ast_{n, \eta})$, $\eta\in \mathcal I$,
 the claim is obvious because the equation does not depend on $c_i^{(j)}$, $j\neq l$.
Therefore, assume that $(S_{l, k})_n$ corresponds to some $(\ast_{n, \eta})$.
If we write $P(\eta) = (j(\eta), m(\eta))$, we have $j(\eta) = l$ and $m(\eta) = k$.
In the equation $(\ast_{n, \eta})$, the part $\tilde{\sigma}_{-q}^{(j)}(\bold c^{(j)})$ depends on 
 $c_i^{(j)}$.

By the construction of $\{\eta_1, \dots, \eta_{\sigma}\}$, if $j>l$, we have 
 $d_l(b_l+a_l-k) < d_j(b_j + a_j-k_j)$, here $k_j = \max\{\kappa \;|\; (j, \kappa)\in psupp(\eta)\}$ 
 if $\{k\;|\; (j, k)\in psupp(\eta)\}$ is non-empty and 
  $k_j=0$ otherwise.
Thus, if we modify $c_i^{(j)}$, $j> l$, by adding a term $\delta_i\in t^{d_ji+\beta}\Bbb C[[t]]$, 
 the part of the term $\tilde{\sigma}_{-q}^{(j)}(\bold c^{(j)})$ which pairs with $\eta$ nontrivially
 will be modified only at the orders at least $d_j(b_j + a_j-k_j) + \beta > d_l(b_l+a_l-k)+\beta$.
This proves the claim.\qed\\

We also remark the following.
\begin{lem}\label{lem:bootsol}
The family $\bold c^{(j)} = \bold c^{(j)}(N)$, $j = 1, \dots, e$,
 gives a solution to all $(S_{l, k})_{d_l(b_l+a_l-k)}$, $l = 1, \dots, e$, $k = 1, \dots, a_l-1$.
\end{lem}
\proof
If the equation $(S_{l, k})_{d_l(b_l+a_l-k)}$ does not come from an element of $\mathcal I$, that is, 
 it is of the form Eq.(\ref{eq:7}), the claim follows from Lemma \ref{lem:ori}.
Assume that the equation $(S_{l, k})_{d_l(b_l+a_l-k)}$ is of the form 
 $(\ast_{d_l(b_l+a_l-k), \eta})$ for some $\eta\in\mathcal I$.
In this case, we have $l = j(\eta)$ and $k = m(\eta)$, where $P(\eta) = (j(\eta), m(\eta))$ as we noted above.
Recall that the term $t^{N+1}(\eta, o_{N+1})$ can be written in the form 
 $o_{b_{j(\eta)}+a_{j(\eta)} - m(\eta)}(\bold c^{(j(\eta))})$.
In particular, if $(\eta, o_{N+1})$ does not vanish, we have
 $N+1\geq d_{j(\eta)}(b_{j(\eta)}+a_{j(\eta)} - m(\eta))+1$.
Thus, over $\Bbb C[t]/t^{d_{j(\eta)}(b_{j(\eta)}+a_{j(\eta)}-m(\eta))+1}$, 
 $t^{N+1}(\eta, o_{N+1})$ is zero.
The claim follows from this observation and Lemma \ref{lem:ori}.\qed\\

Now, we will construct a solution to these systems of equations 
 over $\Bbb C[t]/t^{N+2}$ by a bootstrapping type argument. 
Let us consider the case $j = 1$.
Consider $(S_{1, k})_{d_1(b_1+a_1-k)+1}$, $k = 1, \dots, a_1-1$, as equations with the variables
 $\bold c^{(1)}$, and with $\bold c^{(j)}
  = \bold c^{(j)}(N)$, $j\neq 1$, regarded as constants.
Then, by applying Proposition \ref{prop:perturb} to the solution given in
 Lemma \ref{lem:bootsol}, this system of equations 
 has a solution $\bold c^{(1)}[1]$.

Next, consider the case $j = 2$.
As in Lemma \ref{lem:bootsol},
 $\bold c^{(j)}(N)$, $j = 1, \dots, e$, give a solution to $(S_{2, k})_{d_2(b_2+a_2-k)}$, $k = 1, \dots, a_2-1$.
This is true even after we change the value of $\bold c^{(1)}$ from $\bold c^{(1)}(N)$ to $\bold c^{(1)}[1]$.
Namely, we have $c_i^{(1)}(N)- c_i^{(1)}[1]\in t^{d_1i+1}\Bbb C[[t]]$, and this implies that 
 modifying $\bold c^{(1)}(N)$
  to $\bold c^{(1)}[1]$ possibly 
 changes the equation $(S_{2, k})_{d_2(b_2'+a_2-k)}$ only at the order higher than $t^{d_2(b_2+a_2-k)}$.
It follows that it in fact does not change the equation
 at all, since $(S_{2, k})_{d_2(b_2+a_2-k)}$ is an equation defined over
 $\Bbb C[t]/t^{d_2(b_2+a_2-k)+1}$.
 
Now, consider the equations $(S_{2, k})_{d_2(b_2+a_2-k)+1}$, $k = 1, \dots, a_2-1$, as equations with the variables
 $\bold c^{(2)}$, and with $\bold c^{(j)}
  = \bold c^{(j)}(N)$, $j\neq 1, 2$, 
 and $\bold c^{(1)}
   = \bold c^{(1)}[1]$
   regarded as constants.
Again, by Proposition \ref{prop:perturb}, this system of equations 
 has a solution $\bold c^{(2)}[1]$. 
Here, an important point is that since $c_i^{(2)}[1]$ is obtained by 
 adding some element in $t^{d_2i+1}\Bbb C[[t]]$ to $c_i^{(2)}(N)$, it does not change 
 the equation $(S_{1, k})_{d_1(b_1+a_1-k)+1}$, which is defined over $\Bbb C[t]/t^{d_1(b_1+a_1-k)+2}$,
 by Lemma \ref{lem:beta}.
In particular, $\bold c^{(1)}[1]$ remains as a solution
 to the system $(S_{1, k})_{d_1(b_1+a_1-k)+1}$, $k = 1, \dots, a_1-1$, although the term $\bold c^{(2)}$
 is replaced from 
 $\bold c^{(2)}(N)$ to $\bold c^{(2)}[1]$.
Repeating the same procedure for $j = 1, \dots, e$, we obtain a solution 
 $\bold c^{(j)}[1]$ to all $(S_{j, k})_{d_j(b_j+a_j-k)+1}$, $k = 1, \dots, a_j-1$.

Then, let us return to the case $j = 1$ and construct a solution $\bold c^{(1)}[2]$ to 
 $(S_{1, k})_{d_1(b_j+a_j-k)+2}$ by adding elements in $t^{d_1i+2}\Bbb C[[t]]$ to $c_i^{(1)}[1]$.
Here, the elements $c_i^{(j)}=c_i^{(j)}[1]$, $j\geq 2$, are considered as constants.
This can be achieved through Proposition \ref{prop:perturb}.
Since modifying $c_i^{(1)}[1]$ in this way changes the equations 
 $(S_{j, k})_{d_j(b_j+a_j-k)+1}$, $j\geq 2$, only at the orders
 higher than $d_j(b_j+a_j-k)+1$, it does not affect the validity of the solution 
 $\bold c^{(j)}=\bold c^{(j)}[1]$ to $(S_{j, k})_{d_j(b_j+a_j-k)+1}$, $j\geq 2$.
Then, repeat the procedure for all $j = 1, \dots, e$, and return to $j = 1$ again.
Iterating this $q$ times until 
\[
\min_{j\in\{1, \dots, e\}, k\in \{1, \dots, a_j-1\}}\{d_j(b_j+a_j-k)+q\} =
 \min_{j\in\{1, \dots, e\}}\{d_j(b_j+1)+q\} =  N+1
\]
 holds, we obtain the required solution $\bold c^{(j)}(N+1)=\bold c^{(j)}[q]$.  
When we have $N\geq d_j(b_j+a_j-1)$, instead of $(S_{j, k})_{d_j(b_j+a_j-k)}$, 
 we can start solving the system of equations from 
 $(S_{j, k})_{d_j(b_j+a_j-k) + N-d_j(b_j+a_j-1)}$.
Then, by Proposition \ref{prop:perturb}, the solution satisfies (1) of the claim of Proposition \ref{prop:cal}.

Finally, we show that
 this solution satisfies the claim of Proposition \ref{prop:cal} (2).
We note the following.
\begin{lem}\label{lem:OB}
For those $l$ satisfying $-(a_j-1)\leq l\leq -1$, which are not contained in the set
 $\{-(a_{j(\eta)}-m(\eta))\, | \, \eta\in \mathcal I \ \text{{\rm satisfies}}\ j(\eta) = j\}$, 
 the claim (2) holds.
\end{lem}
\proof
This is clear since for such $l$, we have 
 $\tilde{\sigma}_{-l}^{(j)}(\bold c^{(j)}(N+1)) = \sigma_{-l}^{(j)}(\bold c^{(j)}(N)), 
  \;\; \text{\rm{mod} $t^{N+2}$}$, by construction.\qed\\

Consider the case when $l = -(a_{j(\eta_0)}-m(\eta_0))$, where $P(\eta_0) = (j(\eta_0), m(\eta_0))$ is the smallest among
 $P(\eta)$, $\eta\in\mathcal I$, with respect to the order introduced in Definition \ref{def:eZorder}.
In this case, by substituting $c_i^{(j)} = c_i^{(j)}(N+1)$ to the equation $(*_{N, \eta_0})$, 
 we have
\[
\sum_{j=1}^eRes_{p_{j}}(\eta_0, \sum_{l=-\infty}^{-1}(\sigma_{-l}^{(j)}(\bold c^{(j)}(N)) 
  - \tilde{\sigma}_{-l}^{(j)}(\bold c^{(j)}(N+1)))s_{j}^l\partial_{w_{j}})  
   = 0,\;\; \text{mod $t^{N+1}$}.
\]
Let us write $k_j = \max\{k\;|\; (j, k)\in psupp(\eta_0)\}$.
If there is no positive integer $k$ such that $(j, k)$ belongs to $psupp(\eta_0)$, we put $k_j = 0$ as before.
Then, we have
\[
\begin{array}{l}
Res_{p_{j}}(\eta_0, \sum_{l=-\infty}^{-1}(\sigma_{-l}^{(j)}(\bold c^{(j)}(N)) 
  - \tilde{\sigma}_{-l}^{(j)}(\bold c^{(j)}(N+1)))s_{j}^l\partial_{w_{j}})\\
  = Res_{p_{j}}(\eta_0, \sum_{l=-\infty}^{-(a_j-k_j)}(\sigma_{-l}^{(j)}(\bold c^{(j)}(N)) 
    - \tilde{\sigma}_{-l}^{(j)}(\bold c^{(j)}(N+1)))s_{j}^l\partial_{w_{j}}),\;\; \text{mod $t^{N+1}$}.

\end{array}
\]
Note that we have $k_{j(\eta_0)} = m(\eta_0)$.

\begin{lem}
Any integer $l$ satisfying $-(a_j-1)\leq l\leq -(a_j-k_j)$ is contained in the range of 
 Lemma \ref{lem:OB} for any $j\neq j(\eta_0)$.
That is, for any $j\neq j(\eta_0)$,  the set of integers $\{-(a_j-1), \dots, -1\}$ is disjoint from the set 
 $\{-(a_{j(\eta)}-m(\eta))\, | \, \eta\in \mathcal I \ \text{{\rm satisfies}}\ j(\eta) = j\}$.
\end{lem}
\proof
By definition of $P(\eta_0)$, we have
\[
d_{j(\eta_0)}(b_{j(\eta_0)}+a_{j(\eta_0)}-m(\eta_0))\leq d_j(b_j+a_j-m),
\]
 for any $(j, m)\in psupp(\eta_0)$.
Moreover, the inequality is strict when $j>j(\eta_0)$ holds.
If there is an $\eta_1\in \mathcal I$ satisfying
\[
j(\eta_1)>j(\eta_0)
\]
 and 
\[
m(\eta_1)\leq k_{j(\eta_1)}
\] 
 so that $l = -(a_{j(\eta_1)}-m(\eta_1))$
 is in the range 
\[
-(a_{j(\eta_1)}-1)\leq l\leq -(a_{j(\eta_1)}-k_{j(\eta_1)}),
\] 
 we have 
\[
\begin{array}{ll}
ord(\eta_0) &= d_{j(\eta_0)}(b_{j(\eta_0)}+a_{j(\eta_0)}-m(\eta_0))\\
 &<d_{j(\eta_1)}(b_{j(\eta_1)}+a_{j(\eta_1)}-k_{j(\eta_1)}) \\
 & \leq d_{j(\eta_1)}(b_{j(\eta_1)}+a_{j(\eta_1)}-m(\eta_1)) = ord(\eta_1).
\end{array}
\]
However, this contradicts to the assumption that $P(\eta_0)$ is the smallest among $P(\eta)$, $\eta\in \mathcal I$.
This proves the lemma for $j>j(\eta_0)$.

If there is an $\eta_2\in\mathcal I$ satisfying 
\[
j(\eta_0)>j(\eta_2)
\]
 and 
\[
-(a_{j(\eta_2)}-1)\leq -(a_{j(\eta_2)} - m(\eta_2))\leq -(a_{j(\eta_2)}-k_{j(\eta_2)}),
\] 
 we have 
\[
\begin{array}{ll}
ord(\eta_0) &= d_{j(\eta_0)}(b_{j(\eta_0)}+a_{j(\eta_0)}-m(\eta_0))\\
 &\leq d_{j(\eta_2)}(b_{j(\eta_2)}+a_{j(\eta_2)}-k_{j(\eta_2)}) \\
 & \leq d_{j(\eta_2)}(b_{j(\eta_2)}+a_{j(\eta_2)}-m(\eta_2)) = ord(\eta_2).
\end{array}
\]
If this inequality is strict, we have $P(\eta_2)<P(\eta_0)$.
If the equality $ord(\eta_0)=ord(\eta_2)$ holds, since we have $j(\eta_0)>j(\eta_2)$, 
 again we have $P(\eta_2)<P(\eta_0)$. 
This contradicts to the assumption, and the proof is complete.\qed\\

Thus, we have
\[
\begin{array}{l}
Res_{p_{j}}(\eta_0, \sum_{l=-\infty}^{-1}(\sigma_{-l}^{(j)}(\bold c^{(j)}(N)) 
  - \tilde{\sigma}_{-l}^{(j)}(\bold c^{(j)}(N+1)))s_{j}^l\partial_{w_{j}})\\
  = Res_{p_{j}}(\eta_0, \sum_{l=-\infty}^{-a_j}(\sigma_{-l}^{(j)}(\bold c^{(j)}(N)) 
    - \tilde{\sigma}_{-l}^{(j)}(\bold c^{(j)}(N+1)))s_{j}^l\partial_{w_{j}}),\;\; \text{mod $t^{N+1}$},
\end{array}
\] 
 for $j\neq j(\eta_0)$.
Now, we note the following.

\begin{lem}
For any $j = 1, \dots, e$, we have 
\[
\sigma_{-l}^{(j)}(\bold c^{(j)}(N)) 
     - \tilde{\sigma}_{-l}^{(j)}(\bold c^{(j)}(N+1)) = 0,\;\; \text{\rm{mod} $t^{N+1}$},
\]
 if $l<-(a_{j}-1)$.
\end{lem}
\proof
In the range $N<d_{j}(b_{j}+a_{j}-1)$, 
 since $\sigma_{-l}^{(j)}(\bold c^{(j)}(N))$
 and $\tilde{\sigma}_{-l}^{(j)}(\bold c^{(j)}(N+1))$ 
 has the order at least $d_{j}(b_{j}-l) > d_{j}(b_{j}+a_{j}-1)\geq N+1$, 
 the equality is obvious.
In the range $N\geq d_{j}(b_{j}+a_{j}-1)$, we have 
\[
c_i^{(j)}(N+1)-c_i^{(j)}(N)\in t^{d_{j}i+N+1-d_{j}(b_{j}+a_{j}-1)}\Bbb C[[t]]
\]
 by the claim (1), 
 and it follows that 
\[
\sigma_{-l}^{(j)}(\bold c^{(j)}(N)) 
      - \tilde{\sigma}_{-l}^{(j)}(\bold c^{(j)}(N+1))
      \in t^{N+1-d_{j}(b_{j}+a_{j}-1) + d_{j}(b_{j}-l)}\Bbb C[[t]].
\]
Since we have 
\[
N+1-d_{j}(b_{j}+a_{j}-1) + d_{j}(b_{j}-l)
  = N+1-d_{j}(a_{j}+l-1) > N+1,
\] 
 the equality follows.\qed\\

Thus, we have
\[
\begin{array}{l}
\sum_{j=1}^eRes_{p_{j}}(\eta_0, \sum_{l=-\infty}^{-1}(\sigma_{-l}^{(j)}(\bold c^{(j)}(N)) 
  - \tilde{\sigma}_{-l}^{(j)}(\bold c^{(j)}(N+1)))s_{j}^l\partial_{w_{j}})  \\
  = Res_{p_{j(\eta_0)}}(\eta_0, \sum_{l=-\infty}^{-(a_{j(\eta_0)}-m(\eta_0))}(\sigma_{-l}^{(j(\eta_0))}(\bold c^{(j(\eta_0))}(N)) 
   - \tilde{\sigma}_{-l}^{(j(\eta_0))}(\bold c^{(j(\eta_0))}(N+1)))s_{j(\eta_0)}^l\partial_{w_{j(\eta_0)}})\\
  = 0,\;\; \text{mod $t^{N+1}$}.
\end{array}
\]

Moreover, for those $l$ in $\{-(a_{j(\eta_0)}-1), \dots, -(a_{j(\eta_0)}-m(\eta_0))-1\}$ , 
 we have $\tilde{\sigma}_{-l}^{(j(\eta_0))}(\bold c^{(j(\eta_0))}(N+1)) = \sigma_{-l}^{(j(\eta_0))}(\bold c^{(j(\eta_0))}(N)), 
   \;\; \text{\rm{mod} $t^{N+1}$}$, by Lemma \ref{lem:OB}.
Thus, we have
\[
\begin{array}{l}
Res_{p_{j(\eta_0)}}(\eta_0, \sum_{l=-\infty}^{-(a_{j(\eta_0)}-m(\eta_0))}(\sigma_{-l}^{(j(\eta_0))}(\bold c^{(j(\eta_0))}(N)) 
    - \tilde{\sigma}_{-l}^{(j(\eta_0))}(\bold c^{(j(\eta_0))}(N+1)))s_{j(\eta_0)}^l\partial_{w_{j(\eta_0)}})\\
   = Res_{p_{j(\eta_0)}}(\eta_0, \sigma_{a_{j(\eta_0)}-m(\eta_0)}^{(j(\eta_0))}(\bold c^{(j(\eta_0))}(N)) 
       - \tilde{\sigma}_{a_{j(\eta_0)}-m(\eta_0)}^{(j(\eta_0))}
            (\bold c^{(j(\eta_0))}(N+1))s_{j(\eta_0)}^{-(a_{j(\eta_0)}-m(\eta_0))}\partial_{w_{j(\eta_0)}})\\
   = 0,\;\; \text{mod $t^{N+1}$},
\end{array}
\]
Since $(j(\eta_0), m(\eta_0))\in psupp(\eta_0)$, this equality implies
\[
\sigma_{a_{j(\eta_0)}-m(\eta_0)}^{(j(\eta_0))}(\bold c^{(j(\eta_0))}(N)) 
        - \tilde{\sigma}_{a_{j(\eta_0)}-m(\eta_0)}^{(j(\eta_0))}(\bold c^{(j(\eta_0))}(N+1)) = 0, \;\; \text{mod $t^{N+1}$}.
\]

We can apply the same argument to other elements in $\mathcal I$ in order from the smallest with respect to the 
 order of $P(\eta)$, and conclude that
 $\tilde{\sigma}_{-l}^{(j)}(\bold c^{(j)}(N+1)) = \sigma_{-l}^{(j)}(\bold c^{(j)}(N))$
 mod $t^{N+1}$ holds for all $l\in \{-(a_j-1), \dots, -1\}$. 
This finishes the proof of Proposition \ref{prop:cal}.\qed


\subsubsection{Proof of Theorem \ref{thm:main}}\label{subsec:main}
Finally, we can prove the following.
\begin{prop}\label{prop:qed}
There is an $(N+1)$-th order deformation $\varphi_{N+1}$ of $\varphi$
 which reduces to $\varphi_N$ over $\Bbb C[t]/t^{N'}$ for some $N'\leq N+1$.
Moreover, in the range $N\geq \max_j\{d_j(b_j+a_j-1)\}$, we can take $N' = N+3-\max_j\{d_j(b_j+a_j-1)\}$.
\end{prop}
\proof
First, we will show that using the parameters $c_i^{(j)}(N+1)$ determined in Proposition \ref{prop:cal},
 we can construct a deformation 
 $\bar{\varphi}_N$ of $\varphi$ over $\Bbb C[t]/t^{N+1}$.
On a neighborhood of $p_j$, we take a curve defined by the parameterization
\begin{equation}\label{eq:2}
\begin{array}{l}
(z_j, w_j) = (s_j^{a_j} + c_2^{(j)}(N+1)s_j^{a_j-2} + \cdots + c_{a_j}^{(j)}(N+1), \\
 \hspace{1in}
  \sum_{l=0}^{\infty}\sigma_{-l}^{(j)}(\bold c^{(j)}(N+1))s_j^l 
  + t^M\bar H_N(s_j, t)),
 \;\; \text{mod $t^{N+2}$},
\end{array}
\end{equation}
 where $c_i^{(j)}(N+1)$ is the one 
 satisfying the equations of Proposition \ref{prop:cal},
 and $\bar H_N(s_j, t)$ is determined from it as in Section \ref{subsec:virtdeform}.
 
 \begin{lem}\label{lem:cimage}
The image of the curve given by the above parameterization Eq.(\ref{eq:2})
 coincides with the restriction of the image
 of $\varphi_N$ to a neighborhood of $p_j$ over $\Bbb C[t]/t^{N+1}$. 
 \end{lem} 
\proof  
This follows from 
 the calculation Eq.(\ref{eq:diamond}) in Section \ref{subsec:virtuallocal}
 and Proposition \ref{prop:cal} (2).\qed\\
 
Moreover, for some $N'$ with $N'\leq N+1$, we have $c_i^{(j)}(N+1) = c_i^{(j)}(N)$
 mod $t^{N'}$, $i = 2, \dots, a_j$.
In this case, we have $s_j(N+1) = s_j$, mod $t^{N'}$, and 
 $\bar H_N(s_j, t) = H_N(s_j, t)$, mod $t^{N'}$.
Thus, over $\Bbb C[t]/t^{{N'}}$, the above parameterization and $\varphi_N$ are the same 
 even as maps.
If $N' = N+1$ for all $j = 1, \dots, e$, we can take $\bar{\varphi}_N = \varphi_N$.
So, let us assume we have $N'\leq N$ for some $j$.
We write this map as $\varphi_{N'-1}|_{U_{p_j}}$.

Now, on a neighborhood $U_{p_j}$ of $p_j$, we have two deformations of the map $\varphi_{N'-1}|_{U_{p_j}}$.
The first one $\varphi_{N', p_j}$
 is obtained by restricting $\varphi_N$ to $U_{p_j}$ and reducing it to a map over $\Bbb C[t]/t^{N'+1}$.
The other one $\bar{\varphi}_{N', p_j}$
 is derived from the above parameterization using $c_i^{(j)}(N+1)$ 
 by reducing it 
 to a map over $\Bbb C[t]/t^{N'+1}$.
These two deformations can be different as maps, but they have the same image by Lemma \ref{lem:cimage}.
On the other hand, on the complement of the singular points $\{p_1, \dots, p_e\}$ of $\varphi$, 
 we have a deformation $\varphi_{N', c}$ of 
 $\varphi_{N'-1}$ given by the restriction of $\varphi_N$.
 
The deformations $\varphi_{N', p_j}$ and $\varphi_{N', c}$ clearly glue to a map $\varphi_{N'}$, 
 since they are the restrictions of the given map $\varphi_N$.
On the other hand, since $\varphi_{N', p_j}$ and $\bar{\varphi}_{N', p_j}$ have the same image, 
 the difference between $\bar{\varphi}_{N', p_j}$ and $\varphi_{N', c}$ on the overlap gives a section of 
 the tangent sheaf of the curve $C$.
In particular, it is zero as a section of $\bar{\mathcal N}_{\varphi}$ where the obstruction
 takes value.
Therefore, possibly after deforming the domain curve of $\varphi_{N'}$,
 the local deformations $\bar{\varphi}_{N', p_j}$ and $\varphi_{N', c}$ also glue and
 give a global map $\bar{\varphi}_{N'}$.
Note that the images of $\varphi_{N'}$ and $\bar{\varphi}_{N'}$ are the same.
Note also that the map $\bar{\varphi}_{N'}$ still has the parameterization (\ref{eq:2})
 around $p_j$.
Namely, we have not changed the locally defined maps $\bar{\varphi}_{N', p_j}$ and $\varphi_{N', c}$, 
 up to automorphisms,
 but just changed the gluing of the domain.

Now, we will deform $\bar\varphi_{N'}$.
Although the domain curve of $\bar\varphi_{N'}$ may be different from 
 that of $\varphi_{N'}$,
 the restriction of the maps $\varphi_{N'}$ and $\bar\varphi_{N'}$ 
 to the complement of $\{p_1, \dots, p_e\}$ are equivalent up to automorphisms
 by construction.
Then, again by Lemma \ref{lem:cimage}, the argument in the above paragraph still applies, 
 and we obtain a map $\bar\varphi_{N'+1}$ deforming $\bar\varphi_{N'}$.
By repeating this process, we obtain a map $\bar{\varphi}_{N}$.
The image of $\bar{\varphi}_{N}$ is again the same as that of $\varphi_N$.

Now, we prove that there is a deformation $\varphi_{N+1}$ of $\bar{\varphi}_{N}$.
The above parameterization using $c_i^{(j)}(N+1)$
 gives a local deformation of $\bar{\varphi}_{N}$ on a neighborhood of $p_j$.
By the above construction, the restrictions of 
 $\bar{\varphi}_{N}$ and $\varphi_N$ to the complement of $\{p_1, \dots, p_e\}$ are identical. 
Thus, on the complement of $\{p_1, \dots, p_e\}$, we can take the same local deformation 
 of $\bar{\varphi}_{N}$ as that of $\varphi_N$.
Recall that the obstruction to deforming $\varphi_N$ was given by the class $o_{N+1}$.
The local deformations of $\bar{\varphi}_{N}$ and $\varphi_N$ differ only on a neighborhood of $p_j$, 
 and their difference is given by the calculation $(\ref{eq:diamond})$.
Then, Proposition \ref{prop:cal} (3) shows that the difference between these local deformations
 of $\bar{\varphi}_{N}$ and $\varphi_N$ gives a meromorphic section of $\bar{\mathcal N_{\varphi}}$
 which, through the calculation of Proposition \ref{prop:pair},
 gives a local contribution to the obstruction cancelling the given $o_{N+1}$.
Thus, the obstruction to deforming $\bar{\varphi}_{N}$ vanishes.

Finally, the assertion about $N'$ follows from Proposition \ref{prop:cal} (1)
 concerning the order of $c_i(N+1)-c_i(N)$.\qed\\
 
Assume that we have constructed a deformation $\varphi_N$ of $\varphi$, for some $N\geq M-1$.
Let $\psi_{N'-1}$ be the reduction of $\varphi_N$ to a map over $\Bbb C[t]/t^{N'}$, where
 $N' = 3$ when $N\leq \max_j\{d_j(b_j+a_j-1)\}$, and 
 $N' = N+3-\max_jd_j(b_j+a_j-1)$ when $N\geq \max_j\{d_j(b_j+a_j-1)\}$.
By Proposition \ref{prop:qed}, there is a map $\varphi_{N+1}$ and let $\psi_{N'}$ be the reduction of 
 $\varphi_{N+1}$ to a map over $\Bbb C[t]/t^{N'+1}$.
By Proposition \ref{prop:qed}, $\psi_{N'}$ is a deformation of $\psi_{N'-1}$, though
 $\varphi_{N+1}$ may not be a deformation of $\varphi_N$.
Thus, we obtain a deformations $\psi_n$ of $\varphi$, $n\in \Bbb N$, up to any order.
Applying a suitable algebraization result \cite{Ar}, this finishes the proof of Theorem \ref{thm:main}.

\section{Application of the main theorem}\label{sec:examples}

\subsection{Condition ${\rm (G)}$ and deformations in general}\label{subsec:(G)}
Let $\varphi\colon C\to X$ be a semiregular map from a curve to a surface as before, 
 and let $\{p_1, \dots, p_e\}$ be the set of singular points of $\varphi$. 
Since there is no specific reason that the functions $f_{b_j+i}^{(b_j)}$ 
 expressing the condition in Theorem \ref{thm:main} have non-trivial relations, 
 it would not be too optimistic to expect
 that these and the associated functions $\bar f_{b_j+i}^{(b_j)}$, $F_{-n}$ behave like generic ones.
If this is true, the criterion for the existence of deformations becomes largely independent of 
 the surface $X$ and the map $\varphi$.
Namely, we can formulate the condition for the existence of deformations as follows.
\begin{defn}\label{def:(G)}
Using the notation of Definition \ref{def:perturb}, we say that the polynomials 
 $f_{b+1}^{(b)}, \dots, f_{b+a-1}^{(b)}$ satisfy the condition (${\rm G}$)
 if for any $i$, $1\leq i\leq a-1$, there is a point $\tilde{\bold c}\in \Bbb C^{a-1}$
 such that 
\[
F_{-n}(\tilde{\bold c}) = 0,\;\; \forall n\in \{1, \dots, a-1\}\setminus\{i\}, \;\; F_{-i}(\tilde{\bold c})\neq 0,
\]
 and $f_{b+1}^{(b)}, \dots, f_{b+a-1}^{(b)}$ satisfy the condition (${\rm T}$) at $\tilde{\bold c}$.
Since this condition only depends on the pair of integers $(a, b)$, where $b>a$ and $a$ does not divide $b$, 
 we also say that the pair $(a, b)$ satisfies the condition ({\rm G}).
Recall that for each singular point $p$ of $\varphi$, such a pair of integers $(a, b)$ is determined.
If this pair satisfies the condition ({\rm G}), we will say that the singular point $p$
 satisfies the condition ({\rm G}).
\end{defn}
\begin{defn}\label{def:D}
We will say that a point $p_j\in \{p_1\, \dots, p_e\}$ satisfies the condition (D)
 if the inequality
\[
\dim H^0(C, \varphi^*\omega_X((a_j-1)p_j)) < \dim H^0(C, \varphi^*\omega_X) + a_j-1 
\]
 holds.
Here, $a_j-1$ is the coefficient of $p_j$ in the ramification divisor $Z$ of $\varphi$. 
\end{defn}

\begin{thm}\label{thm:gen2}
If the conditions ${\rm (D)}$ and ${\rm (G)}$ hold at each $p_j$, the map $\varphi$ deforms.
\end{thm}
\proof
First, we define the set of $e$-tuple of positive integers
\[
\mathcal L = \{(l_1, \dots, l_e)\;|\; l_j\in \{1, \dots, a_j-1\}\} 
\]
 in the following way.
Namely, $(l_1, \dots, l_e)\in\mathcal L$ if and only if there is a direct sum decomposition 
\[
H^0(C, \varphi^*\omega_X(Z))=H^0(C, \varphi^*\omega_X) \oplus H_0\oplus H_1
\]
 such that,
\begin{enumerate}
\item[(i)] if $\eta\in H_1$ and $(j, m)\in psupp(\eta)$, we have $m<l_j$, 
\item[(ii)] if $\eta\in H_0$ and $(j, m)\in psupp(\eta)$, we have $m\geq l_j$, and
\item[(iii)] there is no element $\eta\in H_0$ such that $psupp(\eta) = \{(j, l_j)\}$ for some $j\in \{1, \dots, e\}$.
\end{enumerate}
Using this notation, we prove the following.
\begin{lem}\label{lem:setL}
If the condition ${\rm (D)}$ holds at each $p_j$, 
 the set $\mathcal L$ is not empty.
\end{lem}
\proof
Let $\mathcal I_j$ be the subset of $\{1, \dots, a_j-1\}$ such that $k\in \mathcal I_j$ if and only if
 there is a section $\eta\in H^0(C, \varphi^*\omega_X(Z))$ satisfying $psupp(\eta) = \{(j, k)\}$.
By the condition {\rm (D)}, we have $I_j\neq \{1, \dots, a_j-1\}$ for each $j$.
If $\min_{l\in I_j} l > 1$ or $I_j = \emptyset$, take $l_j = 1$.
If $\min_{l\in I_j} l = 1$, take $l_j = \min\{\{1, \dots, a_j-1\}\setminus I_j\}$.
Then, there is a direct sum decomposition 
 $H^0(C, \varphi^*\omega_X(Z)) = H^0(C, \varphi^*\omega_X)\oplus H_0\oplus H_1$
 associated with 
 $(l_1, \dots, l_e)$ which satisfies the conditions (i), (ii) and (iii) above.
Namely, $H_0$ and $H_1$ are uniquely determined, modulo $H^0(C, \varphi^*\omega_X)$,
 by the properties (i) and (ii).\qed\\

We fix an element $(l_1, \dots, l_e)$ of $\mathcal L$.
Also, fix $\tilde{\bold c}^{(j)}\in\Bbb C^{a_j-1}\setminus\{0\}$ for each $j = 1, \dots, e$, which satisfies
\[
F_{-n}(\tilde{\bold c}^{(j)}) = 0,\;\; n\in \{1, \dots, a_j-1\}\setminus\{a_j-l_j\}, \;\; F_{-(a_j-l_j)}(\tilde{\bold c}^{(j)})\neq 0,
\] 
 and $f_{b_j+1}^{(b_j)}, \dots, f_{b_j+a_j-1}^{(b_j)}$ satisfy the condition (${\rm T}$) at $\tilde{\bold c}_j$.
Let $M$ be the least common multiple of $b_j+a_j-l_j$, $j = 1, \dots, e$,  
 and define the integer $d_j$ by $d_j = \frac{M}{b_j+a_j-l_j}$.
Introduce a total order to the set $\{1, \dots, e\}\times \Bbb Z_{> 0}$
 by the rule that $(j, m)>(j', m')$ if and only if
 \begin{enumerate}
 \item $d_j(b_j + a_j - m) < d_{j'}(b_{j'}'+a_{j'} - m')$, or
 \item $d_j(b_j + a_j - m) = d_{j'}(b_{j'}'+a_{j'} - m')$ and $j>j'$.
 \end{enumerate}
For $\eta\in H^0(C, \varphi^*\omega_X(Z))$,
 let $P(\eta) = (j(\eta), m(\eta))$ be the maximal element of $psupp(\eta)$
 with respect to this order.
Using this notation, let us define $ord(\eta) \in \Bbb Z$ by
 \[
 ord(\eta) = d_{j(\eta)}(b_{j(\eta)}+a_{j(\eta)} - m(\eta)).
 \] 
We set $ord(\eta) = \infty$ if $\eta$ belongs to $H^0(C, \varphi^*\omega_X)$.

We take a basis $\{\lambda_1, \dots, \lambda_a, \mu_1, \dots, \mu_b, \nu_1, \dots, \nu_c\}$
 of $H^0(C, \varphi^*\omega_X(Z))$, 
 where $\{\lambda_1, \dots, \lambda_a\}$, $\{\mu_1, \dots, \mu_b\}$ and $\{\nu_1, \dots, \nu_c\}$
 are bases of $H^0(C, \varphi^*\omega_X)$, $H_1$ and $H_0$, respectively.
We take $\{\lambda_1, \dots, \lambda_a\}$ and $\{\mu_1, \dots, \mu_b\}$ arbitrarily.
We take $\{\nu_1, \dots, \nu_c\}$ as in Section \ref{subsec:obbasis}.
Namely, for a positive integer $N\leq M$, let $V_N$ be the subspace of $H^0(C, \varphi^*\omega_X(Z))$
 defined by 
\[
V_N = \{\eta\in H^0(C, \varphi^*\omega_X(Z))\;|\; ord(\eta)\geq N\}.
\]
Let 
\[
H^0(C, \varphi^*\omega_X)\oplus H_1
 \subset V_{i_k}\subset \cdots V_{i_{2}}\subset V_{i_1} = H^0(C, \varphi^*\omega_X(Z))
\]
 be the maximal strictly increasing subsequence, here $i_k = M$.

We have a refinement
\[
V_{i_{j+1}} \subset V_{i_{j}, n_1}\subset V_{i_{j}, n_2}\subset \cdots\subset V_{i_{j}, n_{u_j}} = V_{i_{j}}
\]
 such that 
\[
\dim V_{i_{j}, n_{r+1}} = \dim V_{i_{j}, n_r} + 1,\;\; r = 0, \dots, u_j-1,
\]
 as in Section \ref{subsec:obbasis}.
Then, we define $\{\nu_1, \dots, \nu_c\}$ by successively choosing a general element of $V_{i_j, n_k}$
 as in Section \ref{subsec:obbasis}.

For any $\mu_i$, the condition $(\star_{\mu_i})$ of Definition \ref{def:star} obviously holds, since all the relevant
 $F_{-n}(\tilde{\bold c})$ are zero.
Take any $\nu_i$.
By the conditions (i), (ii) and (iii) above, $\nu_i$ satisfies either
\begin{enumerate}
\item[(a)] $ord(\nu_i) < M$, or
\item[(b)] $ord(\nu_i) = M$, and $\sharp(psupp(\nu_i)\cap \{(1, l_1), \dots, (e, l_e)\})\geq 2$.
\end{enumerate}
In the case (a), the condition $(\star_{\nu_i})$ obviously holds again.
In the case (b), we note that the condition $(\star_{\nu_i})$ is of the form 
\[
\sum_{\{j\;|\; (j, l_j)\in psupp(\nu_i)\}} Res_{p_j}(\nu_i, f_{b_j+a_j-l_j}^{(b_j)}(\tilde{\bold c}^{(j)})s^{-(a_j-l_j)}\partial_{w_j})
 = 0.
\]
Since we have $\sharp(psupp(\nu_i)\cap \{(1, l_1), \dots, (e, l_e)\})\geq 2$,
 we can rescale each $\tilde{\bold c}_j$ 
 so that the condition $(\star_{\nu_i})$ holds.
Moreover, by the way we defined $\nu_i$, we can perform this rescaling simultaneously for all $\nu_i$, $i = 1, \dots, c$.
Then, we can apply the proof of Theorem \ref{thm:main}
 and obtain a deformation of $\varphi$.\qed\\

The condition ${\rm (G)}$ in Theorem \ref{thm:gen2} 
 can be restated as follows.
Namely, 
 the polynomials $f_{b_j+1}^{(b_j)}, \dots, f_{b_j+a_j-1}^{(b_j)}$
 satisfy the condition ${\rm (G)}$ if and only if for each $i = 1, \dots, a_j-1$, 
 the element $F_{-i}\cdot \overline{Jac}_j$ is not 
 contained in $rad(F_{-1}, \dots, \check{F}_{-i}, \dots, F_{-(a_j-1)})$, 
 where $\check{F}_{-i}$ means we remove $F_{-i}$, $rad(\cdots)$ means the radical of the ideal
 generated by $(\cdots)$, 
  and $\overline{Jac}_j$ is the Jacobian of the map
 \[
 (\bar f_{b_j+1}^{(b_j)}, \dots, \bar f_{b_j+a_j-1}^{(b_j)})\colon \Bbb C^{a_j-1}\to \Bbb C^{a_j-1}.
 \]
Recall that $F_{-i}$ is of the form 
$
F_{-i}^{(j)} = \sum_{k = -i}^{-1} \Theta_{k+i}^{(j; k)} f_{b_j-k}^{(b_j)}
$
 and $\Theta_0^{(j;k)} = 1$, $\Theta_1^{(j;k)} = 0$.
Therefore, the ideal $(F_{-1}, \dots, \check{F}_{-i}, \dots, F_{-(a_j-1)})$
 can be written in the form 
\[
(f_{b_j+1}^{(b_j)}, \dots, f_{b_j+i-1}^{(b_j)}, f_{b_j+i+1}^{(b_j)}, 
 f_{b_j+i+2}^{(b_j)} + \Theta_{2}^{(j; -i)} f_{b_j+i}^{(b_j)}, \dots, 
 f_{b_j+a_j-1}^{(b_j)} + \Theta_{a_j-1-i}^{(j; -i)} f_{b_j+i}^{(b_j)}),
\]
 and we can replace $F_{-i}\cdot \overline{Jac}_j$
 by $f_{b_j+i}^{(b_j)}\cdot \overline{Jac}_j$.

Recall that the functions $f_{b+i}^{(b)}$ (and so the functions $\bar f_{b+i}^{(b)}$ and $F_{-i}$, too)
 are determined solely by the pair of positive integers $b>a\geq 2$.
The most optimistic (but not overly optimistic) expectation
 is that the condition ${\rm (G)}$ 
 holds for any pair $b>a\geq 2$ of positive integers.
Computation by Macaulay2 \cite{GS} suggests that the condition (G) holds in most cases,
 see Table \ref{table:1}.
 
\begin{table}[h]
\caption{}\label{table:1}
\begin{tabular}{|r|l|}
\hline
$a$=3 & (G) holds for $4\leq b\leq 30$             \\ \hline
4   & (G) holds for $5\leq b\leq 30$ except $b=6$ \\ \hline
5   & (G) holds for $6\leq b\leq 30$             \\ \hline
6   & (G) holds for $7\leq b\leq 30$             \\ \hline
7   & (G) holds for $8\leq b\leq 20$             \\ \hline
8   & (G) holds for $9\leq b\leq 20$             \\ \hline
9   & (G) holds for $10\leq b\leq 20$            \\ \hline
10  & (G) holds for $11\leq b\leq 15$            \\ \hline
\end{tabular}
\end{table}
 
The only exceptional case occurs when $a = 4$ and $b=6$.
This is attributed to an accidental factorization of a function due to the smallness of the degree.
Namely, when we have $a = 4$ and $b = 6$, the relevant functions are given by the following:
\[
\begin{array}{l}
F_{-1} = -\frac{3}{16}c_2^2c_3 + \frac{3}{4}c_3c_4,\\
F_{-2} = \frac{3}{128}c_2^4 -\frac{3}{16}c_2c_3^2 -\frac{3}{16}c_2^2c_4 + \frac{3}{8}c_4^2,\\
F_{-3} = \frac{3}{64}c_2^3c_3 - \frac{1}{16}c_3^3 - \frac{3}{16}c_2c_3c_4,
\end{array}
\] 
\[
\overline{Jac} = \frac{27}{16384}c_2^6c_3 + \frac{27}{2048}c_2^3c_3^3 + \frac{27}{1024}c_3^5
   -\frac{81}{4096}c_2^4c_3c_4 -\frac{27}{512}c_2c_3^3c_4 +\frac{81}{1024}c_2^2c_3c_4^2 -\frac{27}{256}c_3c_4^3.
\]

The condition $F_{-1} = F_{-3} = 0$ implies
\begin{enumerate}
\item $c_3 = 0$, or
\item $-c_2^2 + 4c_4 = 0$ and $3c_2^3 - 4c_3^2 - 12c_2c_4 = 0$. 
\end{enumerate}
In the first case, $\overline{Jac}$ is also zero, and (G) does not hold.
Moreover, in the second case, the equations implies $c_3 = 0$, too.
Thus, in this case the condition (G) fails.

Note that even in this case, a part of the condition (G) holds.
Namely, for $i=1$ and $3$, the condition of Definition \ref{def:(G)} holds.
This implies that the conclusion of Theorem \ref{thm:gen2} applies to 
 these cases, too.
 
Under the condition ${\rm (G)}$, the problem of whether deformations of $\varphi$ exists is entirely reduced to 
 verifying the cohomological condition ${\rm (D)}$, which is much easier than the obstruction calculation.

\subsection{Deformation of double points}
If $p\in C$ is a double point of the semiregular map $\varphi\colon C\to X$,
 that is, when the multiplicity $a = 2$, the function $f_{b+1}^{(b)} = \bar f_{b+1}^{(b)} = F_{-1}$ takes the simple form 
 $c_2^{\frac{b+1}{2}}$.
In this case, Theorem \ref{thm:gen2} can be significantly strengthened.
Let $\{p_1, \dots, p_l\}$  be the set of points on $C$ at which $\varphi$
 has singularities, and assume each of them is a double point.
Then, the following holds.
\begin{thm}\label{thm:doublepoint}
	The semiregular map $\varphi$ deforms if and only if at least one of the following conditions
	holds.
	\begin{enumerate}
		\item There is at least one $p_i$ such that there is no section of 
		$H^0(C, \varphi^*\omega_X(Z))$ whose polar support is $\{p_i\}$.
		Here, $Z = p_1+\cdots + p_l$ is the ramification divisor of $\varphi$.
		\item The set $H^0(C, \bar{\mathcal N}_{\varphi})$ is not zero.
	\end{enumerate}
\end{thm}
\proof
Assume the condition of (1) is satisfied.
First, assume that for any $p_i$, there is a section of $H^0(C, \varphi^*\omega_X(Z))$
 which contains $p_i$ in its polar support.
By changing the numbering if necessary, we can assume that $p_1$ is a point satisfying the condition of (1).
Also, we assume that among the points $\{p_1, \dots, p_l\}$, $p_m, \dots, p_l$ are the points for which 
 there is a section of $H^0(C, \varphi^*\omega_X(Z))$ whose polar support is $\{p_i\}$, $i = m, \dots, l$.
Here, $1<m\leq l+1$ and when $m = l+1$, the set $\{p_m, \dots, p_l\}$ is empty.
Then, we can take a basis $\{\eta_1, \dots, \eta_k\}$ of $H^0(C, \varphi^*\omega_X(Z))$
 so that $p_1$ is contained in the polar support of $\eta_k$ only when $k = 1$.
We can assume that 
\[
psupp(\eta_{k-l+i}) = \{p_{i}\}, 
\]
 for $i = m, \dots, l$.
Also, following the construction in Section \ref{subsec:obbasis}, 
 we can assume that $psupp(\eta_i)$, $1\leq i\leq k-l+m-1$ is disjoint from $\{p_m, \dots, p_l\}$, 
 and that for each $1\leq i\leq k-l+m-1$, 
\[
psupp(\eta_i)\nsubseteq \bigcup_{j<i}psupp(\eta_j)
\]
 holds.

For each $p_j$, we have a pair of integers $(a_j, b_j) = (2, b_j)$.
Let $M$ be the least common multiple of the integers
\[
\frac{b_1+1}{2}, \dots, \frac{b_l+1}{2}.
\]
Define the integer $d_j$ by $d_j = \frac{2M}{b_j+1}$.
For each $p_j$, we have an element $c_2^{(j)}$, and we take it in 
 $t^{d_j}\Bbb C[[t]]$.
Then, as in Section \ref{subsec:lowerdeformation}, the map $\varphi$ can 
 be deformed up to the order $M-1$ with respect to $t$.

Now, consider the deformation of $\varphi$ of the order $M$.
It is easy to see that
 we can take $c_2^{(j)} = t^{d_j}\tilde c_2^{(j)}$, $j = 1, \dots, m-1$,
 so that they satisfy the conditions $(\star_{\eta_i})$, $1\leq i\leq k-l+m-1$, in the sense of Definition \ref{def:star}.
Here, $\tilde c_2^{(j)}$ is a nonzero complex number.
Also, we take $c_2^{(j)} =0$, $j = m, \dots, l$.
Thus, the conditions $(\star_{\eta_i})$, $k-l+m\leq i\leq k$, also hold obviously.
Then, since the equation $(\star_{\eta_i})$ is itself the condition for the vanishing of the obstruction 
 at the order $t^M$, the map $\varphi$ deforms up to the order $t^M$ for these values of $c_2^{(j)}$.
Let us write this map as $\varphi_M$.
 
Let us consider the deformations of the higher order.
Following the argument in Section \ref{subsec:calculation of o_{N+1}}, 
 by the above choice of $c_2^{(j)}$, we see that extending the map $\varphi_M$ does not have
 an obstruction up to the order $t^{2M-1}$, and let $\varphi_{2M-1}$ be the map obtained in this manner.
At the order $t^{2M}$, there may be an obstruction $o_{2M}$ to deforming $\varphi_{2M-1}$.
When $o_{2M}$ pairs non-trivially with $\eta_i$, $1\leq i\leq k-l+m-1$, 
 we can follow the argument in Section \ref{sec:4} to modify $c_2^{(j)}$, $j = 1, \dots, m-1$, to cancel the 
 obstruction.
If $o_{2M}$ pairs non-trivially with some $\eta_{k-l+i}$, $i = m, \dots, l$,
 we can choose $\tilde c_2^{(i)}$ so that by taking $c_2^{(i)} = t^{2d_i}\tilde c_2^{(i)}$, 
 we cancel the obstruction.

With these values of $c_2^{(j)}$, we can apply the argument of Section \ref{subsec:step3}
 to construct a map $\varphi_{2M}$ which deforms $\varphi_{M'}$ for some $M'\leq 2M-1$.
Repeating this, 
 again by the argument of Section \ref{subsec:step3}, we eventually obtain a formal deformation of 
 $\varphi$, and by applying a suitable algebraization result, we finish the proof.
 
Let us assume that there is a point $p_i$ which is not contained in the polar support of 
 any element of $H^0(C, \varphi^*\omega_X(Z))$.
Let $p_1$ be such a point.
Then, by taking any non-zero $c_2^{(1)} = t^{d_1}\tilde c_2^{(1)}$, and taking all the other 
 $c_2^{(j)}$ to be zero, we obtain a non-trivial deformation $\varphi_M$ of $\varphi$.
Then, we can apply the above argument to $\varphi_M$ to obtain a formal deformation of $\varphi$.

Finally, assume the condition (2) is satisfied.
Then, a non-zero section of $H^0(C, \bar{\mathcal N}_{\varphi})$ gives a first order deformation 
 of the map $\varphi$.
By base change, we can assume that this deformation is defined at the order $t^M$.
Then, we can again apply the above argument to cancel the potential obstructions
 at the higher order by modifying the values of $c_2^{(j)}$, and we obtain a deformation of $\varphi$.
 
Conversely, suppose that neither of the conditions (1) nor (2) holds.
This means that for each $p_j$, there is a section $\eta_j$ 
 of $H^0(C, \varphi^*\omega_X(Z))$ whose polar support is $p_j$, 
 and $H^0(C, \bar{\mathcal N}_{\varphi}) = 0$.
By the latter condition, to have a non-trivial deformation, we need to take some $c_2^{(j)}$
 to be non-zero at some order with respect to $t$.
However, it couples with $\eta_j$ non-trivially and produces a non-trivial obstruction.
Thus, we cannot extend the deformation further.\qed\\

\noindent
{\bf Acknowledgments.}
The author was supported by JSPS KAKENHI Grant Number 18K03313.

\end{document}